\documentclass[12pt]{article}
\usepackage{amsmath,amsfonts,amssymb}
\usepackage[overload]{empheq}
\usepackage[margin=1in]{geometry}
\usepackage{color}
\usepackage{authblk}

\newcommand{\ve}{{\varepsilon}}
\newcommand{\ldg}{\mathrm{LdG}}
\newcommand{\er}{\mathrm{Er}}
\newcommand{\re}{\mathrm{Re}}
\newcommand{\tr}{\mathrm{tr\,}}
\renewcommand{\div}{\mathrm{div}\,}

\newcommand{\bu}{{\bf u}}
\newcommand{\bv}{{\bf v}}
\newcommand{\bvs}{{\bv_*}}
\newcommand{\tbv}{\tilde{{\bf v}}}

\newcommand{\f}{{\bf f}}
\newcommand{\bh}{h}

\newcommand{\bx}{x}
\newcommand{\by}{y}
\newcommand{\bz}{z}
\newcommand{\br}{x}
\newcommand{\hx}{\hat{x}}
\newcommand{\hy}{\hat{y}}
\newcommand{\q}{\mathsf{q}}
\newcommand{\n}{{n_*}}
\newcommand{\Q}{\mathsf{Q}}
\newcommand{\C}{\mathsf{C}}
\newcommand{\A}{\mathsf{A}}
\newcommand{\bbA}{{\mathbb A}}
\newcommand{\W}{\mathsf{W}}
\newcommand{\B}{\mathcal{B}}
\newcommand{\bB}{{\bf B}}
\newcommand{\bbB}{{\mathbb B}}
\newcommand{\D}{\mathsf{D}}
\newcommand{\E}{\mathsf{E}}
\newcommand{\F}{\mathsf{F}}
\newcommand{\bF}{{\bf F}}
\newcommand{\G}{{\mathsf{G}_\gamma}}
\newcommand{\h}{{\mathsf{h}_\gamma}}

\newcommand{\M}{\mathcal{M}}
\newcommand{\R}{\mathbb{R}}
\newcommand{\U}{{\bf U}}
\newcommand{\I}{\mathsf{I}}
\newcommand{\id}{\mathsf{I}}
\newcommand{\T}{\mathsf{T}}
\newcommand{\HH}{\mathsf{H}}
\newcommand{\sR}{\mathsf{R}}
\newcommand{\sS}{\mathsf{S}}
\newcommand{\Qc}{\mathring{\mathsf{Q}}}

\newcommand{\bxi}{\xi}
\newcommand{\ci}{\mathrm{i}}
\newcommand{\wt}{\widetilde}

\newcommand{\cA}{\mathcal{A}}
\newcommand{\cB}{\mathcal{B}}
\newcommand{\cC}{\mathcal{C}}
\newcommand{\cD}{\mathcal{D}}
\newcommand{\cM}{\mathcal{M}}
\newcommand{\cS}{\mathcal{S}}
\newcommand{\cT}{\mathcal{T}}
\newcommand{\cF}{\mathcal{F}}
\newcommand{\cI}{\mathcal{I}}
\newcommand{\cJ}{\mathcal{J}}
\newcommand{\lbar}{\overline}
\newcommand{\ft}[1]{\widetilde{#1}}

\newtheorem{remark}{Remark}

\title{Far Field Asymptotics of Nematic Flows Around a Small Spherical Particle}

\author[1]{Dmitry Golovaty}
\author[2]{Nung Kwan Yip}
\affil[1]{Department of Mathematics, The University of Akron, Akron, OH 44124, USA}
\affil[2]{Department of Mathematics, Purdue University, West Lafayette, IN 47907, USA}
\date{}

\begin{document}

\maketitle
\begin{abstract}
    Given a small spherical particle, we consider flow of a nematic liquid crystal in the corresponding exterior domain. Our focus is on precise far field asymptotic behavior of the flow in a parameter regime when the governing equations can be reduced to a system of linear partial differential equations. We are able to analytically characterize the velocity of the flow and compare it to the classical expression for the Stokes flow. The expression for velocity away from the particle can be computed either numerically or symbolically.    
\end{abstract}
\section{Introduction}
This paper analyzes the flow pattern around a small spherical particle in a nematic liquid crystalline environment. Although a similar phenomenon in a classical isotropic fluid has been well-studied, the understanding in the case of a complex fluid is far from complete, despite the fact that it is of utmost importance in many physical and biological systems. Here we are particularly interested in flows around solid particles immersed in a nematic liquid crystalline medium. A typical nematic liquid crystal consists of molecules which possess a degree of orientational order but no positional order so that this medium can be thought of as an anisotropic Newtonian fluid with additional elastic properties. Hence, in order to understand the flow field around a particle, one needs to analyze the behavior of solutions to a coupled system of PDEs that describes both flow and orientational elasticity.

In a broader context, we are interested in multi-particle systems such as colloidal suspensions of both passive and active particles in complex media. The overall behavior of these systems certainly depends on both the particle-fluid and particle-particle interactions. As a typical initial step, one would consider a dilute system in which particles are far from each other. Therefore, it is crucial to understand far field interactions between the particles.

As a starting point, one then can investigate a single particle in a host medium and determine the far-field asymptotics of the corresponding flow patterns. In the present case, these asymptotics should also include the far-field orientational information about the liquid crystal. This is the situation considered in our work for the parameter regime in which we are able to obtain explicit asymptotics in the same spirit as those known for the classical Stokes flow. Our work seems to be the first to give a precise characterization of the flow around a particle immersed in a nematic liquid crystal.

First, we briefly discuss prior work on liquid crystalline flows with or without particles. The theory of nematic flows has proceeded along two somewhat different but related directions that rely on distinct continuum descriptions of the nematics. The first and most widely used approach, due to Ericksen and Leslie \cite{Leslie_2,Leslie_1}, is based on the $\mathbb{S}^2$-valued nematic {\it director} used to model local orientation of nematic molecules, along with the velocity of the flow. The corresponding system of equations is derived by assuming local balances of mass as well as linear and angular momentum and then specifying the energy dissipation density in the form that conforms with the second law of thermodynamics. As the result, the system consists of a coupled Navier-Stokes and Ginzburg-Landau type equations. The Ericksen-Leslie model has been a subject of numerous studies \cite{MR1860443,MR1329830,MR1784963,MR1740379,Walkington}. 

The alternative approach relies on the Landau-de Gennes description of orientational elasticity that is designed to take into account certain symmetry properties of orientational distribution of nematic molecules. The Landau-de Gennes variational theory replaces the nematic director with a so-called $\Q$-tensor---a symmetric traceless $3\times3$ matrix the eigenvectors of which encode orientational properties of the nematic. The energetics of $\Q$ is modeled by the Landau-de Gennes functional $\mathcal{F}_\ldg(\Q)$ -- see \eqref{eq:ldg}. When applied to nematic flow, the evolution of $\Q$ and the velocity $\bv$ of the nematic is given by a  coupled system of the Landau-de Gennes and Navier-Stokes equations and can be derived following the procedure outlined by Sonnet and Virga (SV) on the basis of the principle of minimum constrained dissipation \cite{Sonnet_dissipative,Sonnet_tensor,SonVir}. Another related but more popular model, due to Beris and Edward (BE) \cite{beris1994thermodynamics}, was obtained by making use of the concept of a Poisson bracket with dissipation--in a later work, we will in fact show that the BE model is a specific example of the SV formulation. The BE model has been studied extensively in recent years, see for example, \cite{abels2016strong, wu2019dynamics, du2020suitable} just to name a few.

Particulate flows in a classical fluid is another widely studied subject, see \cite{galdi2002motion, galdi2011introduction} for overview of relevant literature and mathematical theory. Motion of particles in a complex fluid, and in particular in a nematic host, has also attracted attention of many investigators, especially with the expanding interest in motion of immersed active particles \cite{doi1,doi2,zhou2014living,Zhou_2017}. The experimental and modeling work in this field has so far outpaced the analysis effort \cite{ruhwandl1996friction,stark2001stokes,stark2001physics}. The goal of the present paper is to make inroads into rigorous understanding of a nematic flow around a single particle within the framework of a $\Q$-tensor-based model.

The details of the model are outlined in Section~\ref{sec:model}. Here we highlight some key physically justifiable assumptions that we make in order to obtain precise asymptotics of the flow pattern. These assumptions allow us to reduce the overall dynamics to a linear system of PDEs describing an anisotropic Stokes flow with elastic contribution. We emphasize that the distinction with the standard isotropic Stokes flow is not only due to orientability of the medium (described by $\Q$) that makes the flow anisotropic, but also because the flow is coupled to the elastic properties of the medium (described by $\nabla \Q$).

Our principal assumptions are as follows. We set both Ericksen (Er) and Reynolds (Re) numbers to be small. This means that
the dynamics is dominated by the elastic effects of the nematics and the 
flow is highly dissipative in the sense that inertial effects are negligible.
As the result, our model becomes a partially decoupled system \eqref{eq:qten2s}--\eqref{eq:inc1s} in which $\Q$ satisfies an equation that does not involve $\bv$. The solution of this equation can then be treated as a prescribed function which enters into an inhomogeneous linear Stokes-like equation describing the evolution of $\bv$. 

To obtain an explicit solution for $\Q$, we further take advantage of the small particle limit considered in \cite{QABL}. In this regime, the equation for $\Q$ becomes a Laplace equation with an exact solution given by a harmonic function. 

Our main result is the precise far-field asymptotics for $\bv$ where we are able to analytically characterize the deviation of $\bv$ from the solution for the classical Stokes flow. In particular, we are able to obtain an analytical expression for $\bv$ that can be computed either numerically or symbolically.

The outline of our paper is as follows. In Section 2, we provide a detailed description of our model and relevant parameter regimes. In Section 3, we recall the explicit solutions for a $\Q$-tensor and the classical Stokes flow in the exterior of a small particle that we use in subsequent sections. In Section 4, we analyze the structure of our anisotropic Stokes system. Then in Sections 5 and 6, we prove existence of solutions to the governing system of equations and analyze their far-field asymptotic behavior, in particular their deviation from the isotropic Stokes flow. Finally, in Section 7, we present numerical results. Given the generality of the model considered in this work, the analysis relies on some tedious but routine computations that we summarize in the Appendices A-D, along with verification of our numerical experiments.

We conclude this section with the list of notation and conventions that will be used in this paper. We remind the reader that we work with objects defined on $\R^3$.
\begin{enumerate}
\item 
For any vector $u \in \mathbb{R}^3$, its components are represented by
$u_i$, for $i = 1, 2, 3$ and
$\displaystyle \widehat{u} = \frac{u}{|u|}$ is
the unit vector with the same direction as $u$. 
Hence we have $\displaystyle \widehat{u}_i = \frac{u_i}{|u|}$.
Whenever there is no ambiguity, we will omit
 $\,\,\widehat{}\,\,$. Euclidean inner product between vectors $u,v\in\R^3$
is denoted by $u\cdot v$.
Furthermore, the symbols $e_1,e_2$ and $e_3$ represent
the coordinate vectors of $\R^3$.

\item For any $\mathsf{B},\mathsf{C}\in M^{3\times3}$, $3\times 3$ matrices,
$\mathsf{B}\cdot \mathsf{C}=\tr\left(\mathsf{C}^T\mathsf{B}\right)$ and
$|\mathsf{B}|^2=\mathsf{B}\cdot\mathsf{B}$. 

\item Einstein's convention for repeated indices will be used whenever possible.

\item For any matrix or second order tensor field 
$\T=(\T_{ij})_{1\leq i,j \leq 3}$, $\T_i$ denotes its $i$-th row.
The divergence $\div\T$ 
is to be taken row-wise, i.e. for all $i$, 
$(\div\T)_i := \div\T_i = \partial_{x_j}\T_{ij}$.

\item A symbol such as $x=(x_1, x_2, x_3)^T$ denotes either a generic point or a dummy integration variable in $\mathbb{R}^3$.
For convenience, we will also use $r$ to denote
\[
r=|x|=\sqrt{x_1^2+x_2^2+x_3^2}.
\]

\item The (3-dimensional) volume integration element is denoted by
\[
dx = dx_1dx_2dx_3.
\]
The (2-dimensional) area integration element is denoted by
\[
d\sigma_x,\,\,\,\text{or simply}\,\,\,d\sigma\,\,\,
\text{(if the variable is clear from context).}
\]
The unit outward normal to the domain $\Omega$ occupied by the liquid crystal
is given by $\nu$.

\item	We use the $\partial_i = \partial_{x_i}$
to denote spatial partial derivatives.

\item Sometimes we find it advantageous to use the symbol
$\big\langle\text{matrix}, \text{vector}\big\rangle$
and 
$\big\langle\text{vector}, \text{vector}\big\rangle$
to denote the matrix-vector multiplication and the dot product, respectively.

\item\label{notation} The symbol $``:"$ will be used to denote a generic tensor contraction. It will be explicitly defined whenever necessary.

\item The symbol $``\lesssim"$ will be used to refer to an inequality that holds up to a multiplicative constant. This constant can change from one line to another but is not important to the ultimate conclusion.
\end{enumerate}

\section{Problem formulation}\label{sec:model}
\subsection{A model for the dynamics of nematic liquid crystals}
Suppose that a nematic liquid crystal occupies a region $\Omega\subset\mathbb{R}^3$ and let the $\Q$-tensor 
\[
\Q:\Omega\to\mathcal{M}:=\left\{\Q\in M^{3\times3}(\mathbb{R}):\Q^T=\Q,\ \tr{\Q}=0\right\}\] 
and the velocity $\bv:\Omega\to\mathbb{R}^3$ describe the local state of the nematic. 
We introduce the Landau-de Gennes energy density and the corresponding energy functional as:
\begin{equation}
\label{eq:ldg}
\mathcal{E}_\ldg(\nabla\Q, \Q)=\frac{K}{2}{|\nabla\Q|}^2+f(\Q),\quad
\mathcal{F}_\ldg(\Q)=\int\mathcal{E}_\ldg(\nabla\Q,\Q)\,d^3x.
\end{equation}
The parameter $K$ is the elastic constant of the liquid crystal.
The nonlinear Landau-de Gennes potential $f$ is given by
\begin{equation}
\label{eq:pote}
f(\Q)=-\frac{A}{2}\tr\Q^2+\frac{B}{3}\tr\Q^3+\frac{C}{4}{\left(\tr\Q^2\right)}^2,
\end{equation}
with $A,C > 0$. The minimum of $f$ over $\mathcal{M}$ is attained at the 
nematic states given by
\begin{equation}\label{nematics}
s_*\left(\n\otimes \n - \frac13\I\right),\,\,\,\text{where}\,\,\,
 \n\in{\mathbb S}^2\,\,\,\text{and}\,\,\,
s_* = \left\{
\begin{array}{ll}
-\frac{B+\sqrt{B^2+24AC}}{4C} &\text{for $B > 0$},\\
\frac{-B+\sqrt{B^2+24AC}}{4C} &\text{for $B < 0$}.
\end{array}
\right.
\end{equation}

Let
\begin{equation}
\label{eq:circ}
\Qc=\dot\Q+\Q\W-\W\Q,\quad\A=\frac{1}{2}\left(\nabla\bv+\nabla\bv^T\right),\quad\W=\frac{1}{2}\left(\nabla\bv-\nabla\bv^T\right)
\end{equation}
where $\dot\Q$ denotes the convective derivative
\[
\dot\Q = \partial_t\Q + \bv\cdot\nabla\Q.
\]
The expression $\Qc$ will be used below instead of $\dot\Q$ because it is frame indifferent
\cite[Section 2.1.3, Eq. (2.87), Section 4.1.3]{SonVir}. The $3\times3$-matrices $\A$ and $\W$ are respectively symmetric and skew-symmetric, and the expressions such as $\Q\W$ and so forth refer to matrix multiplications. From the definition of $\Qc$, we infer that it is also a symmetric, traceless matrix.

We consider the incompressible flow of the nematic in $\Omega$ 
described by the following system of equations
\begin{align}[left=\empheqlbrace]
&\frac{\partial\mathcal{E}_\ldg}{\partial\Q}-\div\left[\frac{\partial\mathcal{E}_\ldg}{\partial\nabla\Q}\right]-\Lambda\id+\zeta_1\Qc+\zeta_2\A+\frac{\zeta_3}{2}(\A\Q+\Q\A)+\zeta_9(\A\cdot\Q)\Q=0, \label{eq:qten}\\
&\rho\,\dot\bv+\div\left[p\I-\T^\mathrm{v}_{\text{SV}}
-\T^\mathrm{el}\right]=0, \label{eq:velo}\\
&\div\bv=0, \label{eq:inc} 
\end{align}
where we have 
\begin{eqnarray}
\frac{\partial\mathcal{E}_\ldg}{\partial\Q} & = & 
-A\Q + B\left(\Q^2 - \frac13|\Q|^2\I\right) + C|\Q|^2\Q,
\label{nonlinf}\\
\div\left[\frac{\partial\mathcal{E}_\ldg}{\partial\nabla\Q}\right]
& = & K\triangle\Q.
\label{LDlaplace}
\end{eqnarray}
The parameter $\rho$ denotes the density of the nematic and $\zeta_1$ through $\zeta_{11}$ are the viscosity coefficients of the nematic with $\zeta_8$ being the isotropic viscosity. The pressure $p$ and the function $\Lambda$ are the Lagrange multipliers corresponding, respectively, to the incompressibility and tracelessness constraints for $\bv$ and $\Q$. The viscous and elastic stress tensors are 
\begin{eqnarray}
\T^\mathrm{v}_{\text{SV}}
&=&\zeta_1\left(\Q\Qc-\Qc\Q\right)+\zeta_2\left(\Qc+\Q\A-\A\Q\right)+\frac{\zeta_3}{2}\left(\Q\Qc+\Qc\Q+\Q^2\A-\A\Q^2\right)\nonumber\\
&&+\zeta_4\left(\Q\A+\A\Q\right)+\zeta_5\left(\Q^2\A+\A\Q^2\right)+\zeta_6(\A\cdot\Q)\Q+\zeta_7{|\Q|}^2\A+\zeta_8\A\nonumber\\
&&+\zeta_9(\mathring{{\Q}}\cdot\Q)\Q+\zeta_{10}\left(\left(\Q^2\cdot\A\right)\Q+\left(\Q\cdot\A\right)\Q^2\right)
+\zeta_{11}|\Q|^2(\A\cdot\Q)\Q
\label{eq:visc}
\end{eqnarray}
and
\begin{eqnarray}
\T^\mathrm{el}
&=&-\frac{\partial\mathcal{E}_\ldg}{\partial\left(\partial_j\Q_{mn}\right)}\left(\partial_i\Q_{mn}\right){\mathbf e}_i\otimes{\mathbf e}_j
=-K(\partial_j\Q_{mn})(\partial_i\Q_{mn}){\mathbf e}_i\otimes{\mathbf e}_j
\nonumber\\
&=&-K\left(\partial_i\Q\cdot\partial_j\Q\right){\mathbf e}_i\otimes{\mathbf e}_j
\nonumber\\
&=:&-K\left(\nabla\Q\odot\nabla\Q\right),
\label{eq:elas}
\end{eqnarray}
respectively, so that 
$\left(\nabla\Q\odot\nabla\Q\right)_{ij} = \partial_i\Q\cdot\partial_j\Q$ component-wise.

The model \eqref{eq:qten}--\eqref{eq:inc} is established in 
\cite[Section 4.1.3]{SonVir}, where a dissipation function $R=R(\Q;\A,\Qc)$ is introduced from which the viscous stress tensor \eqref{eq:visc} can be derived. In particular, the following general form of a dissipation function is proposed in \cite[Eq. (4.23)]{SonVir}:
\begin{multline}
\label{eq:diss}
R(\Q;\A,\Qc)=\frac{\zeta_1}{2}\Qc\cdot\Qc+\zeta_2\A\cdot\Qc+\zeta_3(\Qc\Q)\cdot\A+\zeta_4\Q\cdot\A^2+\zeta_5\Q^2\cdot\A^2+\frac{\zeta_6}{2}(\Q\cdot\A)^2\\+\frac{\zeta_7}{2}{|\A|}^2{|\Q|}^2+\frac{\zeta_8}{2}{|\A|}^2+\zeta_9(\Qc\cdot\Q)(\A\cdot\Q)+\zeta_{10}\left(\Q^2\cdot\A\right)\left(\Q\cdot\A\right)
+\frac{\zeta_{11}}{2}|\Q|^2(\A\cdot\Q)^2.
\end{multline}
We remark that \cite{SonVir} considers only the terms 
involving $\zeta_1$ to $\zeta_8$ so that $R$ is exactly quadratic in 
the rates $\Qc$ and $\A$ and at most quadratic in $\Q$ but the idea can
certainly be generalized. In particular, $R$ can involve a linear combination 
of nineteen invariants of the tensor triple $(\Q,\Qc,\A)$ 
\cite[p. 223]{SonVir}.
A simpler version of \eqref{eq:diss}, introduced in \cite{QiaShe}, given by $\zeta_3=\zeta_5=\zeta_7=0$,
also appears in \cite[Eq. (4.25)]{SonVir}.
Here, we include the terms with coefficients $\zeta_9$ to $\zeta_{11}$ because with these terms the model \eqref{eq:qten}-\eqref{eq:inc} subsumes the model in \cite{beris1994thermodynamics}, derived by Beris and Edward. In a forthcoming work, we will show that the Beris-Edward model is, in fact, a particular version of \eqref{eq:qten}-\eqref{eq:inc}, corresponding to a specific choice of dissipation constants $\zeta_i$, $i=1,\ldots,11$. 

The above Sonnet-Virga model is derived using a variational 
framework together with the Principle of Minimum Constrained Dissipation
\cite[Section 2.2, Eqn. (2.172), (2.178), (2.179)]{SonVir}.
Using $R$, we can rewrite 
\eqref{eq:qten}-\eqref{eq:inc} as
-- see also \cite[Eq. (4.21) and (4.22)]{SonVir}:
\begin{align}[left=\empheqlbrace]
&\frac{\partial\mathcal{E}_\ldg}{\partial\Q}-\div\left[\frac{\partial\mathcal{E}_\ldg}{\partial\nabla\Q}\right] + \frac{\partial R}{\partial\Qc}
=0, 
\label{eq:qten0}\\
&\rho\,\dot\bv-\div\left[\T\right]=0, 
\label{eq:velo0}\\
&\div\bv=0, \label{eq:inc0} 
\end{align}
where
\begin{equation}
\T = -p\I + \T_{\text{SV}}^\text{v} + \T^\text{el}
\label{TR1}
\end{equation}
and
\begin{equation}
\T^\mathrm{v}_{\text{SV}} 
= 
\frac{\partial R}{\partial\A}
+\Q\frac{\partial R}{\partial\Qc}
-\frac{\partial R}{\partial\Qc}\Q
\,\,\,\,\,\,\text{and}\,\,\,\,\,\,
\T^\text{el}
= -\nabla\Q\odot
\frac{\partial\mathcal{E}_\ldg}{\partial\left(\nabla\Q\right)}.
\label{TR2}
\end{equation}

To complete the description of the model, we point out that the dissipation function $R$ 
has to be positive semidefinite. 
The set of conditions satisfied by $\Q$ and 
$\zeta_i,\ i=1,\ldots,11$ to ensure that this property holds is too 
complicated to list here but some of these conditions can be found in \cite[pp. 221-222]{SonVir}. However, if $R$ only has the terms associated with
$\zeta_1, \zeta_2$ and $\zeta_8$, i.e.
\[
R(\Q;\A,\Qc)=\frac{\zeta_1}{2}\Qc\cdot\Qc+\zeta_2\A\cdot\Qc
+\frac{\zeta_8}{2}{|\A|}^2,
\]
then $R$ is positive semidefinite if and only if
\begin{equation}
\label{eq:disineq}
\zeta_1 > 0\,\,\,\text{and}\,\,\,\zeta_2^2 \leq \zeta_1\zeta_8.
\end{equation}


\subsection{Non-dimensionalization}
Suppose $L$ is the characteristic length of the system and
$v$ is the characteristic velocity. With these, we introduce
\[\tilde{x}=\frac{x}{L},\quad
\tilde{\bv} = \frac{\bv}{v},\quad
\tilde{t}=\frac{vt}{L},\quad 
\tilde{f}(\Q)=\frac{f(\Q)}{C},\,\,\,\text{and}\,\,\,
\tilde{\mathcal{E}}_\ldg={\frac{L^2}{K}\mathcal{E}_\ldg},
\]
so that
\[
\partial_t = \frac{v}{L}\partial_{\tilde t},\quad
\partial_x = \frac{1}{L}\partial_{\tilde x},\quad
\Q=\tilde\Q,\quad
\A = \frac{v}{L}\tilde\A,\quad
\tilde\W = \frac{v}{L}\tilde\W,\quad
\dot{\Q} = \frac{v}{L}\dot{\widetilde{\Q}},\quad
\mathring{\Q} = \frac{v}{L}\mathring{\widetilde{\Q}},
\]
\[
\tilde{f}(\Q) =
-\frac{A}{2C}\tr\Q^2+\frac{B}{3C}\tr\Q^3+\frac{1}{4}{\left(\tr\Q^2\right)}^2,
\quad\text{and}\quad
\tilde{\mathcal{E}}_\ldg(\tilde\Q)
=\frac{1}{2}{|\tilde\nabla\tilde Q|}^2+\frac{CL^2}{K}\tilde f(\tilde\Q).
\]
The forms of $\tilde{f}$ and $\tilde{\mathcal{E}}_\ldg$ are consistent with 
the following change of variable relation for the Ginzburg-Landau functional,
\[
\mathcal{F}_\ldg(\Q)
= \int \frac{K}{2}|\nabla\Q|^2 + f(\Q)\,dx
= 
KL\int \frac{1}{2}|\tilde\nabla\tilde\Q|^2 + \frac{CL^2}{K}\tilde f(\tilde\Q)\,d\tilde x
= 
KL\int\tilde{\mathcal{E}}_\ldg(\tilde{Q})\,d\tilde x.
\]
If we further introduce
$\displaystyle \ve=\frac{1}{L}\sqrt{\frac{K}{C}}$ as the
ratio between the nematic correlation length 
$\displaystyle \sqrt{\frac{K}{C}}$ and $L$,
we can then write the energy density in the following non-dimensional form,
\[
\tilde{\mathcal{E}}_\ldg(\tilde\nabla\tilde\Q,\tilde\Q)
=\frac{1}{2}{|\tilde\nabla\tilde Q|}^2+\frac{1}{\ve^2}\tilde f(\tilde\Q).
\]
If $\ve \ll 1$, then the potential function imposes a heavy penalty on deviations of $\Q$ from 
the nematic states \eqref{nematics}.

The nondimensional system of the governing equations \eqref{eq:qten}-\eqref{eq:inc} is then
\begin{align}
&\frac{\partial\tilde{\mathcal{E}}_\ldg}{\partial\tilde\Q}
-\widetilde{\div}\left[\frac{\partial\tilde{\mathcal{E}}_\ldg}{\partial\tilde\nabla\tilde\Q}\right]-\tilde\Lambda\id
+\er
\left(\gamma_1\mathring{\tilde\Q}+\gamma_2\tilde\A+\frac{\gamma_3}{2}(\tilde\A\tilde\Q+\tilde\Q\tilde\A)+\gamma_9(\tilde\A\cdot\tilde\Q)\tilde\Q\right)=0, \label{eq:qtend}\\
&\re\,\dot{\tilde{\bv}}+\widetilde{\div}\left[
\tilde{p}\id
-\tilde{\T}^\mathrm{v}_{\text{SV}}
-\frac{1}{\er}\tilde\T^\mathrm{el}
\right]=0, \label{eq:velnd}\\
&\widetilde{\div}\tilde\bv=0, \label{eq:incnd} 
\end{align}
where the nondimensional viscous and elastic stress tensors are 
\begin{eqnarray}
\tilde{\T}^\mathrm{v}_{\text{SV}}
&=&\gamma_1\left(\tilde\Q\mathring{\tilde{\Q}}-\mathring{\tilde{\Q}}\tilde\Q\right)
+\gamma_2\left(\mathring{\tilde{\Q}}+\tilde\Q\tilde\A-\tilde\A\tilde\Q\right)
+\frac{\gamma_3}{2}\left(\tilde\Q\mathring{\tilde{\Q}}+\mathring{\tilde{\Q}}\tilde\Q+\tilde\Q^2\tilde\A-\tilde\A\tilde\Q^2\right)\nonumber\\
&&+\gamma_4\left(\tilde\Q\tilde\A+\tilde\A\tilde\Q\right)+\gamma_5\left(\tilde\Q^2\tilde\A+\tilde\A\tilde\Q^2\right)+\gamma_6(\tilde\A\cdot\tilde\Q)\Q
+\gamma_7{|\tilde\Q|}^2\tilde\A+\tilde\A\nonumber\\
&&+\gamma_9(\mathring{\tilde{\Q}}\cdot\tilde\Q)\tilde\Q+\gamma_{10}\left(\left(\tilde\Q^2\cdot\tilde\A\right)\tilde\Q+\left(\tilde\Q\cdot\tilde\A\right)\tilde\Q^2\right)
+\gamma_{11}|\tilde\Q|^2(\tilde\A\cdot\tilde\Q)\tilde\Q
\label{eq:visc5}
\end{eqnarray}
and
\begin{equation}
\label{eq:elas2}
\tilde{\T}^\mathrm{el}
= -\tilde\nabla\tilde\Q\odot\tilde\nabla\tilde\Q,
\end{equation}
respectively. In the above, the nondimensional groups
\begin{equation}
\er=\frac{\zeta_8vL}{K}\,\,\,\mbox{ and }\,\,\,\re=\frac{\rho v L}{\zeta_8}
\end{equation}
are respectively, the Ericksen and Reynolds number with
\begin{equation}
\label{eq:gag}
\gamma_i=\frac{\zeta_i}{\zeta_8},\ i=1,\ldots,7,9,\ldots,11.
\end{equation} 
We note that the constant $\er$ gives the ratio between 
the viscous and elastic forces while $\re$ gives the ratio between the 
inertial and viscous forces.

For the rest of the paper, we will refer to our nondimensionalized system
\eqref{eq:qtend}--\eqref{eq:incnd} but for simplicity omit the tildes in all variables.

\subsection{Decoupling of the governing equations}\label{sec:decouple}
In this section, we introduce physical regimes in which
the far-field behavior of $\Q$ and $\bv$ can be explicitly characterized.
In particular, we will assume that the Ericksen number is small so that
-- as we will see shortly -- the equation for the $\Q-$tensor decouples from the equation for the velocity $\bv$.
To this end, observe that 
\begin{eqnarray}
\div\T^\mathrm{el}&=&-\div\left(\left(\partial_i\Q\cdot\partial_j\Q\right){\mathbf e}_i\otimes{\mathbf e}_j\right)=-\partial_j(\partial_i\Q\cdot\partial_j\Q){\mathbf e}_i\nonumber\\
&=&-\nabla\left(\frac12|\nabla\Q|^2\right)-\left(\partial_i\Q\cdot\Delta\Q\right){\mathbf e}_i,
\label{eq:eldi}
\end{eqnarray}
from which $\div\T^\mathrm{el}$ can be eliminated. More precisely,
using \eqref{eq:qtend}, we have 
 \begin{equation}
\label{eq:wert}
\partial_i\Q\cdot\Delta\Q=\frac{1}{\ve^2}\frac{\partial f}{\partial\Q}\cdot\partial_i\Q+\er\left(\gamma_1\Qc+\gamma_2\A+\frac{\gamma_3}{2}(\A\Q+\Q\A)+\gamma_9(\A\cdot\Q)\Q\right)\cdot\partial_i\Q,
\end{equation}
where the $\Lambda\I$ disappears as $\tr{\Q}=0$:
$\Lambda\I\cdot\partial_iQ = \Lambda\partial_i\text{tr}{\Q} = 0$.
Combining \eqref{eq:eldi} and \eqref{eq:wert}, we have
\begin{eqnarray*}
&&\div\T^\mathrm{el}\\
&=&-\nabla\left(\frac12|\nabla\Q|^2+\frac{1}{\ve^2}f(\Q)\right)
-\er\left(\gamma_1\Qc+\gamma_2\A+\frac{\gamma_3}{2}(\A\Q+\Q\A)+\gamma_9(\A\cdot\Q)\Q\right)\cdot\partial_i\Q{\mathbf e}_i\\
&=&-\nabla\mathcal{E}_\ldg(\Q)
-\er\left(\gamma_1\Qc+\gamma_2\A+\frac{\gamma_3}{2}(\A\Q+\Q\A)+\gamma_9(\A\cdot\Q)\Q\right)\cdot\partial_i\Q{\mathbf e}_i.
\end{eqnarray*}
Hence
\begin{align}
\div\left(-p\I + \T^\mathrm{v}_{\text{SV}} + \frac{1}{\er}\T^\mathrm{el}
\right)
\,\,\,=\,\,\,&
\div\left(-p\I - \frac{1}{\er}\mathcal{E}_\ldg(\Q)\I 
+ \T^\mathrm{v}_{\text{SV}}\right)\nonumber\\
&-\left(\gamma_1\Qc+\gamma_2\A+\frac{\gamma_3}{2}(\A\Q+\Q\A)+\gamma_9(\A\cdot\Q)\Q\right)\cdot\partial_i\Q{\mathbf e}_i.
\end{align}
Absorbing the gradient of the Landau-de Gennes energy density 
into the pressure field leads to the following form of our system
\eqref{eq:qtend}--\eqref{eq:incnd}:
\begin{align*}[left=\empheqlbrace]
&\frac{\partial\mathcal{E}_\ldg}{\partial\Q}-\div\left[\frac{\partial\mathcal{E}_\ldg}{\partial\nabla\Q}\right]-\Lambda\id+
\er\left(\gamma_1\Qc+\gamma_2\A+\frac{\gamma_3}{2}(\A\Q+\Q\A)+\gamma_9(\A\cdot\Q)\Q\right)=0,\\
&\re\,\dot\bv+\div\left[p\id-\T^\mathrm{v}_\text{SV}\right]
+\left(\gamma_1\Qc+\gamma_2\A+\frac{\gamma_3}{2}(\A\Q+\Q\A)+\gamma_9(\A\cdot\Q)\Q\right)\cdot\partial_i\Q{\mathbf e}_i=0, \\
&\div\bv=0. 
\end{align*}

In this paper, we consider a specific regime of the above system
so that the far-field spatial behavior can be revealed explicitly.
This is described as follows.
\begin{enumerate}
\item The Ericksen number is small, i.e. $\er\to0$,
so that the elastic stress in the liquid crystal dominates the viscous stress. Formally, this leads to
\begin{align*}[left=\empheqlbrace]
&\frac{\partial\mathcal{E}_\ldg}{\partial\Q}-\div\left[\frac{\partial\mathcal{E}_\ldg}{\partial\nabla\Q}\right]=0,\\
&\re\,\dot\bv+\div\left[p\id-\T^\mathrm{v}_\text{SV}\right]
+\left(\gamma_1\Qc+\gamma_2\A+\frac{\gamma_3}{2}(\A\Q+\Q\A)+\gamma_9(\A\cdot\Q)\Q\right)\cdot\partial_i\Q{\mathbf e}_i=0,\\
&\div\bv=0
\end{align*}
where there is no need for $\Lambda\I$ as the tracelessness condition for $\Q$
is already incorporated in \eqref{nonlinf}.
Note that the first equation is the Euler-Lagrange equation for the Landau-de 
Gennes energy and it is \emph{decoupled} from the equation for the velocity.
In other words, the tensor field $\Q$ serves as the 
inhomogeneous source term for the velocity field $\bv$.

\item The characteristic length of the problem is much smaller than the 
nematic correlation length, i.e., $\varepsilon\to\infty$ so that
$\displaystyle \frac{\partial\mathcal{E}_\ldg}{\partial\Q}
= \frac{1}{\varepsilon^2}\nabla f(\Q) \approx 0$.
Hence we are led to the following system:
\begin{align}[left=\empheqlbrace]
&-\Delta\Q=0,\label{eq:qten2}\\
&\re\,\dot\bv+\div\left[p\id-\T^\mathrm{v}_\text{SV}\right]
=-\left(\gamma_1\Qc+\gamma_2\A+\frac{\gamma_3}{2}(\A\Q+\Q\A)+\gamma_9(\A\cdot\Q)\Q\right)\cdot\partial_i\Q{\mathbf e}_i, \label{eq:velo1}\\
&\div\bv=0\label{eq:inc1},
\end{align}
Note that now $\Q$ is harmonic. 
We remark that the positive semidefiniteness of the 
dissipation function is ensured if among others, the inequality 
\eqref{eq:disineq} is satisfied.
However, as $\Q$ is fixed or actually ``prescribed" by \eqref{eq:qten2}
in our asymptotic regime, these inequalities can simply be replaced by the 
condition that all the
coefficients, $\gamma_1$ through $\gamma_{11}$ are sufficiently small.

\item
The Reynolds number is small and the system has reached stationarity,
i.e. $\re\rightarrow 0$ and $\partial_t\Q = \partial_t\bv = 0$.
Hence the system \eqref{eq:qten2}--\eqref{eq:inc1} becomes
\begin{align}[left=\empheqlbrace]
&-\Delta\Q=0,\label{eq:qten2s}\\
&\div\left[p\id-\T^\mathrm{v}_\text{SV}\right]=-\left(\gamma_1\Qc+\gamma_2\A+\frac{\gamma_3}{2}(\A\Q+\Q\A)+\gamma_9(\A\cdot\Q)\Q\right)\cdot\partial_i\Q{\mathbf e}_i, \label{eq:velo1s}\\
&\div\bv=0\label{eq:inc1s},
\end{align}
Note that in this case, we have $\Qc=\bv\cdot\nabla\Q + \Q\W-\W\Q$.
We record 
${\T}^\mathrm{v}_{\text{SV}}$ \eqref{eq:visc5} here again for convenience,
\begin{eqnarray}
{\T}^\mathrm{v}_{\text{SV}}
&=&\gamma_1\left(\Q\mathring{{\Q}}-\mathring{{\Q}}\Q\right)
+\gamma_2\left(\mathring{{\Q}}+\Q\A-\A\Q\right)
+\frac{\gamma_3}{2}\left(\Q\mathring{{\Q}}+\mathring{{\Q}}\Q+\Q^2\A-\A\Q^2\right)\nonumber\\
&&+\gamma_4\left(\Q\A+\A\Q\right)+\gamma_5\left(\Q^2\A+\A\Q^2\right)+\gamma_6(\A\cdot\Q)\Q
+\gamma_7{|\Q|}^2\A+\A\nonumber\\
&&+\gamma_9(\mathring{{\Q}}\cdot\Q)\Q+\gamma_{10}\left(\left(\Q^2\cdot\A\right)\Q+\left(\Q\cdot\A\right)\Q^2\right)
+\gamma_{11}|\Q|^2(\A\cdot\Q)\Q.
\label{eq:visc5-2}
\end{eqnarray}
\end{enumerate}

To conclude, the current paper analyzes the stationary system 
\eqref{eq:qten2s}-\eqref{eq:inc1s} in the domain exterior to a sphere of 
radius $a>0$. 
We are particularly interested in the far-field spatial behavior of the flow.

To complete
the description of the above system, we need to incorporate boundary
conditions for $\bv$ and $\Q$ which are discussed next.

\subsection{Boundary conditions}\label{sec:BCs}

We will solve the above system in the exterior domain 
$\Omega=\mathbb{R}^3\backslash\bB_a(0)$ in a {\em moving frame}.
The following boundary conditions will be imposed for $\Q$ and $\bv$:
\begin{equation}
\label{eq:ext1}
\Q \longrightarrow \Q_*,\quad
\bv \longrightarrow \bv_*
\quad\text{as $|x|\longrightarrow\infty$.}
\end{equation}
\begin{equation}
\label{eq:ext2}
\text{on $\partial\Omega$}:\quad
\left\{
\begin{array}{rl}
\frac{1}{w}\frac{\partial\Q}{\partial\nu} &= \Q_b-\Q,\\
\bv &= \bv_b.
\end{array}
\right.
\end{equation}
In the above, $\Q_*, \bv_*$ are the far-field states for $\Q$ and $\bv$, 
$\Q_b$ and $\bv_*$ are to be specified, $w$ is some positive number,
and $\nu$ is the outward unit normal for $\partial\Omega$ 
(or inward to $\mathbf{B}_a(0)$). We make the
following remarks about the above boundary conditions.

\begin{remark} 
The boundary condition for $\Q$ is associated with the following
surface anchoring energy:
\[
\mathcal{F}_s(\Q) = \frac{w}{2}\int_{\partial\Omega}|\Q-\Q_b|^2\,d\sigma.
\]
We will choose $\Q_b$ so as to have an explicit solution for $\Q$ --
see \bf{Section \ref{Q.stat}}.
\end{remark}

\begin{remark}\label{vBC.rmk}
The problem \eqref{eq:qten2}-\eqref{eq:inc1}, \eqref{eq:ext1}-\eqref{eq:ext2} describes the flow of a nematic liquid crystal in the exterior of a colloidal particle under various scenarios. We emphasize that the quantity
$\bv_b$ is defined in the frame associated with the moving particle.

\begin{enumerate}
\item
For a passive particle, the condition $\bv_b = 0$ describes no-slip boundary 
conditions on the surface of a particle that is stationary with respect to an 
inertial frame. The second condition in \eqref{eq:ext1} imposes the constant 
velocity $\bv_*$ of the flow at infinity. 

\item
Now suppose the passive particle 
moves in the nematic fluid with an externally imposed velocity $\bv_c$ subject to the no-slip boundary condition, while the nematic is stationary at infinity, then the velocity on the boundary of the particle and at infinity equal 
$\bv_b=\bv_c$ and $\bv_*=0$, respectively. 
In this case, if we go to a frame moving with the particle, then the velocity 
of the nematic liquid crystal at infinity will equal $\bv_*=-\bv_c$ and the 
velocity on the surface of the particle will vanish, that is, $\bv_b=0$. 
However, this change of frame will induce an additional forcing term in 
\eqref{eq:velo1s} due to the presence of the 
convective derivative $\bv\cdot\nabla\Q$ of $\Q$. This term will be 
described more explicitly in Section \ref{DecompSys}, in particular,
equation \eqref{eq:A}.

\item
For a general active particle, $\bv_b$ is typically prescribed and nonconstant on the surface of the particle. Similar to the previous paragraph, 
if the particle moves with constant velocity $\bv_c$, then changing to a frame 
moving with the particle, we can replace $\bv_b$ and $\bv_*$ by 
$\bv_b-\bv_c$ and $\bv_*=-\bv_c$. The value of $\bv_c$ can be determined by 
solving the problem \eqref{eq:qten2}-\eqref{eq:inc1}, \eqref{eq:ext1}-\eqref{eq:ext2} and choosing $\bv_c$ so that the total stress on the surface of the 
particle vanishes.
\end{enumerate}
\end{remark}

From the mathematical point of view, in dimensons three or higher, 
the problem \eqref{eq:qten2}-\eqref{eq:inc1}, \eqref{eq:ext1}-\eqref{eq:ext2} is uniquely solvable in exterior domains
for any $\bv_b$ and $\bv_*$.
This is in contrast to the situation in bounded domain for which, due to the 
incompressibility condition, the total flux 
$\int\bv_b\cdot\nu\,d\sigma$ at the boundary must vanish. In dimension two,
in general, there is no solution in the exterior domain with a bounded velocity field. This is the origin of Stokes' paradox.

As a final remark, we point out that it is certainly advantageous to 
consider our problem in the frame of the (moving) particle so that the 
domain does not change in time.
Thus, in the case of passive particle which is the emphasis of the current paper, we set 
$\bv_b=0$ and $\bv_*$ to be some prescribed value. 
(As mentioned above, the extension to active particle in the current framework 
is achieved by setting $\bv_b$ to be some general non-constant function.
See \cite{Lighthill} for examples of such a function.)
Our goal is then to compute and analyze the flow pattern of the nematic fluid. 

\section{Preliminary information}
We will make use of a known explicit stationary solution for $\Q$ in 
the exterior domain $\Omega$ and investigate its role 
in determining the flow pattern. 

\subsection{Stationary state for $\Q$.}\label{Q.stat}
Under the physical regime and boundary conditions considered in the
Sections \ref{sec:decouple} and \ref{sec:BCs}, we are looking for a 
$\Q$-tensor function satisfying
\begin{eqnarray*}
\triangle \Q & = & 0,
\quad\text{in $\Omega$}\\
\frac{1}{w}\frac{\partial \Q}{\partial\nu} & = & \Q_b - \Q
\quad\text{on $\partial\Omega$,}\\
\Q & = & \Q_*\quad\text{at $|x|=\infty$.}
\end{eqnarray*}
The work Alama-Bronsard-Lamy \cite[Theorem 1]{QABL} gives the following 
explicit solution:
\begin{eqnarray}
\Q &=&
\left(1 - \frac{w}{1+w}\frac{1}{r}\right)\Q_* +\frac{w}{3+w}\frac{1}{r^3}\Q_b,
\quad r>1,
\label{QABL}
\\
\text{where}\quad
\Q_* & = & s_*\left(\n\otimes \n - \frac{\id}{3}\right),
\quad\text{with a given $\n\in\mathbb{S}^2$}
\label{QABL2}\\
\text{and}\quad
\Q_b & = &
s_*\left(\widehat{x}\otimes \widehat{x} - \frac{\id}{3}\right),
\quad \widehat{x} = \frac{x}{|x|},
\label{QABL3}
\end{eqnarray}
with parameters $s_*,w > 0$ and the particle radius $a$ is taken to be one.
Note that $\Q$ is harmonic 
and the boundary function $\Q_b$ is of ``hedgehog'' type.
The function $\Q$ has the following far-field spatial asymptotics
\begin{equation}\label{Q.farfield}
\Q \sim \Q_* + O\left(\frac{1}{r}\right)
\,\,\,\text{and}\,\,\,
\nabla \Q \sim O\left(\frac{1}{r^2}\right),
\,\,\,\text{for $r\gg 1$}.
\end{equation}

We remark that \cite{QABL} derives the above equation in the small
particle regime $a^2 \ll K$, corresponding to $\varepsilon \gg 1$. For large particle, 
$a^2 \gg K$, corresponding to $\varepsilon \ll 1$, the solution $\Q$ tends to a 
harmonic map into ${\mathbb S}^2$ taking the form, 
$\Q(x) = s_*\left(n(x)\otimes n(x) - \frac{\I}{3}\right)$
\cite[Theorem 2]{QABL}. From the work \cite{ABLV}, it is also shown that 
$n(x)$ has comparable spatial decay as the harmonic $\Q$.

\subsection{Green's function for classical isotropic Stokes 
system}
We introduce here the fundamental solution $(\E, \q)$ of the 
Stokes system which solves the system of equations,
\begin{equation}
\left\{
\begin{array}{rcl}
-\Delta \E_{ij}(x) - \frac{\partial}{\partial x_i}\q_j(x) & = & \delta_{ij}\delta_0(x),\\
\frac{\partial}{\partial x_i} \E_{ij}(x) & = & 0.
\end{array}
\right.
\end{equation}
Following \cite[Chapter 4.2, p. 238]{galdi2011introduction}, $(\E, \q)$ is given as:
\begin{eqnarray}
\E(\br) &=& \frac{1}{8\pi}\left[
\frac{\id}{r} + \frac{\br\otimes\br}{r^3}
\right],
\quad\text{i.e.}\quad
\E_{ij}(\br) = \frac{1}{8\pi}\left[
\frac{\delta_{ij}}{r} + \frac{x_ix_j}{r^3}
\right],\,\,\,i,j=1,2,3,\label{GreenE}\\
\q(\br) & = & -\frac{1}{4\pi}\frac{\br}{r^3},
\quad\text{i.e.}\quad
\q_i(\br) = -\frac{x_i}{r^3},
\,\,\,i=1,2,3.
\end{eqnarray}
Note that 
\[
\E(\br) \sim \frac{1}{r},
\quad\text{and}\quad
\q(\br) \sim \frac{1}{r^2}
\quad\text{for}\quad
r\gg 1.
\]

The above fundamental solution can be used to produce solutions of Stokes 
system on the whole $\mathbb{R}^3$. More precisely, if $(\bu, p)$ solves
\begin{eqnarray*}
- \Delta \bu + \nabla p &=& {\bf f},
\quad\text{on $\mathbb{R}^3$},\\
\div\bu &=& 0,
\quad\text{on $\mathbb{R}^3$},\\
\bu & = & 0, \quad\text{at $|x|=\infty$},
\end{eqnarray*}
then it is given by:
\begin{equation}
\bu(x) = \int_{\mathbb{R}^3} \E(x-y){\bf f}(y)\,dy
\quad\text{and}\quad
p(x) = \int_{\mathbb{R}^3} -\q(x-y)\cdot {\bf f}(y)\,dy.
\end{equation}
The above integrals are well-defined for $\f$ with sufficient spatial
decay, for example, $\f(\bx) \lesssim 1\wedge\frac{1}{\br^3}$.
For general existence theorems in $L^p(\R^3)$, we refer to 
\cite{galdiARMA, galdi2011introduction}.

Similarly, in the case of a Stokes system in an exterior domain, for example,
\begin{eqnarray*}
- \Delta \bu + \nabla p &=& {\bf f},
\quad\text{on $\Omega$},\\
\div\bu & = & 0,
\quad\text{on $\Omega$,}\\
\bu & = & 0, \quad\text{at $|x|=\infty$},
\end{eqnarray*}
then $\bu$ and $p$ can be represented by:
\begin{eqnarray}
\bu(x)
&=& \int_\Omega \E(x-y)\f(y)\,dy\nonumber\\
&&+ \int_{\partial\Omega}
\big \langle \E(x-y), \T(\bu, p)(y)\nu_y\big\rangle\, d\sigma_y
- \int_{\partial\Omega}\big \langle \T(\E, \q)(x-y)\nu_y, \bu(y)\big \rangle
\,d\sigma_y,\label{u-ext}\\
p(x) 
&=& \int_\Omega - \q(x-y)\cdot\f(y)\,dy\nonumber\\
&&-\int_{\partial\Omega} 
\big\langle \T(\bu, p)(y)\nu_y,\q(x-y)\big\rangle\,d\sigma_y
+ 2\int_{\partial\Omega}\big\langle \nabla_x\q(x-y)\nu_y,\bu(y)\big\rangle
\,d\sigma_y\label{p-ext}.
\end{eqnarray}
In the above, for any given vector and scalar fields $\bf w$ and $\pi$, we 
define the stress tensor as
\begin{equation}\label{Tten1}
\T({\bf w}, \pi) = \nabla {\bf w} + (\nabla{\bf w})^T - \pi\id,
\end{equation}
and for any matrix and vector fields $\mathbf{M}$ and $\mathbf{g}$
and vector $v$, $\T(\mathbf{M}, \mathbf{g})v$ is a matrix with its
$i$-th row given by
\[
\big(\T(\mathbf{M},\mathbf{g})v\big)_i = 
\T(\mathbf{M}_i,\mathbf{g}_i)v.
\]
We also recall the convention as stated in item \ref{notation}
at the end of the Introduction.

In the above and the bulk of this paper, we are dealing with solutions $\bu$
that converge to their far-field limits with rate $\displaystyle \frac{1}{r}$. With this
in mind, we expect the following estimates for the boundary integrals,
\begin{eqnarray}
\left|\int_{\partial\Omega}\big \langle \E(x,-y), \T(\bu, p)(y)\nu_y
\big \rangle\, d\sigma_y\right|
&\lesssim& \frac{1}{r},\label{bdry.int.key}\\
\left|\int_{\partial\Omega}\big \langle \T(\E, \q_j)(x-y)\nu_y, \bu(y)\big \rangle
\,d\sigma_y\right| 
&\lesssim& \frac{1}{r^2},\\
\left|\int_{\partial\Omega}
\big\langle \T(\bu, p)(y)\nu_y,\q(x-y)\big\rangle\,d\sigma_y\right|
& \lesssim & \frac{1}{r^2},\\
\left|\int_{\partial\Omega}\big\langle \nabla_x\q(x-y)\nu_y,\bu(y)\big\rangle
\,d\sigma_y\right|
&\lesssim & \frac{1}{r^3}.\label{bdry.int.4}
\end{eqnarray}
Clearly \eqref{bdry.int.key} gives the dominating far-field behavior. It can be
decomposed as:
\begin{eqnarray}
&&\int_{\partial\Omega}\big \langle \E(x-y), \T(\bu, p)(y)\nu_y
\big \rangle\, d\sigma_y\nonumber\\
& = & 
\Big\langle\E(x), \int_{\partial\Omega} \T(\bu, p)(y)\nu_y d\sigma_y
\Big\rangle + 
\int_{\partial\Omega} 
\big\langle \E(x-y)-\E(x), \T(\bu, p)(y)\nu_y\big\rangle
d\sigma_y
\nonumber\\
& = & 
\langle\E(x), {\mathcal F}\rangle + O\left(\frac{1}{r^2}\right)
\label{bdry.stress.asym}
\end{eqnarray}
where
\begin{equation}\label{Bdry.Stress}
{\mathcal F} := \int_{\partial\Omega} \T(\bu, p)(y)\nu_y d\sigma_y
\end{equation}
denotes the \emph{boundary stress} or \emph{drag force} on the particle. 
The decomposition \eqref{bdry.stress.asym}
is the same as \cite[Theorem 1, Eq. (4.2a)]{ChangFinn}.

Precise far-field asymptotics of the bulk integrals and the
method of solution for the Stokes system 
\eqref{eq:qten2s}--\eqref{eq:inc1s} will be given in Sections \ref{asym.solns.model} and Appendix \ref{AppFFZ}.

\subsection{Uniform Stokes Flow}
Here we provide the solution of a uniform Stokes flow $\U$ with 
far-field velocity $\U_*$, passing a sphere of radius $a$. It solves the following system of 
equations:
\begin{eqnarray}
-\Delta\U + \nabla p &=& 0,\quad\text{for $|x| > a$};
\label{classicStokes0}\\
\div\U &=& 0,\quad\text{for $|x| > a$};\\
\U & = & 0,\quad\text{at $|x|=a$};\\
\U & = & \U_*,\quad\text{at $|x|=\infty$}.
\end{eqnarray}
Using the spherical coordinates (following the physicists' convention)
with $\theta$ being the polar angle (measured from the polar axis) and 
$\phi$ being the azimuthal angle (measured from the meridian plane), we
can write the velocity flowing along the polar axis as $\U=(u_r, u_\theta,
u_\phi=0)$.
Following \cite[Section 7.2]{acheson}, we have
\begin{equation}\label{potential.flow}
u_r = \frac{1}{r^2\sin\theta}\frac{\partial\Psi}{\partial\theta},
\quad\text{and}\quad
u_\theta = -\frac{1}{r\sin\theta}\frac{\partial\Psi}{\partial r},
\end{equation}
where 
\begin{align}
\Psi(r,\theta) &= f(r)\sin^2\theta,
\quad
f(r)=\frac{V}{4}\left(2r^2 -3ar + \frac{a^3}{r}\right),
\quad V=|\U_*|,\label{potential.form}\\
u_r  &=\frac{V}{4}\left(2-\frac{3a}{r}+\frac{a^3}{r^3}\right)2\cos\theta,
\label{ur.form}\\
u_\theta  &=-\frac{V}{4}\left(4-\frac{3a}{r}-\frac{a^3}{r^3}\right)\sin\theta.
\label{uth.form}
\end{align}
We note that the form of $f$ is found by
solving $\displaystyle
\left(\frac{\partial^2}{\partial r^2} - \frac{2}{r^2}\right)^2 f(r) = 0
$. The method of finding $f$ in a bounded (annular) domain will be
presented in Appendix \ref{StokesBdDom}.

If the flow is in the direction of the $x_1$-axis, i.e., $\U_*=Ve_1$, 
then we have in Cartesian coordinates that
\begin{eqnarray}
\U & = & V\left(
{e_1} - \frac{3a}{4}\left(\frac{x_1\br+r^2{e_1}}{r^3}\right)
- \frac{a^3}{4}\left(\frac{r^2{e_1} - 3x_1\br}{r^5}\right)
\right)\nonumber\\
& = & V\left(
{e_1} - \frac{3a}{4}\left(\frac{2x_1^2 + x_2^2 + x_3^2, x_1x_2, x_1x_3}{r^3}\right)
\right.\nonumber\\
&&\left.- \frac{a^3}{4}\left(\frac{-2x_1^2 + x_2^2 + x_3^2, -3x_1x_2, -3x_1x_3}{r^5}\right)
\right)\\
p 
& = & p_\infty - \frac{3}{2}\frac{Va}{r^2}\cos\theta
= p_\infty - \frac{3}{2}\frac{Va}{r^2}\frac{x_1}{\sqrt{x_1^2 + x_2^2}}.
\end{eqnarray}
More generally, in vector form, we have 
\begin{eqnarray}
\U &=& \U_* - \frac{3a}{4}\left(\frac{x\otimes x \U_*+r^2\U_*}{r^3}\right)
- \frac{a^3}{4}\left(\frac{r^2\U_* - 3x\otimes x \U_*}{r^5}\right)\nonumber\\
&=& \left[\I - \frac{3a}{4}\left(\frac{x\otimes x +r^2\I}{r^3}\right)
+ \frac{a^3}{4}\left(\frac{3x\otimes x - r^2\I}{r^5}\right)\right]\U_*
\end{eqnarray}
Making use of the Green's function $\E$ \eqref{GreenE} and upon introducing
\begin{equation}\label{GreenF}
\F(x) = \frac{3x\otimes x - r^2\I}{r^5},
\end{equation}
we have
\begin{eqnarray*}
\U &=& \U_* - 6\pi a\E(x)\U_*
+ \frac{a^3}{4}\F(x)\U_*,
\end{eqnarray*}
or more compactly,
\begin{equation}\label{StokesFormGenVec}
\U = \E_S\U_*,\,\,\,\text{where}\,\,\,
\E_S = \I - 6\pi a\E(x)
+ \frac{a^3}{4}\F(x).
\end{equation}
We note the following asymptotics,
\[
\E(x) \sim O\left(\frac{1}{r}\right),\,\,\,
\F(x) \sim O\left(\frac{1}{r^3}\right),\,\,\,
\]
so that
\begin{equation}
\U \sim \U_* + O\left(\frac{1}{r}\right),
\,\,\,\text{and}\,\,\,
\nabla \U,\,\A,\,\W \sim O\left(\frac{1}{r^2}\right),
\,\,\,\text{for $r\gg 1$.}
\end{equation}

Next we compute the drag force 
${\mathcal F}$ \eqref{Bdry.Stress} on the moving particle.
For this purpose, the components of $\T(\bu, p)$ on $\partial\bB_a$, expressed in spherical coordinates are:
\begin{eqnarray*}
\T_{rr} & = & -p + 2\frac{\partial u_r}{\partial r}
= -p_\infty + \frac{3}{2}\frac{V}{a}\cos\theta,\\
\T_{r\theta} & = & r\frac{\partial}{\partial r}\left(\frac{u_\theta}{r}\right)
+ \frac{1}{r}\frac{\partial u_r}{\partial\theta}
= -\frac{3}{2}\frac{U}{a}\sin\theta,\\
\T_{r\phi} & = & 0.
\end{eqnarray*}
With the above, $\mathcal F$,
in the direction of $\U_*$, is given by
\begin{eqnarray}
{\mathcal F} = \int_0^{2\pi}\int_0^\pi\Big(
T_{rr}\cos\theta - \T_{r\theta}\sin\theta
\Big)a^2\sin\theta d\theta d\phi
= 6\pi aV
\end{eqnarray}
or in vector form, written as,
\begin{equation}
{\mathcal F} = 6\pi a\U_*
\end{equation}
which is the celebrated Stokes Law.

\begin{remark}
The classical isotropic \emph{Stokeslet} fundamental solution given by
\[\bv_s=\E\boldsymbol\alpha,\quad p_s=\q\cdot\boldsymbol\alpha,\] 
for any $\boldsymbol\alpha\in\mathbb{R}^3$ corresponds to Stokes flow driven by a point source of strength $\boldsymbol\alpha$ located at the origin. Note that the far-field behavior of the Stokeslet is the same as the leading order asymptotics of classical Stokes flow in the exterior of a spherical domain. See \cite{Chwang} for a explanation and application of such a notion.
\end{remark}

\section{Structure of the anisotropic Stokes equation \eqref{eq:velo1s}}
\label{Structure.Eqn}
Here we will start our analysis for the stationary system
\eqref{eq:qten2s}--\eqref{eq:inc1s} with boundary conditions
\eqref{eq:ext1}--\eqref{eq:ext2}. We emphasize again the feature that 
the tensor field $\Q$ given by \eqref{QABL}--\eqref{QABL3} 
acts as an inhomogeneous term for the Stokes equation 
\eqref{eq:velo1s} for $\bv$. 
The key is to understand the far-field spatial behavior 
($r\gg 1$) of $\bv$.

\subsection{Decomposition of equation \eqref{eq:velo1s}}
\label{DecompSys}
Here we decompose system \eqref{eq:qten2}--\eqref{eq:inc1} into a form 
amenable for asymptotic analysis. 
We first substitute form \eqref{QABL} for $\Q$ into \eqref{eq:velo1s} for 
$\bv$ and analyze the resulting system. 
\begin{enumerate}
\item Recall the decay property \eqref{Q.farfield} of $\Q$:
$\Q = \Q_* + O\left(\frac{1}{r}\right)$ and
$\nabla\Q = O\left(\frac{1}{r^2}\right)$.
We look for a solution $\bv$ satisfying
$\bv = \bv_* + O\left(\frac{1}{r}\right)$ and 
$\nabla \bv = O\left(\frac{1}{r^2}\right)$. Hence, we expect
\begin{eqnarray*}
\A,\,\,\,\W& \sim & O\left(\frac{1}{r^2}\right),\\
\Qc 
& = & \bv\cdot\nabla \Q + \Q\W - \W\Q\\
& = & 
\underset{O\left(\frac{1}{r^2}\right)}
{\underbrace{\bv_*\cdot\nabla\Q +\Q_*\W - \W\Q_*}}
+
\underset{O\left(\frac{1}{r^3}\right)}
{\underbrace{(\bv-\bv_*)\cdot\nabla\Q +(\Q-\Q_*)\W - \W(\Q-\Q_*)}}\\
& \sim & O\left(\frac{1}{r^2}\right).
\end{eqnarray*}

\item	With the above, the right hand side of \eqref{eq:velo1s} becomes
\begin{equation}
\cD:=
-\left(\gamma_1\Qc+\gamma_2\A+\frac{\gamma_3}{2}(\A\Q+\Q\A)+\gamma_9(\A\cdot\Q)\Q\right)\cdot\partial_i\Q{\mathbf e}_i
\label{D}
\sim O\left(\frac{1}{r^4}\right).
\end{equation}

\item Here we analyze the terms constituting the viscous stress 
$\T^\mathrm{v}_\text{SV}$ given by \eqref{eq:visc5-2}.
\begin{enumerate}
\item\label{item1} The $\gamma_1$-term:
\begin{eqnarray}
\left[\Q\Qc - \Qc\Q\right]
& = & 
\left[
\Q\left(\bv\cdot\nabla \Q + \Q\W - \W \Q\right)
-\left(\bv\cdot\nabla \Q + \Q\W - \W\Q\right)\Q
\right]\nonumber\\
& = & 
\left[
\Q\left(\bv\cdot\nabla \Q\right) -\left(\bv\cdot\nabla \Q\right)\Q
\right]+
\left[
\Q^2\W - 2\Q\W\Q  + \W\Q^2
\right].\label{g12}
\end{eqnarray}
For the first bracketed term in \eqref{g12}, we have the following decomposition
\begin{eqnarray}
\label{eq:ag1}
\left[
\Q\left(\bv\cdot\nabla\Q\right) -\left(\bv\cdot\nabla\Q\right)\Q
\right]
&=& 
\left[
\Q_*(\bv_*\cdot\nabla\Q) - (\bv_*\cdot\nabla\Q)\Q_*\right]\nonumber \\
& & 
+ 
\left[\Q_*(\bv-\bv_*)\cdot\nabla\Q - (\bv-\bv_*)\cdot(\nabla\Q)\Q_*
\right]\nonumber\\
& & 
+\left[(\Q-\Q_*)\bv\cdot\nabla\Q - (\bv\cdot \nabla\Q)(\Q-\Q_*)
\right]\nonumber\\
& = & 
\left[
\Q_*(\bv_*\cdot\nabla\Q) - (\bv_*\cdot\nabla\Q)\Q_*
\right] + O\left(\frac{1}{r^3}\right)\nonumber\\
&=:&\mathcal{A}_1(\br)+ O\left(\frac{1}{r^3}\right).
\end{eqnarray}
For the second bracketed term in \eqref{g12}, we have
\begin{eqnarray}
\label{eq:ag2}
\left[\Q^2\W - 2\Q\W\Q  + \W\Q^2\right]
&=&\left[\Q_*^2\W - 2\Q_*\W\Q_*  + \W\Q^2_*\right]
+ O\left(\frac{1}{r^3}\right)\nonumber\\&=:&\B_1[\nabla\bv]+ O\left(\frac{1}{r^3}\right).
\end{eqnarray}
In the above, $\mathcal{A}_1$ and $\B_1$ are 
linear in $\bv_*\cdot\nabla\Q$ and $\nabla\bv$ respectively, 
both having coefficients depending only on $\Q_*$. 
Furthermore both $\cA_1$ and $\B_1$ decay as $r^{-2}$.
A similar structure exists for all the remaining terms in the stress tensor 
as it will become clear while we proceed through the rest of this computation.
\item The $\gamma_2$-term:
\begin{eqnarray*}
&&\left[\Qc + \Q\A-\A\Q\right]\\
& = & 
\left[\bv\cdot\nabla\Q + \Q\nabla\bv-(\nabla\bv)\Q\right]\\
& = & \left[
\bv_*\cdot\nabla\Q + \Q_*\nabla\bv - (\nabla\bv)\Q_*
\right]
+ \left[
(\bv-\bv_*)\cdot\nabla\Q + (\Q-\Q_*)\nabla\bv - \nabla\bv (\Q-\Q_*)
\right]\\
& = &\mathcal{A}_2(\br)+\B_2[\nabla\bv]+ O\left(\frac{1}{r^3}\right),
\end{eqnarray*}
where
\begin{equation}
\label{eq:bg2}
\mathcal{A}_2(\br):=\bv_*\cdot\nabla\Q
\,\,\mbox{ and }\,\,\B_2[\nabla\bv]:=\Q_*\nabla\bv - (\nabla\bv)\Q_*.
\end{equation}
\item	The $\gamma_3$-term:
\begin{eqnarray*}
& & \frac{1}{2}\left(\Q\Qc+\Qc\Q +\Q^2\A-\A\Q^2\right)\nonumber\\
& = & 
\frac{1}{2}\Big(\Q_*\big(\bv_*\cdot\nabla\Q + \Q_*\W-\W\Q_*\big)
+ \big(\bv_*\cdot\nabla\Q + \Q_*\W-\W\Q_*\big)\Q_*
\nonumber\\&&+\Q_*^2\A - \A\Q_*^2\Big) + O\left(\frac{1}{r^3}\right)\nonumber\\
& = & 
\frac{1}{2}\left(\Q_*(\bv_*\cdot\nabla\Q) + (\bv_*\cdot\nabla\Q)\Q_*+\Q_*^2(\nabla\bv) - (\nabla\bv)\Q_*^2\right)+ O\left(\frac{1}{r^3}\right)\nonumber \\
& = &\mathcal{A}_3(\br)+\B_3[\nabla\bv]+ O\left(\frac{1}{r^3}\right),
\end{eqnarray*}
where
\begin{equation}
\label{eq:bg3}
\begin{aligned}
\mathcal{A}_3(\br)&:=\frac{1}{2}\left(\Q_*(\bv_*\cdot\nabla\Q) + (\bv_*\cdot\nabla\Q)\Q_*\right),\\\B_3[\nabla\bv]&:=\frac{1}{2}\left(\Q_*^2(\nabla\bv) - (\nabla\bv)\Q_*^2\right).
\end{aligned}
\end{equation}
\item	The $\gamma_4$-term:
\begin{equation}
\label{eq:bg4}
\Q\A+\A\Q = \Q_*\A + \A\Q_* + O\left(\frac{1}{r^3}\right)=:\B_4[\nabla\bv]+ O\left(\frac{1}{r^3}\right).
\end{equation}
\item	The $\gamma_5$-term:
\begin{equation}
\label{eq:bg5}
\Q^2\A+\A\Q^2 = \Q_*^2\A + \A\Q_*^2 + O\left(\frac{1}{r^3}\right)=:\B_5[\nabla\bv]+ O\left(\frac{1}{r^3}\right).
\end{equation}

\item	The $\gamma_6$-term:
\begin{equation}
\label{eq:bg6}
(\A\cdot\Q)\Q = (\A\cdot\Q_*)\Q_* + O\left(\frac{1}{r^3}\right)=:\B_6[\nabla\bv]+ O\left(\frac{1}{r^3}\right).
\end{equation}

\item\label{item7}The $\gamma_7$-term:
\begin{equation}
\label{eq:bg7}
|\Q^2|\A = |\Q^2_*|\A + O\left(\frac{1}{r^3}\right)
=:\B_7[\nabla\bv]+ O\left(\frac{1}{r^3}\right).
\end{equation}

\item	The $\gamma_9$-term:
\begin{eqnarray}
(\Qc\cdot\Q)\Q 
& = &
\Big((\bv\cdot\nabla\Q + \Q\W - \W\Q)\cdot\Q\Big)\Q
=
\Big((\bv\cdot\nabla\Q)\cdot\Q\Big)\Q
\nonumber\\
& = & 
\Big((\bv_*\cdot\nabla\Q)\cdot\Q_*\Big)\Q_*
+ O\left(\frac{1}{r^3}\right)\nonumber \\
& = &\mathcal{A}_9(\br)+ O\left(\frac{1}{r^3}\right),
\end{eqnarray}
where we have used the fact that
$(\Q\W - \W\Q)\cdot\Q=0$ as $(\C\D-\D\C)\cdot \C=0$ for any 
$\C,\D\in M^{3\times3}$ with $\C$ symmetric.
\item	The $\gamma_{10}$-term:
\begin{eqnarray}
\label{eq:bg10}
(\Q^2\cdot\A)\Q + (\Q\cdot\A)\Q^2
&=& 
(\Q_*^2\cdot\A)\Q_* + (\Q_*\cdot\A)\Q^2_*
+ O\left(\frac{1}{r^3}\right)\nonumber\\&=:&\B_{10}[\nabla\bv]+ O\left(\frac{1}{r^3}\right).
\end{eqnarray}

\item\label{item11}The $\gamma_{11}$-term:
\begin{eqnarray}
\label{eq:bg11}
|\Q|^2(\A\cdot\Q)\Q 
&=& |\Q_*|^2(\A\cdot\Q_*)\Q_* + O\left(\frac{1}{r^3}\right)
=:\B_{11}[\nabla\bv]+ O\left(\frac{1}{r^3}\right).
\end{eqnarray}
\end{enumerate}
From the above, we note that
\begin{equation}
\label{eq:azero}
\mathcal{A}_4=\mathcal{A}_5=\mathcal{A}_6=\mathcal{A}_7=\mathcal{A}_{10}=\mathcal{A}_{11} = 0\,\,\text{and}\,\,
\B_9 = 0.
\end{equation}
\end{enumerate}

Taking into account the items (a)-(c) above, 
we can write equation \eqref{eq:velo1s} as
\begin{equation}
\label{eq:nolabel}
- \left(
\Delta\bv 
+\text{div}\Big[\B_\gamma[\nabla\bv] + \mathcal{A}_\gamma(\br) + \mathcal{C}_\gamma(\br)\Big]
\right) 
+ \nabla p 
= 
\mathcal{D}_\gamma(\br)
\end{equation}
where
\begin{eqnarray}
\label{eq:I}
\B_\gamma[\nabla \bv] & =& 
\B_{\gamma, \Q_*}[\nabla \bv] := 
\sum_{i=1,\neq8,9}^{11}\gamma_i\B_i[\nabla\bv],\\
\mathcal{A}_\gamma(\br) &=&
\mathcal{A}_{\gamma, \Q_*}[\bv_*\cdot\nabla\Q]
:=\sum_{i=1,2,3,9}\gamma_i\mathcal{A}_i(x) = O\left(\frac{1}{r^2}\right),\\
\cC_\gamma(\br) &= & \cC_\gamma(\bv, \br) :=  
\T_{\text{SV}}^\text{v} - \A - \B_\gamma[\nabla\bv] - \cA_\gamma(\br)
= O\left(\frac{1}{r^3}\right),\label{C_gamma_form}\\
\cD_\gamma(\br) & = & \cD_\gamma(\bv, \br) := 
-\left(\gamma_1\Qc+\gamma_2\A+\frac{\gamma_3}{2}(\A\Q+\Q\A)+\gamma_9(\A\cdot\Q)\Q\right)\cdot\partial_i\Q{\mathbf e}_i\\
&=&O\left(\frac{1}{r^4}\right).\nonumber
\label{D_gamma_asym}
\end{eqnarray}
In the above and what follows, we will use the symbol $\gamma$ to 
denote an expression or quantity that genuinely depends on 
$\gamma_1,\ldots, \gamma_{11}$. Furthermore,
we set
\begin{equation}
|\gamma| := \max\{|\gamma_1|,\ldots,|\gamma_{11}|\}.
\end{equation}
We might omit $\gamma$ if the dependence 
is clear from the context.

The advantage of representation \eqref{eq:nolabel}--\eqref{D_gamma_asym}
is highlighted as follows.
\begin{enumerate}
\item
The linear form $\B_\gamma[\cdot]$ in $\nabla\bv$ corresponds to the leading 
order 
far-field contribution to the diffusivity matrix originating from the stress 
tensor. In particular, $\div[\B_\gamma[\nabla v]]$ is linear in $D^2\bv$ with
constant coefficients depending only on $\Q_*$. 
\item
The term $\div[\cA_\gamma(\br)]$ decays as $r^{-3}$ for $r \gg 1$. 
It can be treated as a purely inhomogeneous 
forcing term involving $D^2\Q$, $\bv_*$ and $\Q_*$. Note that it is linear in 
the expression $\bv_*\cdot\nabla\Q$ so it vanishes if $\bv_*=0$.
\item	The terms $\div[\cC_\gamma(x)]$ and $\cD_\gamma(\br)$ decay as
$r^{-4}$ for $r\gg 1$. They are integrable in the exterior domain $\Omega$,
\begin{equation}
\int_\Omega \big|\div[\cC_\gamma(\br)]\big|\,dx,\,\,\,
\int_\Omega \big|\cD_\gamma(\br)\big|\,dx < \infty.
\end{equation}
\end{enumerate}
We point out that the presence of $\cA_\gamma$ and $\cD_\gamma$ is due to 
the dependence of the interacting potential function $R$ on $\Qc$ so that
only $\gamma_1, \gamma_2, \gamma_3$ and $\gamma_9$ appear in the expressions
for $\cA_\gamma$ and $\cD_\gamma$. (See the form \eqref{eq:diss} of $R$.)

For the purpose of analyzing equation \eqref{eq:nolabel},
we will keep $\div\big[\B_\gamma[\nabla\bv]\big]$ in the left hand
side but move
$\cA_\gamma, \cC_\gamma, \cD_\gamma$ to the right of that equation and
re-write it as:
\begin{equation}
\label{uEqnM}
- \Delta\bv 
-\div\big[\B_\gamma[\nabla\bv]\big]
+ \nabla p
=
\text{div}\big[\mathcal{A}_\gamma(\br)\big]
+ \text{div}\big[\mathcal{C}_\gamma(\br)\big]
+\mathcal{D}_\gamma(\br).
\end{equation}
Note that $-\Delta\bv -\div\big[\B_\gamma[\nabla\bv]\big]$ 
is a second order differential operator in $\bv$ with constant coefficients.

Before proceeding further, we will express $\cC_\gamma$ and $\cD_\gamma$ in more transparent forms.
\begin{enumerate}
\item For $\cC_\gamma$, from \eqref{C_gamma_form}, we compute
\begin{align*}
\cC_\gamma
\,\,= \,\,& \T_{\text{SV}}^\bv - \A - \B_\gamma[\nabla\bv]] - \cA_\gamma(\br)\\
\,\,= \,\,& \gamma_1
\Big[
\big(\Q\left(\bv\cdot\nabla \Q\right) -\left(\bv\cdot\nabla \Q\right)\Q\big) 
+\big(\Q^2\W - 2\Q\W\Q  + \W\Q^2\big)\\
& \hspace{15pt}
- \big(\Q_*(\bv_*\cdot\nabla\Q) - (\bv_*\cdot\nabla\Q)\Q_*\big)
- \big(\Q_*^2\W - 2\Q_*\W\Q_*  + \W\Q^2_*\big)
\Big]\\
&+\gamma_2\Big[
\big(\bv-\bvs\big)\cdot\nabla\Q 
+ \big(\Q-\Q_*\big)\nabla\bv
-(\nabla\bv)\big(\Q-\Q_*\big)
\Big]\\
&+\frac{\gamma_3}{2}
\Big[
\left(\Q(\bv\cdot\nabla\Q)+ (\bv\cdot\nabla\Q)\Q
+\Q^2\nabla\bv-(\nabla\bv)\Q^2\right)\\
&\hspace{25pt}
-
\left(\Q_*(\bv_*\cdot\nabla\Q) + (\bv_*\cdot\nabla\Q)\Q_*\right)
-
\left(\Q_*^2(\nabla\bv) - (\nabla\bv)\Q_*^2\right)
\Big]\\
& +\gamma_4\Big[(\Q-\Q_*)\A+\A(\Q-\Q_*)\Big]\\
& +\gamma_5\Big[(\Q^2-\Q_*^2)\A+\A(\Q^2-\Q_*^2)\Big]\\
& +\gamma_6\Big[
(\A\cdot\Q)\Q - (\A\cdot\Q_*)\Q_*
\Big]\\
& +\gamma_7\Big[\left(|\Q^2| - |\Q^2_*|\right)\A\Big]\\
& + \gamma_9\Big[
\big((\bv\cdot\nabla\Q)\cdot\Q\big)\Q 
-\Big((\bv_*\cdot\nabla\Q)\cdot\Q_*\Big)\Q_*
\Big]\\
& + \gamma_{10}\Big[
(\Q^2\cdot\A)\Q + (\Q\cdot\A)\Q^2
-(\Q_*^2\cdot\A)\Q_* - (\Q_*\cdot\A)\Q^2_*
\Big]\\
& + \gamma_{11}\Big[
|\Q|^2(\A\cdot\Q)\Q - |\Q_*|^2(\A\cdot\Q_*)\Q_*
\Big]\\
\,\,\sim \,\,&
O\left(\frac{1}{r^3}\right).
\end{align*}

\item
For $\cD_\gamma$, from \eqref{D}, we have,
\begin{eqnarray*}
\cD_\gamma&=&
-\Big[\gamma_1\big(
\bv_*\cdot\nabla\Q +\Q_*\W - \W\Q_*
+(\bv-\bv_*)\cdot\nabla\Q +(\Q-\Q_*)\W - \W(\Q-\Q_*)
\big)\nonumber\\
&&\hspace{12pt}+\gamma_2\A
+ \frac{\gamma_3}{2}\big(\A\Q_* + \Q_*\A + \A(\Q-\Q_*) + (\Q-\Q_*)\A\big)
\nonumber\\
&&\hspace{12pt}+
\gamma_9\big(
(\A\cdot\Q_* + \A\cdot(\Q-\Q_*))\Q_*
+(\A\cdot\Q_* + \A\cdot(\Q-\Q_*))(\Q-\Q_*)
\big)
\Big]
\cdot\partial_i\Q{\mathbf e}_i\\
&\sim&O\left(\frac{1}{r^4}\right).
\end{eqnarray*}
\end{enumerate}
Note that given $\Q$, both $\cC_\gamma$ and $\cD_\gamma$ are linear in 
$\bv$.

Combining the above with the explicit expression \eqref{QABL} of $\Q$, 
the terms on the right hand side of \eqref{uEqnM} take the following forms:
\begin{eqnarray}
\div\cA_\gamma(\br) 
& = & \gamma\frac{{\bf F}(\hat{\br})}{r^3} + O\left(\frac{1}{r^4}\right),\label{eq:A}\\
\div\cC_\gamma(\br) + \cD_\gamma(\br) 
&=&
\gamma\left(\frac{{\bf G}(\hat{\br})}{r^4} 
+ \frac{{\bf H}(\hat{\br}):(\bv-\bv_*)}{r^3} 
+ \frac{{\bf I}(\hat{\br}):D\bv}{r^2}
+ \frac{{\bf J}(\hat{\br}):D^2\bv}{r}
\right)\nonumber\\
& & 
+ O\left(\frac{1}{r^5}\right)
\label{eq:GIJH}
\end{eqnarray}
for some spatially bounded vector or tensor fields $\bf F, G, H, I, J$ 
defined on $\mathbb{S}^2$ which depend on $\Q$ and $\Q_*$ but not on
$\bv$ and $\bv_*$. 
We again recall the convention about the symbol $``:"$ stated in 
item \ref{notation} at the end of Introduction. The precise
forms of the contractions between tensors are not too important for our 
analysis.
The key is the homogeneity in $r$ leading to appropriate spatial decays.
The presence of $\gamma$ on the right hand side refers to the fact that
the terms are multiplied by $\gamma_i$'s.
In particular, $\cA_\gamma, \cB_\gamma, \cC_\gamma, \cD_\gamma$ are all 
bounded in magnitude by $O(|\gamma|)$.
Hence, if $|\gamma| = 0$, then \eqref{uEqnM} simply
becomes the classical isotropic Stokes equation \eqref{classicStokes0}.

Relating back to Section \ref{sec:BCs}, Remark \ref{vBC.rmk},
note that the system \eqref{eq:qten2s}--\eqref{eq:inc1s} with boundary conditions
\eqref{eq:ext1}--\eqref{eq:ext2} is solvable for any $\bv_b$ and $\bv_*$.
From the form of $\cA_\gamma$ written in \eqref{eq:A}, it seems the system
becomes simpler by setting $\bv_*=0$. This can be achieved by a change of 
frame or simply consider the new vector field $\bv' = \bv-\bv_*$. However, 
due to the presence of the convective derivative $\bv\cdot\nabla\Q$ of $\Q$,
either of these procedures will necessarily give rise to the term 
$\bv_*\nabla\Q$ which at leading order is embedded in $\cA_\gamma$.

\subsection{Identification of $\div\big[\B_\gamma[\nabla\bv]\big]$}
\label{IdentifyM}
In this section, we express $\div\big[\B_\gamma[\nabla\bv]\big]$ in 
\eqref{uEqnM} in more explicit form. More precisely, we will write
\begin{equation}\label{uEqnM1}
\div\big[\B_\gamma[\nabla\bv]\big]
= \cM_\gamma:D^2\bv
\end{equation}
for some constant fourth order tensor $\M_\gamma$. Upon 
introducing the coordinates $\cM_\gamma = (\cM_{i,j;k,l})$, the 
right hand side of the above is understood as
\begin{equation}\label{uEqnM2}
\Big(\cM_\gamma:D^2\bv\Big)_i
= \M_{i,j;k,l}\partial_{kl}\bv_j,\quad
\text{for $i=1,2,3$.}
\end{equation}
The main purpose of this section is to identify $\cM_{i,j;k,l}$ 
explicitly. We note the following symmetry action of $\M_{i,j;k,l}$ 
with respect to $k$ and $l$, 
\begin{equation}\label{M.sym}
\M_{i,j;k,l}\partial_{kl}\bv_j = \M_{i,j;l,k}\partial_{kl}\bv_j.
\end{equation}
Furthermore, we emphasize that $\cM_\gamma$ will only act on 
\emph{incompressible} vector fields $\bv$: $\partial_j\bv_j = 0$.

Before proceeding, using \eqref{QABL2}, we record that
\begin{equation}
{\Q_*}_{ij} = s_*\left(\n_i\n_j - \frac13\delta_{ij}\right),\,\,\,
{\Q_*^2}_{ij} = \frac{s_*^2}{3}\left(\n_i\n_j + \frac13\delta_{ij}\right),\,\,\,
|\Q_*|^2 = \frac{2s_*^2}{3}.
\end{equation}

With that, we compute.
\begin{enumerate}
\item $\mathrm{div}(\B_1[\nabla\bv]):$
\begin{eqnarray*}
&&\div\left(\left[\Q_*^2\W+\W\Q_*^2-2\Q_*\W\Q_*\right]_i\right)
=\partial_j\Big[\Q_*^2\W+\W\Q_*^2-2\Q_*\W\Q_*\Big]_{ij}\\
&=&\partial_j\Big[
\frac{1}{2}\left(\Q_*^2\right)_{ik}(\partial_j\bv_k - \partial_k\bv_j)
+\frac{1}{2}\left(\Q_*^2\right)_{kj}\big(\partial_k\bv_i - \partial_i\bv_k\Big)
- {\left(\Q_*\right)}_{ik}
(\partial_l\bv_k - \partial_k\bv_l)
{\left(\Q_*\right)}_{lj}
\Big]
\\
&=&\frac{1}{2}\left(\Q_*^2\right)_{ik}\partial_{jj}\bv_k+\frac{1}{2}\left(\Q_*^2\right)_{kj}\Big(\partial_{kj}\bv_i - \partial_{ij}\bv_k\Big)-{\left(\Q_*\right)}_{ik}{\left(\Q_*\right)}_{lj}
\big(\partial_{jl}\bv_k - \partial_{jk}\bv_l\big)\\
&=&\left(\frac{1}{2}\delta_{k,l}\left(\Q_*^2\right)_{ij}+\frac{1}{2}\delta_{i,j}\left(\Q_*^2\right)_{kl} - \frac{1}{2}\delta_{i,l}\left(\Q_*^2\right)_{jk}-{\left(\Q_*\right)}_{ij}{\left(\Q_*\right)}_{lk} + {\left(\Q_*\right)}_{ik}{\left(\Q_*\right)}_{lj}\right)\partial_{kl}\bv_j\\
&=:&\M^1_{i,j;k,l}\partial_{kl}\bv_j,
\end{eqnarray*}
where
\begin{eqnarray}
\M^1_{i,j;k,l}&=&s_*^2\left[\frac{1}{6}\delta_{k,l}\left(\n_i\n_j+\frac13\delta_{i,j}\right)+\frac{1}{6}\delta_{i,j}\left(\n_k\n_l+\frac13\delta_{k,l}\right) - \frac{1}{6}\delta_{i,l}\left(\n_k\n_j+\frac13\delta_{k,j}\right)\right.\nonumber\\&&\left.-{\left(\n_i\n_j-\frac13\delta_{i,j}\right)}{\left(\n_k\n_l-\frac13\delta_{k,l}\right)}+ {\left(\n_k\n_i-\frac13\delta_{k,i}\right)}{\left(\n_j\n_l-\frac13\delta_{j,l}\right)}\right]\nonumber\\
&=&\frac{s_*^2}{2}\left(
\delta_{kl}\n_i\n_j + \delta_{ij}\n_k\n_l -\delta_{ik}\n_l\n_j
\right).
\label{eq:m1}
\end{eqnarray}
In the above, we have used the symmetry property \eqref{M.sym}
of $\M$ and the incompressibility of $\bv$:
\begin{equation}
\delta_{il}\n_k\n_j\partial_{kl}\bv_j = \delta_{ik}\n_l\n_j\partial_{kl}\bv_j
\,\,\,\text{and} \,\,\,
\delta_{il}\delta_{kj}\partial_{kl}\bv_j 
= \delta_{ik}\delta_{lj}\partial_{kl}\bv_j
= \partial_{ij}\bv_j = 0.
\end{equation}
The above property will be used in several places in what follows.

\item $\mathrm{div}(\B_2[\nabla\bv]):$
\begin{eqnarray*}
\div\left[\Q_*\nabla\bv-(\nabla\bv)\Q_*\right]_i
& = & \partial_j\left[{\left(\Q_*\right)}_{ik}\partial_j\bv_k-\partial_k\bv_i{\left(\Q_*\right)}_{kj}\right]\\
& = & {\left(\Q_*\right)}_{ik}\partial_{jj}\bv_k- {\left(\Q_*\right)}_{kj}\partial_{kj}\bv_i\\& = & \left(\delta_{k,l}{\left(\Q_*\right)}_{ij}- \delta_{i,j}{\left(\Q_*\right)}_{kl}\right)\partial_{kl}\bv_j\\&=:&\M^2_{i,j;k,l}\partial_{kl}\bv_j.
\end{eqnarray*}
Here
\begin{equation}
\label{eq:m2}
\M^2_{i,j;k,l}=s_*\left(\delta_{k,l}\n_i\n_j- \delta_{i,j}\n_k\n_l\right).
\end{equation}
\item $\mathrm{div}(\B_3[\nabla\bv]):$
\begin{eqnarray*}
\frac{1}{2}\div\left[\Q_*^2\nabla\bv-(\nabla\bv)\Q_*^2\right]_i
& = & \frac{1}{2}\left({\left(\Q_*^2\right)}_{ik}\partial_{jj}\bv_k- \left({\Q_*^2}\right)_{kj}\partial_{kj}\bv_i\right)\\& = & \frac{1}{2}\left(\delta_{k,l}{\left(\Q_*^2\right)}_{ij}- \delta_{i,j}\left({\Q_*^2}\right)_{kl}\right)\partial_{kl}\bv_j\\&=:&\M^3_{i,j;k,l}\partial_{kl}\bv_j.
\end{eqnarray*}
Here
\begin{equation}
\label{eq:m3}
\M^3_{i,j;k,l}=\frac{s_*^2}{6}\left(\delta_{k,l}\n_i\n_j- \delta_{i,j}\n_k\n_l\right).
\end{equation}

\item $\mathrm{div}(\B_4[\nabla\bv]):$
\begin{eqnarray*}
\div\left[\Q_*\A+\A\Q_*\right]_i
& = & \frac{1}{2}\left(\Q_*\right)_{ik}\partial_{jj}\bv_{k}+\frac{1}{2}\left(\partial_{ij}\bv_k+\partial_{kj}\bv_i\right)\left(\Q_*\right)_{kj}\\& = & \frac{1}{2}\left(\delta_{k,l}\left(\Q_*\right)_{ij}+\delta_{i,l}\left(\Q_*\right)_{jk}+\delta_{i,j}\left(\Q_*\right)_{kl}\right)\partial_{kl}\bv_j\\&=:&\M^4_{i,j;k,l}\partial_{kl}\bv_j.
\end{eqnarray*}
Here
\begin{eqnarray}
\M^4_{i,j;k,l}&=&\frac{s_*}{2}\left(\delta_{k,l}\left(\n_i\n_j-\frac13\delta_{i,j}\right)+\delta_{i,l}\left(\n_j\n_k-\frac13\delta_{j,k}\right)+\delta_{i,j}\left(\n_k\n_l-\frac13\delta_{k,l}\right)\right)\nonumber\\
&=&
\frac{s_*}{2}\left(
\delta_{kl}\n_i\n_j + \delta_{ij}\n_k\n_l + \delta_{ik}\n_l\n_j
-\frac23\delta_{kl}\delta_{ij}\right).
\label{eq:m4}
\end{eqnarray}

\item $\mathrm{div}(\B_5[\nabla\bv]):$
\begin{eqnarray*}
\div\left[\Q_*^2\A+\A\Q_*^2\right]_i
& = & \frac{1}{2}\left(\Q_*^2\right)_{ik}\partial_{jj}\bv_{k}+\frac{1}{2}\left(\partial_{ij}\bv_k+\partial_{kj}\bv_i\right)\left(\Q_*^2\right)_{kj}\\& = & \frac{1}{2}\left(\delta_{k,l}\left(\Q_*^2\right)_{ij}+\delta_{i,l}\left(\Q_*^2\right)_{jk}+\delta_{i,j}\left(\Q_*^2\right)_{kl}\right)\partial_{kl}\bv_j\\&=:&\M^5_{i,j;k,l}\partial_{kl}\bv_j.
\end{eqnarray*}
Here
\begin{eqnarray}
\M^5_{i,j;k,l}
&=&\frac{s_*^2}{6}\left(\delta_{k,l}\left(\n_i\n_j+\frac13\delta_{i,j}\right)+\delta_{i,l}\left(\n_j\n_k+\frac13\delta_{j,k}\right)+\delta_{i,j}\left(\n_k\n_l+\frac13\delta_{k,l}\right)\right)\nonumber\\
&=&\frac{s_*^2}{6}\left(
\delta_{kl}\n_i\n_j + \delta_{ij}\n_k\n_l 
+ \delta_{ik}\n_l\n_j
+\frac23\delta_{ij}\delta_{kl} 
\right).
\label{eq:m5}
\end{eqnarray}

\item $\mathrm{div}(\B_6[\nabla\bv]):$
\begin{eqnarray*}
\div\left[(\A\cdot\Q_*)\Q_*\right]_i
& = & \frac{1}{2}\left(\Q_*\right)_{ij}\left(\Q_*\right)_{kl}\left(\partial_{lj}\bv_{k}+\partial_{kj}\bv_l\right)\\& = & \frac{1}{2}\left(\left(\Q_*\right)_{ik}\left(\Q_*\right)_{jl}+\left(\Q_*\right)_{il}\left(\Q_*\right)_{kj}\right)\partial_{kl}\bv_j\\&=:&\M^6_{i,j;k,l}\partial_{kl}\bv_j.
\end{eqnarray*}
Here
\begin{eqnarray}
\M^6_{i,j;k,l}&=&\frac{s_*^2}{2}\left(\left(\n_i\n_k-\frac13\delta_{i,k}\right)\left(\n_j\n_l-\frac13\delta_{j,l}\right)\right.\nonumber\\&&\left.+\left(\n_i\n_l-\frac13\delta_{i,l}\right)\left(\n_j\n_k-\frac13\delta_{j,k}\right)\right)\nonumber\\
\label{eq:m6}
&=&
s_*^2\left(
\n_i\n_j\n_k\n_l - \frac13\delta_{ik}\n_l\n_j
\right).
\end{eqnarray}

\item $\mathrm{div}(\B_7[\nabla\bv]):$
\begin{eqnarray*}
\div\left[|\Q_*|^2\A\right]_i
=\div\left[\frac{2s_*^2}{3}\A\right]_i
=:\M^7_{i,j;k,l}\partial_{kl}\bv_j.
\end{eqnarray*}
Here
\begin{equation}
\label{eq:m7}
\M^7_{i,j;k,l}=\frac{s_*^2}{3}\delta_{k,l}\delta_{i,j}.
\end{equation}

\item $\mathrm{div}(\B_{10}[\nabla\bv]):$
\begin{eqnarray*}
\div\left[(\Q_*^2\cdot\A)\Q_* + (\Q_*\cdot\A)\Q^2_*\right]_i
& = & \frac{1}{2}\left(\left(\Q_*^2\right)_{kl}\left(\Q_*\right)_{ij}+\left(\Q_*\right)_{kl}\left(\Q_*^2\right)_{ij}\right)\left(\partial_{lj}\bv_{k}+\partial_{kj}\bv_{l}\right)\\& = & \left(\left(\Q_*^2\right)_{jk}\left(\Q_*\right)_{il}+\left(\Q_*\right)_{jk}\left(\Q_*^2\right)_{il}\right)\partial_{kl}\bv_j\\&=:&\M^{10}_{i,j;k,l}\partial_{kl}\bv_j.
\end{eqnarray*}
Here
\begin{eqnarray}
\M^{10}_{i,j;k,l}
&=&\frac{s_*^3}{3}\left(\left(\n_j\n_k+\frac13\delta_{j,k}\right)\left(\n_i\n_l-\frac13\delta_{i,l}\right)\right.\nonumber\\&&\left.+\left(\n_j\n_k-\frac13\delta_{j,k}\right)\left(\n_i\n_l+\frac13\delta_{i,l}\right)\right)\nonumber\\
&=&\frac{2s_*^3}{3}
\n_i\n_j\n_k\n_l.
\label{eq:m10}
\end{eqnarray}

\item $\mathrm{div}(\B_{11}[\nabla\bv]):$
\begin{eqnarray*}
\div\left[|\Q_*|^2(\A\cdot\Q_*)\Q_*\right]_i = \frac{2s_*^2}{3}\left(\Q_*\right)_{il}\left(\Q_*\right)_{jk}\partial_{kl}\bv_j=:\M^{11}_{i,j;k,l}\partial_{kl}\bv_j.
\end{eqnarray*}
Here
\begin{eqnarray}
\M^{11}_{i,j;k,l}&=&\frac{2s_*^4}{3}\left(\n_i\n_l-\frac13\delta_{i,l}\right)\left(\n_j\n_k-\frac13\delta_{j,k}\right)\nonumber\\
& = & 
\frac{2s_*^4}{3}\left(
\n_i\n_j\n_k\n_l 
- \frac13\delta_{ik}\n_l\n_j
\right)
.
\label{eq:m11}
\end{eqnarray}

\end{enumerate}

With the above, the fourth order tensor $\cM_\gamma$ can be decomposed as
\begin{equation}\label{uEqnM3}
\cM_{i,j;k,l} = \sum_{p=1,\ldots, 7, 10,11}\gamma_p\cM^p_{i,j;k,l}.
\end{equation}
From \eqref{eq:m1} to \eqref{eq:m11}, the $\M^p_{i,j;k,l}$'s are effectively
given by a linear combination of the following tensors:
\begin{equation}\label{MLinComb}
\delta_{kl}\delta_{ij},\,\,\,
\delta_{kl}\n_i\n_j,\,\,\,
\delta_{ij}\n_k\n_l,\,\,\,
\delta_{ik}\n_l\n_j\,\,\,(\text{or equivalently, $\delta_{il}\n_k\n_j$}),\,\,\,
\n_i\n_j\n_k\n_l.
\end{equation}

\section{Analysis of the anisotropic Stokes equation \eqref{eq:velo1s}--\eqref{eq:inc1s}}
Collecting the forms from \eqref{eq:A} and \eqref{eq:GIJH}, we write again
here the governing system as:
\begin{eqnarray}
\mathcal{L}_\gamma\bv + \nabla p &=& \f_\gamma(\bv),
\quad\text{for}\,\,\,|\bx| > a,
\label{eq:concise}\\
\quad\div\bv &=& 0,\quad\text{for}\,\,\,|\bx| > a,\\
\bv & =& \bv_b\,\,(=0),\quad \text{on $\partial\mathbf{B}_a$}\\
\bv  &=& \bv_*,\quad \text{at $|\bx|=\infty$.}\label{bdry.infty}
\end{eqnarray}
where 
\begin{eqnarray}
\mathcal{L}_\gamma\bv & := & 
-\Big(\Delta\bv + \M_\gamma:D^2\bv\Big),\\
\f_\gamma(\bu) & := & 
\div\cA_\gamma(x) 
+ \div\cC_\gamma(\bu, x) 
+ \cD_\gamma(\bu, x),\nonumber\\
&=& \gamma\left(\frac{{\bf F}(\hat{\br})}{r^3} +
\frac{{\bf G}(\hat{\br})}{r^4}
+ \frac{{\bf H}(\hat{\br}):(\bu-\bv_*)}{r^3}
+ \frac{{\bf I}(\hat{\br}):D\bu}{r^2}
+ \frac{{\bf J}(\hat{\br}):D^2\bu}{r}\right)\nonumber\\
& & + O\left(\frac{1}{r^5}\right).
\label{f:explicit}
\end{eqnarray}
Note that we have emphasized the dependence of $\cC_\gamma$ and
$\cD_\gamma$ on $\bu$. Given $\Q$, this dependence is in fact
linear in $\bu$ so that
\begin{equation}\label{f.lin}
\f_\gamma(\bu_1) - \f_\gamma(\bu_2) = \f_\gamma(\bu_1 - \bu_2)
\end{equation}
making the above a linear system 
which we designate as our anisotropic Stokes system.

We remark again that the above system describes the flow 
in the moving frame attached to the particle. 
The far-field velocity $\bv_*$ is prescribed. 
For passive particle which is the case in this paper, we take $\bv_b=0$ 
while for active particle, it is
in general some prescribed, non-constant function. Note also that
if $\bv_b = 0$ and $\bv_* =0$, then we have only the trivial solution
$\bv=0$.

Upon introducing
\begin{equation}\label{Tten2}
\T_\gamma({\bf w}, \pi) 
= \nabla {\bf w} + (\nabla{\bf w})^T
+ \cB_\gamma(\nabla\bv)
- \pi\id,
\end{equation}
we can write \eqref{eq:concise} as
\begin{equation}
-\div \T_\gamma(\bv, p) = \f_\gamma(\bv).
\label{eq:concise2}
\end{equation}
In terms of $\cM_{i,j;k,l}$, we have
$\big(\cB_\gamma(\nabla\bv)\big)_{ij} = \cM_{i,m;k,j}\partial_k\bv_m$
so that it is consistent with
$\big(\div\cB_\gamma(\nabla\bv)\big)_i = \cM_{i,j;k,l}\partial_{kl}\bv_j$.

The main purpose of the next few sections is to prove the existence
and uniqueness of a solution for \eqref{eq:concise}--\eqref{bdry.infty} in 
a suitable function space when $|\gamma|\ll 1$.

\subsection{Computation of the Green's Function for $\mathcal L$.}
Given a function or vector field $f$ defined on $\R^3$, its Fourier transform 
and inverse are given by
\[
\ft{f}(\bxi) = \int_{\mathbb{R}^n}f(\br)e^{-\ci\br\cdot\bxi}\,d\br,
\quad
f(\br) = (2\pi)^{-n}\int_{\mathbb{R}^n}\ft{f}(\bxi)e^{\ci\br\cdot\bxi}\,d\bxi,
\]
where $\bxi=(\xi_1, \xi_2, \xi_3)^T$.
Based on \eqref{uEqnM1} and \eqref{uEqnM2}, we introduce the 
$3\times 3$ matrix $[\M_\gamma](\xi)$:
\begin{equation}
[\M_\gamma]_{ij}(\xi) = \M_{i,j;k,l}\xi_k\xi_l.
\end{equation}
Then taking the Fourier transform of \eqref{uEqnM} gives
\begin{equation}
{|\bxi|}^2\ft{\bv} +[\M_\gamma](\xi)\ft{\bv}
+\ci\xi\ft{p}=\ft{\f_\gamma}
\end{equation}
where $\f_\gamma$ is the right hand side of \eqref{uEqnM}.
Using the incompressibility condition $\div\bv = 0$ written in Fourier mode
$\langle \ft{\bv}, \bxi\rangle = 0$, we have,
\begin{eqnarray*}
\langle [\M_\gamma](\xi)\ft{\bv}, \xi\rangle
+\ci|\bxi|^2\ft{p} 
= \langle \ft{\f_\gamma},\bxi\rangle,
\quad\text{so that}\quad
\ft{p}  =  
\frac{-\big\langle [\M_\gamma](\xi)\ft{\bv}, \xi\big\rangle + 
\langle \ft{\f_\gamma},\bxi\rangle}
{\ci|\bxi|^2}.
\end{eqnarray*}
Hence, 
\begin{eqnarray*}
{|\bxi|}^2\ft{\bv} + [\M_\gamma](\xi)\ft{\bv}
+\frac{-\big\langle[\M_\gamma](\xi)\ft{\bv}, \bxi\big\rangle\xi
+\langle\ft{\f_\gamma},\bxi\rangle\bxi}{|\bxi|^2}
=\ft{\f_\gamma},
\end{eqnarray*}
i.e.
\begin{eqnarray*}
{|\bxi|}^2\left(
\ft{\bv} + [\M_\gamma](\hat{\xi})\ft{\bv}
-
\big\langle[\M_\gamma](\hat{\xi})\ft{\bv}, \hat{\bxi}\big\rangle
\hat{\bxi}
\right)
&=&\ft{\f_\gamma}
-\langle\ft{\f_\gamma},\hat{\bxi}\rangle\hat{\bxi},
\,\,\,
\hat{\bxi} = \frac{\bxi}{|\bxi|}, 
\end{eqnarray*}
or equivalently,
\begin{eqnarray*}
{|\bxi|}^2\left(
\ft{\bv} + [\M_\gamma](\hat{\xi})\ft{\bv}
- (\hat{\xi}\otimes\hat{\xi})[\M_\gamma](\hat{\xi})\ft{\bv}
\right)
&=&(\I - \xi\otimes\xi)\ft{\f_\gamma}.
\end{eqnarray*}

To conclude, the solution $\bv$ and $p$ is given by
\begin{eqnarray}
\ft{\bv} &=& 
\frac{1}{|\bxi|^2}\Big(
\I + (\I-\hat{\bxi}\otimes\hat{\bxi})[\M_\gamma](\hat{\bxi})
\Big)^{-1}
(\I-\hat{\bxi}\otimes\hat{\bxi})\ft{\f_\gamma}
\end{eqnarray}
and
\begin{eqnarray}
\ft{p} = 
\frac{1}{\ci|\xi|}\left\langle
\left(
\I
- (\I-\hat{\bxi}\otimes\hat{\bxi})
\Big(
\I + (\I-\hat{\bxi}\otimes\hat{\bxi})[\M_\gamma](\hat{\bxi})
\Big)^{-T}
[\M_\gamma]^T(\hat{\xi})
\right)\hat{\xi}, \ft{\f_\gamma}
\right\rangle.
\end{eqnarray}
Note that the matrix inverse in the above is well-defined 
if $|\gamma| \ll 1$.

Let $\G(\bx)$ and $\h(x)$ be the inverse Fourier transforms respectively
of
\begin{equation}
\frac{1}{|\bxi|^2}\Big(
\I + (\I-\hat{\bxi}\otimes\hat{\bxi})[\M_\gamma](\hat{\bxi})
\Big)^{-1}
(\I-\hat{\bxi}\otimes\hat{\bxi})
\end{equation}
and
\begin{equation}
\frac{1}{\ci|\xi|}
\left(
\I
- (\I-\hat{\bxi}\otimes\hat{\bxi})
\Big(
\I + (\I-\hat{\bxi}\otimes\hat{\bxi})[\M_\gamma](\hat{\bxi})
\Big)^{-T}
[\M_\gamma]^T(\hat{\xi})
\right)\hat{\xi}.
\end{equation}
Then we have that $\G$ and $\h$ are homogeneous with degrees
$-1$ and $-2$ so that
\begin{equation}\label{G.Fct}
\G(\bx) = \frac{\G(\hat{\bx})}{|\bx|}
\,\,\,\text{and}\,\,\,
\h(\bx) = \frac{\h(\hat{\bx})}{|\bx|^2}.
\end{equation}
The above leads to the following properties: 
for $k=0, 1,2\ldots$, 
\begin{equation}\label{DkG.est}
(D_\bx^k\G)(\bx) = \frac{(D_\bx^k\G)(\hat{\bx})}{|\bx|^{k+1}},
\,\,\,
\Big|D^k_{\bx}\G(\bx)\Big| \lesssim \frac{1}{|\bx|^{k+1}},
\end{equation}
and
\begin{equation}\label{Dkh.est}
(D_\bx^k\h)(\bx) = \frac{(D_\bx^k\h)(\hat{\bx})}{|\bx|^{k+2}},
\,\,\,
\Big|D^k_{\bx}\h(\bx)\Big| \lesssim \frac{1}{|\bx|^{k+2}}.
\end{equation}

Using the $\G$ and $\h$ above, we are looking for a solution $\bv$ of 
\eqref{eq:concise} with suitable decay property at infinity such that the 
following representation holds for $\bv$:
\begin{eqnarray}
\bv(x)
&=& \int_\Omega \G(x-y)\f_\gamma(\bv(y))\,dy
+\int_{\partial\Omega}\big \langle 
\G(x-y), \T_\gamma(\bv, p)(y)\nu_y
\big \rangle\, d\sigma_y 
\nonumber\\
&&- \int_{\partial\Omega}\big \langle \T_\gamma(\G, \h)(x-y)\nu_y, \bv(y)\big \rangle
\,d\sigma_y + \bv_*.
\label{bv.rep}
\end{eqnarray}
Similar to \eqref{bdry.int.key}--\eqref{bdry.int.4}, we have
\begin{equation}
\left|\int_{\partial\Omega}\big\langle 
\G(x-y), \T_\gamma(\bv, p)(y)\nu_y
\big\rangle\, d\sigma_y\right| \lesssim \frac{1}{r},
\,\,\,\text{and}\,\,\,
\left|\int_{\partial\Omega}\big \langle \T_\gamma(\G, \h)(x-y)\nu_y, \bv(y)\big \rangle
\,d\sigma_y\right|
\lesssim \frac{1}{r^2},
\end{equation}
so that $\bv-\bv_* \approx r^{-1}$. 
More precise far-field asymptotics of $\bv$ will be given in 
Section \ref{asym.solns.model}. However, in preparation for the proof 
of existence of $\bv$, we will first analyze the
dominating term of $\f_\gamma(\bv)$ given by $\div\cA_\gamma(x)$.

\subsection{Far-field behavior of the inhomogeneous term}
\label{InhomogFarField}
We recall form \eqref{eq:A} for the dominating term of
$\f_\gamma(\bv)$:
\[
\div\cA_\gamma(\br) 
= \gamma\frac{{\bf F}(\hat{\br})}{r^3} + O\left(\frac{1}{r^4}\right).
\]
Hence, in order for the property $\bv-\bv_* \approx \frac{1}{r}$ to hold, 
from the representation formula \eqref{bv.rep}, it is necessary to have
$\displaystyle \int_\Omega \G(x-y)\div\cA_\gamma(y)\,dy 
\approx \frac{1}{r}$. By \eqref{bulk.no.log}, this is true only if
\begin{equation}\label{mean.zero.A}
\int_{{\mathbb S}^2}{\bf F}(\hat{\br})\,d\sigma = 0.
\end{equation}

To verify the above for our system, we write $\cA_\gamma$ explicitly as:
\begin{eqnarray}
\cA_\gamma
&=& \gamma_1\cA_1 + \gamma_2\cA_2 + \gamma_3\cA_3 + \gamma_9\cA_9
\nonumber\\
&=&
\gamma_1\Big(\Q_*(\bv_*\cdot\nabla\Q) - (\bv_*\cdot\nabla\Q)\Q_* \Big)
+\gamma_2\Big(\bv_*\cdot\nabla\Q\Big)\nonumber\\
&&+
\frac{\gamma_3}{2}\left(\Q_*(\bv_*\cdot\nabla\Q) + (\bv_*\cdot\nabla\Q)\Q_*\right)
+\gamma_9\Big((\bv_*\cdot\nabla\Q)\cdot\Q_*\Big)\Q_*.
\label{A.explicit}
\end{eqnarray}
By considering only the dominating term $\displaystyle \frac{1}{r}$ in the
expression \eqref{QABL} of $\Q$, we have
\begin{eqnarray*}
\Q - \Q_*& = & \left(-\frac{w\Q_*}{1+w}\right)\frac{1}{r} + O\left(\frac{1}{r^2}\right)\\
\bv_*\cdot\nabla \Q 
& = & \left(-\frac{w\Q_*}{1+w}\right)
\left\langle \bv_*, \nabla\frac{1}{r} \right\rangle 
+ O\left(\frac{1}{r^3}\right),\\
\Q_*(\bv_*\cdot\nabla \Q),\,\,\,
(\bv_*\cdot\nabla \Q)\Q_*
& = & \left(-\frac{w\Q_*^2}{1+w}\right)
\left\langle \bv_*, \nabla\frac{1}{r} \right\rangle
+ O\left(\frac{1}{r^3}\right), \\
\big((\bv_*\cdot\nabla \Q)\cdot \Q_*\big)\Q_*
& = & \left(-\frac{w|\Q_*|^2\Q_*}{1+w}\right)
\left\langle \bv_*, \nabla\frac{1}{r} \right\rangle 
+ O\left(\frac{1}{r^3}\right).
\end{eqnarray*}
Hence we can write
for some constant $3\times3$ matrix ${\bf M}_* = (m_{ij})$ that, 
\begin{equation}
\div\cA_\gamma = 
\text{div}\left(
\left\langle \bv_*, \nabla\frac{1}{r} \right\rangle
{\bf M}_*
\right)
+ O\left(\frac{1}{r^4}\right)
\label{divAform0}
\end{equation}
Note that the term in $\cA_\gamma$ multiplied by $\gamma_1$ completely 
vanishes.

We next claim that
\begin{equation}\label{mean.zero.int}
\int_{\mathbb S^2}
\left.\text{div}\left(
\left\langle \bv_*, \nabla\frac{1}{r} \right\rangle
{\bf M}_*
\right)\right|_{r=1}\,d\sigma = 0.
\end{equation}
To see this, let $\displaystyle g(r) = \frac{1}{r}$, and we compute for 
$i=1,2,3$,
\begin{equation}
\text{div}\left(
\left\langle \bv_*, \nabla\frac{1}{r} \right\rangle
{\bf M}_*
\right)_i
= \partial_{j}\big({\bv_*}_k \partial_{k}g(r)m_{ij}\big)
={\bv_*}_km_{ij}\partial^2_{kj}g(r),
\label{divAform}
\end{equation}
and
\begin{eqnarray}
\partial^2_{kj}g(r)
&=& g''(r)\partial_{k}r\partial_{j}r + g'(r)\partial^2_{kj}r
\nonumber\\
&=& g''(r)\frac{x_kx_j}{r^2} 
+ \frac{g'(r)}{r}\left(\delta_{kj} - \frac{x_kx_j}{r^2}\right)
\nonumber\\
&=& \frac{g'(r)}{r}\delta_{kj}
+ \left(g''(r)-\frac{g'(r)}{r}\right)\frac{x_kx_j}{r^2}
\nonumber\\
& = & -\frac{1}{r^3}\left(
\delta_{kj} - 3\frac{x_kx_j}{r^2}
\right)\label{DDgr}
\end{eqnarray}
where we have used the facts that
$\displaystyle g'(r) = -\frac{1}{r^2}$ and 
$\displaystyle g''(r) = \frac{2}{r^3}$.
As
\[
\int_{\mathbb{S}^2}3\hat{\br}\otimes\hat{\br}\,d\sigma = 4\pi\I,
\]
we have for all $k,j$ that
\[
\int_{\mathbb{S}^2}\partial^2_{kj}g(r)\,d\sigma = 0.
\]
Thus \eqref{mean.zero.int} holds and so does \eqref{mean.zero.A}.

\subsection{Existence of solution of \eqref{eq:concise} in Schauder spaces}
\label{Exist.Soln}
We recall the form of equation \eqref{eq:concise}, its inhomogeneous
term \eqref{f:explicit}, and the Green's function $\G$ \eqref{G.Fct}.
Estimate \eqref{DkG.est} for $\G$ plays an important role in our analysis.

We will show the existence and uniqueness of a solution by means of
the Banach Fixed Point Theorem in a suitable weighted Schauder space $\cS$.
For this purpose, we define
\begin{equation}
\cS =\big\{
\bu: \Omega\rightarrow\R^3,\,\,
\|\bu\|_{\cS} < \infty
\big\},
\end{equation}
where
\begin{equation}
\|\bu\|_{\cS} = 
\sup_{|\bx|\in\Omega}
\big\{
|\bx||\bu(\bx)-\bv_*|,\,\,\,
|\bx|^2|D\bu(\bx)|,\,\,\,
|\bx|^3|D^2\bu(\bx)|,\,\,\,
|\bx|^{3+\alpha}[D^2\bu]_\alpha(\bx)
\big\}.
\end{equation}
In the above, $0 < \alpha < 1$, and
\begin{equation}
[D^2\bu]_\alpha(\bx) = \sup_{\by\in\Omega, |\bx-\by|\leq 1}
\frac{|D^2\bu(\bx) - D^2\bu(\by)|}{|\bx-\by|^\alpha}.
\end{equation}

Now, given a $\bu\in\cS$, let $\bv = \cT(\bu)$ be the solution of 
\eqref{eq:concise}--\eqref{bdry.infty} with $\f_\gamma = \f_\gamma(\bu)$.
We will find a fixed point $\bv$ of $\cT$: $\bv = \cT(\bv)$. 
To achieve this, we will show the following two properties of $\cT$:
\begin{enumerate}
\item	$\cT$ maps $\cS$ into $\cS$, in particular, there is a $C> 0$ 
such that for any $\bu\in\cS$,
\begin{equation}
\|\cT(\bu)\|_\cS \leq C\Big[
\|\bu\|_\cS + \|\bv_b\|_{C^{2,\alpha}(\partial\Omega)} + |\bv_*|
\Big].
\label{T.bd}
\end{equation}
\item	For $|\gamma|\ll 1$, 
there exists a $0 < c < 1$ such that for any $\bu_1,\bu_2\in\cS$, 
we have 
\begin{equation}
||\cT(\bu_1)-\cT(\bu_2)||_\cS \leq c||\bu_1-\bu_2||_{\cS}.
\label{T.contract}
\end{equation}
\end{enumerate}
The proof can be obtained via the following steps.\\

\noindent
{\bf (I) Well-posedness of $\cT$.}
General existence and uniqueness theory for $\f_\gamma\in L^p(\Omega)$ can be
found in \cite{galdiSimaderARMA}, also described in the encyclopedic reference
\cite[Chapter V]{galdi2011introduction}. 
But as $\f_\gamma$ has specific spatial decay 
property and we are looking for classical solutions, we find it convenient
to follow the classical approach outlined in \cite[Chapter 3]{lady} using
the theory of single and double layer potentials, and Schauder estimates. 
This theory is also outlined in \cite{MRS1999}. The recent survey \cite{Russo2010} covers the existence and uniqueness of solutions for the Stokes equation in Schauder spaces.
For the convenience of the reader, we now outline the approach that demonstrates
existence.

Given $\bu\in{\mathcal S}$, we write $\f(\br) := \f_\gamma(\bu(\br))$ as
\begin{eqnarray}
{\f}(\br) &=& {\bf h}_1(\br) + {\bf h}_2(\br),\label{h1+h2}\\
{\bf h}_1(\br) & := & \div\cA_\gamma(\br) = \gamma\frac{\bF(\hat{\br})}{r^3}
+ O\left(\frac{1}{r^4}\right)
\,\,\,\left(
\text{recall}\,\,\,\int_{\mathbb{S}^2}\bF(\hat{\br})\,d\sigma =0
\right),\\
{\bf h}_2(\br) & := & \div\cC_\gamma(\bu, x) + \cD_\gamma(\bu, x).
\label{h2}
\end{eqnarray}
From \eqref{eq:A} and \eqref{eq:GIJH}, we have
$\left|{\bf h}_1\right| \lesssim \frac{1}{r^3}$ and 
$\left|{\bf h}_2\right| \lesssim \frac{1}{r^4}$ for $r \gg 1$.
The solution $\bv$ is found by the following steps.

\begin{enumerate}
\item
We extend $\f$ smoothly from $\Omega$ to $\tilde{\f}$ defined on the 
whole space $\R^3$. Given the decay property of $\tilde{\f}$ inherited from 
$\bf f$, we can set
\[
\tilde{\bv}_1(x) = \int_{\R^3}\G(x-y)\tilde{\f}(y)\,dy
\,\,\,\text{and}\,\,\,
\tilde{p}_1(x) = \int_{\R^3}\h(x-y)\cdot\tilde{\f}(y)\,dy.
\]
By \eqref{h1+h2}--\eqref{h2} and \eqref{bulk.no.log}, we have that 
$\tilde{\bv}_1$ and $p_1$ solves the following equation
\[
-\Delta\tilde{\bv}_1 + \nabla\tilde{p}_1 = \tilde{\f}(\br),\,\,\,
\div\tilde{\bv}_1 = 0,\,\,\,\text{with}\,\,\,
|\tilde{\bv}_1(\br)| \lesssim \frac{1}{r}\,\,\,\text{and}\,\,\,
|\tilde{p}_1| \lesssim \frac{1}{r^2},
\]
on the whole of $\R^3$.
\item
We find $\bv_2$ and $p_2$ on $\Omega$ that solve the
following equation:
\begin{eqnarray*}
-\Delta\bv_2 + \nabla p_2 &=& 0,\,\,\,\text{in}\,\,\,\Omega\\
\div\bv_2 & = & 0,\,\,\,\text{in}\,\,\,\Omega\\
\bv_2  & = & \bv_b - \tilde{\bv}_1 - \bv_*,\,\,\,\text{on}\,\,\,\partial\Omega,\\
\bv_2  & = & 0,\,\,\,\text{at}\,\,\,|\bx|=\infty.
\end{eqnarray*}
Following \cite[Theorem 6.1, eqn. (6.3)]{MRS1999}, 
we can have that in the class
$\bv_2 = o(r)$, the solution exists, is unique, and it satisfies
$|\bv_2(\br)|\lesssim\frac{1}{r}$.

\item
Define 
\[
\bv = \tilde{\bv}_1 + \bv_2 + \bv_*,\,\,\,
p = \tilde{p}_1 + p_2.
\]
Then $\bv$ and $p$ solve \eqref{eq:concise} on 
$\Omega$. 
\end{enumerate}

The spatial decay property of $\bv$ shows that it can have the following 
representation,
\begin{multline}
\bv(x)
= \bv_* + \int_\Omega \G(x-y)\f(y)\,dy\\
+ \int_{\partial\Omega}
\big \langle \G(x-y), \T_\gamma(\bv, p)(y)\nu_y\big\rangle\,d\sigma_y 
- \int_{\partial\Omega}\big \langle \T_\gamma(\G, \h)(x-y)\nu_y, \bv(y)\big \rangle
\,d\sigma_y.
\label{bulk.int.soln}
\end{multline}

Next we use \eqref{bulk.int.soln} to estimates $\bv$.
For estimates near the boundary, we can invoke
\cite[Chapter 3, Theorem 5]{lady} to deduce that
\begin{equation}\label{Bdry.Holder}
||\bv||_{C^{2,\alpha}(\overline{\Omega}_{a+1})},\,\,
||p||_{C^{1,\alpha}(\overline{\Omega}_{a+1})}
\lesssim
||\f||_{C^{0,\alpha}(\Omega)}+||\bv_b||_{C^{2,\alpha}(\partial\Omega)} + |\bv_*|
\end{equation}
where $\Omega_{b} = \Omega\cap\{|\bx| \leq b\}$.
In particular, on the boundary, we have
\begin{equation}\label{Bdry.C0}
||\bv||_{C^0(\partial\Omega)},\,\,
||D\bv||_{C^0(\partial\Omega)},\,\,
||p||_{C^0(\partial\Omega)}
\lesssim
||\f||_{C^{0,\alpha}(\Omega)}+||\bv_b||_{C^{2,\alpha}(\partial\Omega)} + |\bv_*|.
\end{equation}
Furthermore, from the form of $\bf f$ \eqref{f:explicit}, we have
for some constant $C_1, C_2$ that
\begin{equation}
||{\bf f}||_{C^{0,\alpha}(\Omega)}
\leq
\gamma \big[C_1\|\bu\|_{\cS} + C_2\big].
\end{equation}

Now we proceed with the following two proofs.\\

\noindent
{\bf (II) Proof of \eqref{T.bd}}.
With boundary estimates given by \eqref{Bdry.Holder},
we will just concentrate here on {\em interior weighted} estimates, 
i.e., for $\bx$ such that $|\bx| > a + 1$. Again, we will utilize the
representation \eqref{bulk.int.soln}.

For the boundary integrals in \eqref{bulk.int.soln}, we have
\begin{align*}
\left|\int_{\partial\Omega}\big \langle \G(x-y), \T_\gamma(\bv, p)(y)\nu_y\big \rangle\, d\sigma_y\right|
&\lesssim \frac{||D\bv||_{C^0(\partial\Omega)}+||p||_{C^0(\partial\Omega)}}{r}
\end{align*}
and more generally, for $k=1, 2, \ldots$,
\begin{align*}
\left|D^k_\bx\int_{\partial\Omega}\big \langle \G(x-y), 
\T_\gamma(\bv, p)(y)\nu_y\big \rangle\, 
d\sigma_y\right|
 & \lesssim
\int_{\partial\Omega}
\left|\big \langle 
D^k_{\bx}\G(x-y),
\T_\gamma(\bv, p)(y)\nu_y
\big\rangle\right|
\,d\sigma_y \\
& \lesssim \frac{||D\bv||_{C^0(\partial\Omega)}+||p||_{C^0(\partial\Omega)}}{r^{1+k}}.
\end{align*}
Furthermore, for $\bh\in\R^n$ with $|\bh|\leq 1$,
\begin{eqnarray*}
&&\frac{|\bx|^{3+\alpha}}{|h|^\alpha}
\left|\int_{\partial\Omega}\big \langle 
D_\bx^2\G(x+h-y) - D_\bx^2\G(x-y),
\T_\gamma(\bv, p)(y)\nu_y
\big \rangle\, d\sigma_y\right|\\
& \lesssim &
\big(\|\bv\|_{C^1(\partial\Omega)}+\|p\|_{C^0(\partial\Omega)}\big)
\frac{|\bx|^{3+\alpha}}{|h|^\alpha}
\left|\int_{\partial\Omega}
\frac{\big||x+h-y|^3-|x-y|^3\big|}{|x+h-y|^3|x-y|^3}
\, d\sigma_y\right|\\
& \lesssim & 
\big(\|\bv\|_{C^1(\partial\Omega)}+\|p\|_{C^0(\partial\Omega)}\big)
\left(\frac{|\bx|^{3+\alpha}}{|h|^\alpha}\right)
\frac{|\bx|^2|h|}{|\bx|^6}\\
& \lesssim & 
\big(\|\bv\|_{C^1(\partial\Omega)}+\|p\|_{C^0(\partial\Omega)}\big)
\left(\frac{|h|^{1-\alpha}}{|\bx|^{1-\alpha}}\right)\\
& \lesssim & 
\big(\|\bv\|_{C^1(\partial\Omega)}+\|p\|_{C^0(\partial\Omega)}\big).
\end{eqnarray*}
Hence, 
\begin{eqnarray*}
\left\|
\int_{\partial\Omega}\big \langle 
\G(x-y),
\T_\gamma(\bv, p)(y)\nu_y
\big \rangle
\, d\sigma_y
\right\|_{\cal S}
&\lesssim &
\|\bv\|_{C^1(\partial\Omega)}+\|p\|_{C^0(\partial\Omega)}\\
& \lesssim &
||{\bf f}||_{C^{0,\alpha}(\Omega)}+||\bv_b||_{C^{2,\alpha}(\partial\Omega)}\\
& \lesssim &
|\gamma|\Big[\|\bu\|_{\cS} 
+ 1\Big] +||\bv_b||_{C^{2,\alpha}(\partial\Omega)} + |\bv_*|.
\end{eqnarray*}

Using the same technique, we can similarly have,
\begin{eqnarray*}
\left\||\bx|\left|
\int_{\partial\Omega}\big 
\langle \T_\gamma(\G, \h)(x-y)\nu_y, \bv(y)\big \rangle
d\sigma_y
\right|
\right\|_{\cal S}
&\lesssim &
\|\bv\|_{C^1(\partial\Omega)}+\|p\|_{C^0(\partial\Omega)}\\
& \lesssim &
||{\bf f}||_{C^{0,\alpha}(\Omega)}+||\bv_b||_{C^{2,\alpha}(\partial\Omega)}\\
& \lesssim &
|\gamma|\Big[\|\bu\|_{\cS} 
+ 1\Big] + ||\bv_b||_{C^{2,\alpha}(\partial\Omega)} + |\bv_*|.
\end{eqnarray*}

For the bulk integrals in \eqref{bulk.int.soln}, the proof is very similar to 
the classical estimates for Newtonian potential -- see \cite[Chapter 4]{GT}. 
We will produce the proof for our weighted space $\cal S$ in Appendix 
\ref{PotEst}.\\

\noindent
{\bf (III) Proof of \eqref{T.contract}}.
With \eqref{T.bd}, estimate \eqref{T.contract} follows easily, due to the
linearity of the equation. Given $\bu_1, \bu_2\in\cS$, we have
solutions $\bv_1 = \cT(\bu_1)$ and $\bv_2 = \cT(\bu_2)$. 
Upon subtracting the corresponding
equations, we deduce that
\begin{eqnarray}
{\mathcal L}_\gamma(\bv_1-\bv_2) + \nabla (p_1-p_2)
= \f_\gamma(\bu_1) - \f_\gamma(\bu_2),
\end{eqnarray}
with $\bv_1-\bv_2=\bf 0$ on $\partial\Omega$ and at $|\bx|=\infty$.
Hence
\[
||\bv_1-\bv_2||_\cS
\lesssim
||\f_\gamma(\bu_1) - \f_\gamma(\bu_2)||_{C^{0,\alpha}(\Omega)}
\lesssim
|\gamma| \|\bu_1-\bu_2\|_\cS,
\]
where we have used \eqref{f.lin}.
Now choosing $|\gamma|\ll 1$ gives the result.

\section{Properties of anisotropic Stokes flows}
\label{asym.solns.model}

Here we make more precise the far-field behavior of the solution $\bv$.
As an application of our analysis, we will analyze the symmetry property
of the solution and give a decomposition formula for the Stokes drag.
They will be validated and illustrated by numerical simulations.

We start from representation \eqref{bulk.int.soln} for the solution $\bv$:
\begin{eqnarray}
\bv(x)
&=& \int_\Omega \G(x-y)\f_\gamma(\bv(y))\,dy
+ \int_{\partial\Omega}\big \langle 
\G(x-y), \T_\gamma(\bv, p)(y)\nu_y
\big \rangle\, d\sigma_y\nonumber\\
&&- \int_{\partial\Omega}\big \langle
\T_\gamma(\G, \h)(x-y)\nu_y, \bv(y)\big \rangle
\,d\sigma_y + \bv_* \nonumber\\
& =: & \sS_1(\bx) + \sS_2(\bx) + \sS_3(\bx) + \bv_*
\end{eqnarray}
where
\begin{eqnarray*}
\sS_1(\bx) 
&:=& \int_\Omega \G(x-y)\f_\gamma(\bv(y))\,dy\nonumber\\
&=& \int_\Omega \G(x-y)\div\cA_\gamma(y)\,dy
 + \int_\Omega\G(x-y)\big(\div\cC_\gamma(\bv) + \cD_\gamma(\bv)\big)(y)\,dy
\nonumber\\
& = & \int_\Omega \G(x-y)\div\cA_\gamma(y)\,dy
+ \G(x)\left(\int_\Omega\big(\div\cC_\gamma(\bv) + \cD_\gamma(\bv)\big)(y)\,dy
\right)\\
& & +
\int_\Omega(\G(x-y)-\G(x))\big(\div\cC_\gamma(\bv) + \cD_\gamma(\bv)\big)(y)
\,dy,\\
\sS_2(\bx) 
&:=&
\int_{\partial\Omega}\big \langle \G(x-y), \T_\gamma(\bv, p)(y)\nu_y
\big \rangle\, d\sigma_y\nonumber\\
& = & 
\left\langle\G(x), \int_{\partial\Omega}\T_\gamma(\bv, p)(y)\nu_y\,d\sigma_y
\right\rangle
+ 
\int_{\partial\Omega}\Big\langle \G(x-y)-\G(x),
\T_\gamma(\bv, p)(y)\nu_y
\Big\rangle
\,d\sigma_y,\\
\sS_3(\bx) &:=&
- \int_{\partial\Omega}\big \langle
\T_\gamma(\G, \h)(x-y)\nu_y, \bv(y)\big \rangle
\,d\sigma_y.
\end{eqnarray*}
We re-arrange the above terms in the following way,
\begin{eqnarray}
\bv-\bv_*&=&
\left\langle\G(x), \int_{\partial\Omega}\T_\gamma(\bv, p)(y)\nu_y\,d\sigma_y
\right\rangle
\label{1.0}\\
&&+
\int_\Omega \G(x-y)\div\cA_\gamma(y)\,dy
+ \G(x)\left(\int_\Omega\big(\div\cC_\gamma(\bv) + \cD_\gamma(\bv)\big)(y)\,dy
\right)\label{1.1}\\
&&
+\int_{\partial\Omega}\Big\langle
\G(x-y)-\G(x), \T_\gamma(\bv, p)(y)\nu_y\Big\rangle
\,d\sigma_y\label{2.0.1}\\
&&
- \int_{\partial\Omega}\big \langle
\T_\gamma(\G, \h)(x-y)\nu_y, \bv(y)\big \rangle
\,d\sigma_y\label{2.0.2}
\\
&&+\int_\Omega(\G(x-y)-\G(x))\big(\div\cC_\gamma(\bv) + \cD_\gamma(\bv)\big)(y)
\,dy\label{2.1}.
\end{eqnarray}
Note that for $|\bx|\gg 1$, we have
\begin{eqnarray*}
\text{\eqref{1.0}} & = & \frac{O(1)}{|\bx|},\\
\text{\eqref{1.1}} & = & \frac{O(\gamma)}{|\bx|},\\
\text{\eqref{2.0.1} and \eqref{2.0.2}} & = & \frac{O(1)}{|\bx|^2},\\
\text{\eqref{2.1}} & = & \frac{O(\gamma)\log|x|}{|\bx|^2}
\quad\text{(by \eqref{logx2})}
\end{eqnarray*}
so that up to $O(\gamma)$, \eqref{1.0} and \eqref{1.1} are the dominant terms in the expression for the anisotropic flow. 

We next describe asymptotically \eqref{1.0} and \eqref{1.1} for $\gamma\ll 1$.

\subsection{Precise asymptotics: deviation from isotropic Stokes flow}
The purpose of this section is to reveal more clearly the difference between 
$\bv$ and the classical Stokes flow $\bv_0$ which is set to satisfy
\begin{eqnarray}
-\Delta\bv_0 + \nabla p_0 &=& 0,
\quad\text{for}\,\,\,|\bx| > a,
\label{eq:concise0}\\
\div(\bv_0) & = & 0,\quad\text{for}\,\,\,|\bx| > a,\nonumber\\
\bv_0 &=& 0,\,\,\,\text{on $|x|=a$,}\nonumber\\
\bv_0 &=& \bv_*,\,\,\,\text{at $|x|=\infty$.}\nonumber
\end{eqnarray}
Define $\varphi_\gamma :=\bv-\bv_0$. Then we have
\begin{eqnarray*}
-\Delta \varphi_\gamma + \nabla (p-p_0) 
&=& -\M_\gamma:D^2\bv + \f_\gamma(\bv)\\
&=& 
-\M_\gamma:D^2\bv_0 + \f_\gamma(\bv_0) + 
\big(\f_\gamma(\bv) - \f_\gamma(\bv_0)
-\M_\gamma:(D^2\bv - D^2\bv_0)\big).
\end{eqnarray*}
Note that $|\M_\gamma|,\,\,|{\bf f}_\gamma(\bv_0)| \lesssim |\gamma|$
and 
$|\f_\gamma(\bv)-\f_\gamma(\bv_0)|=|\f_\gamma(\bv-\bv_0)| 
\lesssim |\gamma|\|\bv-\bv_0\|_{\cS}$.
Now let $\lbar{\varphi}_\gamma$ solve
\begin{eqnarray}
-\Delta \lbar{\varphi}_\gamma + \nabla (\lbar{p}_\gamma)
&=& -\M_\gamma:D^2\bv_0 + \f_\gamma(\bv_0),\label{bar.varphi.g}\\
\text{div}(\lbar{\varphi}_\gamma) &=& 0, \quad
\text{for $|\bx| > a$,}\\
\lbar{\varphi}_\gamma &=& 0,\quad\text{on $|\bx|=0$ and at $|\bx|=\infty$.}
\end{eqnarray}
Then the same approach in deriving estimates for $\bv$ gives
\begin{equation}
\|\lbar{\varphi}_\gamma\|_{\cS} \lesssim O(\gamma),
\quad\text{and}\quad
\|\lbar{\varphi}_\gamma - \varphi_\gamma\|_{\cS} \lesssim O(\gamma^2).
\end{equation}
Hence we have
\begin{equation}\label{v-v0-varphig}
\bv = \bv_0 + \lbar{\varphi}_\gamma + O(\gamma^2).
\end{equation}

Finally, using the Green's function $\E$ \eqref{GreenE} for the classical 
Stokes flow, we have the following representation of $\lbar{\varphi}_\gamma$,
\begin{eqnarray}
\lbar{\varphi}_\gamma (x) 
&=& \int_\Omega \E(x-y)\big(-\M_\gamma:D^2\bv_0 + \f_\gamma(\bv_0(y))\big)\,dy
+ \int_{\partial\Omega}\big \langle 
\E(x-y), \T(\lbar{\varphi}_\gamma, \lbar{p}_\gamma)(y)\nu_y
\big \rangle\, d\sigma_y.\nonumber\\
&=& \int_\Omega \E(x-y)\big(-\M_\gamma:D^2\bv_0 + \f_\gamma(\bv_0(y))\big)\,dy
+ \E(x)\int_{\partial\Omega} \T(\lbar{\varphi}_\gamma, \lbar{p}_\gamma)(y)
\nu_y\,d\sigma_y
\nonumber\\
& & + 
\int_{\partial\Omega}\big \langle 
\E(x-y)-\E(x),
\T(\lbar{\varphi}_\gamma, \lbar{p}_\gamma)(y)\nu_y
\big \rangle\, d\sigma_y
\end{eqnarray}
where $\T$ is the Stokes stress tensor \eqref{Tten1}.
Note that the last term in the above decays as $\frac{1}{|x|^2}$.
Hence we have,
\begin{eqnarray}
\lbar{\varphi}_\gamma(x) & = &
\cI_\gamma(x) + \E(x)\cJ_\gamma + O\left(\frac{1}{|x|^2}\right),\\
\text{where}\quad
\cI_\gamma(x) &:=&
\int_\Omega \E(x-y)\big(-\M_\gamma:D^2\bv_0 + \f_\gamma(\bv_0(y))\big)\,dy,\\
\cJ_\gamma & := &
\int_{\partial\Omega} \T(\lbar{\varphi}_\gamma, \lbar{p}_\gamma)(y)
\nu_y d\sigma_y.
\end{eqnarray}

The rest of this section will describe more explicitly the 
bulk and boundary integrals $\cI_\gamma$ and $\cJ_\gamma$.

\subsubsection{Analysis of $\cI_\gamma$}
For this purpose, recalling from \eqref{StokesFormGenVec}
that $\bv_0 = \E_S\bv_*$, we then have
\begin{eqnarray}
\cI_\gamma(x) & = &
\int_\Omega \E(x-y)\big[-\M_\gamma:D^2\bv_0 + \f_\gamma(\bv_0(y))\big]\,dy
\nonumber\\
&=&\int_\Omega \E(x-y)\Big[6\pi a\M_\gamma:D^2\E(y)\bv_*
-\frac{a^3}{4}\M_\gamma:D^2\F(y)\bv_*
\nonumber\\
& & 
\quad\quad\quad\quad\quad+ \div\cA_\gamma(y) + \div\cC_\gamma(y) + \cD_\gamma(y)\Big]\,dy
\nonumber\\
&=&
\int_\Omega 
\E(x-y)\Big[6\pi a\M_\gamma:D^2\E(y)\bv_* + \div\cA_\gamma(y) \Big]\,dy\\
& & + \E(x)\int_\Omega \left[
-\frac{a^3}{4}\M_\gamma:D^2\F(y)\bv_*
+ \div\cC_\gamma(y) + \cD_\gamma(y)\right]\,dy
\nonumber\\
& & + \int_\Omega (\E(x-y)-\E(x))\left[
-\frac{a^3}{4}\M_\gamma:D^2\F(y)\bv_*
+ \div\cC_\gamma(y) + \cD_\gamma(y)\right]\,dy
\nonumber\\
&=:& I_1 + I_2 + I_3.
\end{eqnarray}
Note that $I_2$ decays as $\E(x) \sim \frac{1}{|x|}$ because the integrand is
integrable. By \eqref{logx2}, the integral $I_3$ decays as 
$\frac{\log|x|}{|x|^2}$.

For $I_1$, we will show that $D^2\E(x)$ satisfies the mean zero condition. 
(Such a condition is already verified for $\div\cA_\gamma$ in 
Section \ref{InhomogFarField}.) By \eqref{bulk.no.log}, we can then conclude that
$I_1(x)$ decays as $\frac{1}{|x|}$.
To this end, by \eqref{DDE}, we have
\begin{eqnarray}
\int_{\mathbb{S}^2}
\partial_{kl}\E_{ij}(x)\,d\sigma
& = & 
4\pi\left(-\delta_{ij}\delta_{kl} + \delta_{ik}\delta_{jl} 
+ \delta_{jk}\delta_{il}\right)
-4\pi\left(
-\delta_{ij}\delta_{kl} +\delta_{ik}\delta_{jl} +\delta_{jk}\delta_{il}
\right)\nonumber\\
& & -4\pi\left(
\delta_{il}\delta_{jk} +\delta_{jl}\delta_{ik} +\delta_{kl}\delta_{ij}
\right)
+ 15\int_{\mathbb{S}^2}\hat{x}_i\hat{x}_j\hat{x}_k\hat{x}_l\,d\sigma
\nonumber\\
& = & -4\pi\left(
\delta_{il}\delta_{jk} +\delta_{jl}\delta_{ik} +\delta_{kl}\delta_{ij}
\right)
+ 15\int_{\mathbb{S}^2}\hat{x}_i\hat{x}_j\hat{x}_k\hat{x}_l\,d\sigma.
\end{eqnarray}
Using spherical coordinates and considering symmetry, we check
\begin{enumerate}
\item	$i\neq j, k\neq l$, $(k,l) = (i,j)$ or $(j,i)$
\begin{eqnarray*}
\int_{\mathbb{S}^2}
\partial_{kl}\E_{ij}(x)\,d\sigma 
& =&  -4\pi + 15\int_{\mathbb{S}^2}\hat{x}_1^2\hat{x}_2^2\,d\sigma\\
&=& -4\pi 
+ 15\int_0^\pi\int_0^{2\pi} (\sin\phi\cos\theta)^2(\sin\phi\sin\theta)^2
\sin\phi\,d\phi\,d\theta\\
& = & -4\pi 
+ 15\int_0^\pi \sin^5\phi \,d\phi
\int_0^{2\pi} \cos^2\theta\sin^2\theta \,d\theta = 0
\end{eqnarray*}

\item	$i\neq j, k\neq l$, $(k,l) \neq (i,j), (j,i)$
\begin{eqnarray*}
\int_{\mathbb{S}^2} \partial_{kl}\E_{ij}(x)\,d\sigma = 0
\end{eqnarray*}

\item	$i\neq j, k=l=i$ or $k=l=j$
\begin{eqnarray*}
\int_{\mathbb{S}^2} \partial_{kl}\E_{ij}(x)\,d\sigma 
= + 15\int_{\mathbb{S}^2}\hat{x}_1^3\hat{x}_2\,d\sigma = 0.
\end{eqnarray*}

\item	$i\neq j, k=l,\neq i, j$
\begin{eqnarray*}
\int_{\mathbb{S}^2} \partial_{kl}\E_{ij}(x)\,d\sigma = 0
\end{eqnarray*}

\item	$i=j, k\neq l$
\begin{eqnarray*}
\int_{\mathbb{S}^2} \partial_{kl}\E_{ij}(x)\,d\sigma = 0
\end{eqnarray*}

\item	$i=j, k=l,\neq i,j$
\begin{eqnarray*}
\int_{\mathbb{S}^2} \partial_{kl}\E_{ij}(x)\,d\sigma 
= -4\pi + 15\int_{\mathbb{S}^2}\hat{x}_1^2\hat{x}_2^2\,d\sigma = 0,
\quad\text{as in (a)}.
\end{eqnarray*}

\item	$i=j=k=l$
\begin{eqnarray*}
\int_{\mathbb{S}^2} \partial_{kl}\E_{ij}(x)\,d\sigma 
= -12\pi + 15\int_{\mathbb{S}^2}\hat{x}_3^4\,d\sigma 
= -12\pi + 15\int_0^\pi\int_0^{2\pi}\cos^4\phi\sin\phi\,d\theta\,d\phi= 0.
\end{eqnarray*}
\end{enumerate}

Hence we can invoke \eqref{bulk.no.log} to conclude that
\begin{eqnarray}
\cI_\gamma(x)
& = &
\frac{1}{|x|}
\int_{\mathbb{S}^2}
\HH(\hat{x},\omega)
\Big[6\pi a\M_\gamma:D^2\E(\omega)\bv_* 
+ \div\cA_\gamma(\omega)
\Big]\,d\sigma_\omega
\label{cB1}\\
& & + \E(x)\int_\Omega \left[
-\frac{a^3}{4}\M_\gamma:D^2\F(y)\bv_*
+ \div\cC_\gamma(y) + \cD_\gamma(y)\right]\,dy
\label{cB2}\\
&& + O\left(\frac{\log|x|}{|x|^2}\right)\nonumber
\end{eqnarray}

\subsubsection{Analysis of $\cJ_\gamma$}
To obtain an explicit formula for the boundary stress $\cJ_\gamma$
associated with 
$\lbar{\varphi}_\gamma$, we multiply \eqref{bar.varphi.g} 
by a test function given by $\tbv_0-\tbv_*$ where 
$\tilde{\bv}_0$ solves the following equation:
\begin{eqnarray*}
-\Delta\tbv_0 + \nabla \tilde{p}_0 &=& 0,
\quad\text{for}\,\,\,|\bx| > a,\\
\div(\tbv_0) & = & 0,\quad\text{for}\,\,\,|\bx| > a,\nonumber\\
\tbv_0 &=& 0,\,\,\,\text{on $|x|=a$,}\nonumber\\
\tbv_0 &=& \tbv_*,\,\,\,\text{at $|x|=\infty$.}\nonumber
\end{eqnarray*}
with an arbitrary $\tbv_*$. Using integration 
by parts -- see \cite[p.53 (10), (11)]{lady}, we obtain
\begin{multline*}
-\int_{\partial\Omega}
\Big\langle\T(\lbar{\varphi}_\gamma, \lbar{p}_\gamma)\nu_y,
(\tbv_0-\tbv_*)\Big\rangle\,d\sigma_y + 
\int_{\partial\Omega}
\Big\langle\T(\tbv_0, \tilde{p}_0)\nu_y,
\lbar{\varphi}_\gamma\Big\rangle\,d\sigma_y \\
=\int_\Omega
\Big\langle -\M_\gamma:D^2\bv_0 + \f_\gamma(\bv_0(y)),\tbv_0-\tbv_*\Big\rangle
\,dy.
\end{multline*}
The above is justified by the decay estimates for 
$\lbar{\varphi}_\gamma$ and $\tbv_0-\tbv_*$:
\begin{eqnarray*}
\left|\int_{\{|\bx|=R\}}
\Big\langle\T(\lbar{\varphi}_\gamma, \lbar{p}_\gamma)\nu_y,
(\tbv_0-\tbv_*)\Big\rangle\,d\sigma_y\right|
\lesssim \frac{1}{R^2}\frac{1}{R}R^2 = \frac{1}{R}\\
\left|\int_{\{|\bx|=R\}}
\Big\langle\T(\tbv_0, \tilde{p}_0)\nu_y,
\lbar{\varphi}_\gamma\Big\rangle\,d\sigma_y\right|
\lesssim \frac{1}{R^2}\frac{1}{R}R^2 = \frac{1}{R}
\\
\int_{\{|\bx|>R\}} 
\Big|-\M_\gamma:D^2\bv_0 + \f(\bv_0(y))\Big|\Big|\tbv_0-\bv_*\Big|\,dy
\lesssim \int_R^\infty\frac{1}{r^3}\frac{1}{r}r^2\,dr \lesssim \frac{1}{R}
\end{eqnarray*}
which all vanish as $R\longrightarrow0$.
Hence we have
\begin{eqnarray*}
&&\int_{\partial\Omega}
\Big\langle\T(\lbar{\varphi}_\gamma, \lbar{p}_\gamma)\nu_y,
(\tbv_0-\tbv_*)\Big\rangle\,d\sigma_y\nonumber\\
&=&
\int_{\partial\Omega}
\Big\langle\T(\tbv_0, \tilde{p}_0)\nu_y,
\lbar{\varphi}_\gamma\Big\rangle\,d\sigma_y
-\int_\Omega
\Big\langle -\M_\gamma:D^2\bv_0 + \f(\bv_0(y)),\tbv_0-\tbv_*\Big\rangle
\,dy.
\end{eqnarray*}
As $\tbv_0=\lbar{\varphi}_\gamma=0$ on $\partial\Omega$, the
following holds for any $\tbv_*$,
\begin{eqnarray*}
\int_{\partial\Omega}
\Big\langle\T(\lbar{\varphi}_\gamma, \lbar{p}_\gamma)\nu_y, \tbv_*\Big\rangle\,d\sigma_y
&=& \int_\Omega
\Big\langle -\M_\gamma:D^2\bv_0 + \f_\gamma(\bv_0(y)),\tbv_0-\tbv_*\Big\rangle
\,dy\nonumber\\
&=& \int_\Omega
\Big\langle -\M_\gamma:D^2\E_S(y)\bv_* + \f_\gamma(\E_S(y)\bv_*),(\E_S-\I)\tbv_*\Big\rangle
\,dy.
\end{eqnarray*}
Hence we have
\begin{eqnarray}
\cJ_\gamma & = & \int_{\partial\Omega}
\T(\lbar{\varphi}_\gamma, \lbar{p}_\gamma)\nu_y \,d\sigma_y\nonumber\\
&=& \int_\Omega
(\E_S-\I) \Big[ -\M_\gamma:D^2\E_S(y)\bv_* 
+ \f_\gamma\big(\E_S(y)\bv_*\big) \Big] \,dy\nonumber\\
&=& \int_\Omega
(\E_S-\I) \Big[ -\M_\gamma:D^2\E_S(y)\bv_* \nonumber\\
&&\hspace{60pt}
+ \div\cA_\gamma(y) 
+ \div\cC_\gamma(\E_S(y)\bv_*, y) 
+ \cD_\gamma(\E_S(y)\bv_*, y)\Big]\,dy.
\label{cJ}
\end{eqnarray}

\subsection{Symmetry properties of solution and drag force}
Here we investigate the symmetries of $\bv$ with respect to our 
far-field data $\Q_*$ and $\bv_*$. Note that the interacting 
potential $R$ in \eqref{eq:diss} is chosen to be frame 
indifferent. Hence our solution is naturally invariant with respect to 
orthogonal transformations. For the reader's convenience, 
we outline the derivation.

Let $\bbB=(b_{ij})$ and $\bbA=(a_{ij})=\bbB^T (=\bbB^{-1})$ be an orthogonal 
matrix and its transpose (and inverse). Upon introducing
\begin{eqnarray*}
y=\bbB x,\,\,\,\bv(x,t) = \bbA\widetilde{\bv}(y,t),\,\,\,
p(x,t) = \wt{p}(y,t),\,\,\,
\partial_i = \partial_{x_i},\,\,\,
\widetilde{\partial}_i = \partial_{y_i}
\,\,\,(\text{so that}\,\,\,\partial_i = a_{ij}\wt{\partial}_j),
\end{eqnarray*}
we have,
\begin{eqnarray*}
&\bv\cdot\nabla = \widetilde{\bv}\cdot\widetilde{\nabla},\,\,\,
\bv\cdot\nabla\bv = \bbA\left(\wt{\bv}\cdot\wt{\nabla}\wt{\bv}\right),\,\,\,
\div\bv = \wt{\div}\wt{\bv},\,\,\,
\triangle\bv = \bbA\,\wt{\triangle}\wt{\bv},\\
&\nabla p = \bbA\wt{\nabla}\wt{p},\,\,\,
\bv_t = \bbA\widetilde{\bv}_t.
\end{eqnarray*}
The above shows that the Navier-Stokes equation
\[
\bv_t + \bv\cdot\nabla \bv + \nabla p = \triangle \bv,\,\,\,\div\bv = 0,
\]
is equivalent to
\[
\wt{\bv}_t + \wt{\bv}\cdot\nabla \wt{\bv} + \wt{\nabla} \wt{p} = 
\triangle \wt{\bv},
\,\,\,\wt{\div}\wt{\bv} = 0.
\]

To take into account of our model \eqref{eq:qten}--\eqref{eq:inc}, 
we first note that for any order two tensors (or tensor fields) 
$\sR$ and $\sS$, suppose $\wt{\sR}$ and $\wt{\sS}$ are such that
$\sR(x) = \bbA\wt{\sR}(y)\bbB$ and $\sS(x) = \bbA\wt{\sS}(y)\bbB$, 
then it holds that
\[
\sR\cdot\sS = \wt{\R}\cdot\wt{\sS},\,\,\,
\text{and}\,\,\,
\div\sR = \bbA\,\wt{\div}\wt{\sR}.
\]
As application, consider $\A=\nabla\bv + \nabla\bv^T$,
$\W=\nabla\bv - \nabla\bv^T$, 
$\wt{\A} = \wt{\nabla}\wt{\bv} + \wt{\nabla}\wt{\bv}^T$ and 
$\wt{\W} = \wt{\nabla}\wt{\bv} - \wt{\nabla}\wt{\bv}^T$, then we have
\[
\A = \bbA \wt{\A} \bbB,\,\,\,\text{and}\,\,\,
\W = \bbA \wt{\W} \bbB.
\]
For $\Q$, let $\wt{\Q}$ be such that
$\Q(x) = \bbA\widetilde{\Q}(y)\bbB$.
In particular, if
$\Q = s_*(\n\otimes\n - \frac{\I}{3})$, then
$\wt{\Q} = s_*((\bbB\n)\otimes(\bbB\n) - \frac{\I}{3})$.
In general, the following hold,
\[
\mathcal{E}_\ldg(\nabla\Q, \Q)
=\mathcal{E}_\ldg(\wt{\nabla}\wt{\Q}, \wt{\Q}),\,\,\,
\mathcal{F}_\ldg(\Q)=\mathcal{F}_\ldg(\wt{\Q}),
\]
and 
\[
\mathring{\Q}(x) 
= \Q_t + \bv\cdot\nabla\Q + \Q\W - \W\Q
= \bbA\mathring{\widetilde{\Q}}(y)\bbB.
\]
Hence for
$\T_{\text{SV}}^\text{v}=\T_{\text{SV}}^\text{v}(\mathring{\Q}, \Q, \A)$
and
$\T^{\text{el}}=\T^{\text{el}}(\nabla\Q)$,
upon defining
$\wt{\T}_{\text{SV}}^\text{v}
=\T_{\text{SV}}^\text{v}(\mathring{\wt{\Q}}, \wt{\Q}, \wt{\A})$ and 
$\wt{\T}^{\text{el}}
=\T^{\text{el}}(\wt{\nabla}\wt{\Q})$,
from \eqref{eq:visc} and \eqref{eq:elas}, it holds that
\[
\T_{\text{SV}}^\text{v} = \bbA\wt{\T}_{\text{SV}}^\text{v}\bbB,
\,\,\,\,\,\,
\T^\mathrm{el} = \bbA\wt{\T}^\mathrm{el}\bbB,
\,\,\,\,\,\,
\div\T_{\text{SV}}^\text{v} = \bbA\,\wt{\div}\wt{\T}_{\text{SV}}^\text{v},
\,\,\,\,\,\,
\div\T^\text{el} = \bbA\,\wt{\div}\wt{\T}^\text{el}.
\]
Furthermore, let 
$\wt{\mathcal{E}}_\ldg(\wt{\nabla}\wt{\Q}, \wt{\Q})
=\mathcal{E}_\ldg(\wt{\nabla}\wt{\Q}, \wt{\Q})$,
from \eqref{nonlinf} and \eqref{LDlaplace}, then
\begin{eqnarray*}
\frac{\partial\mathcal{E}_\ldg}{\partial\Q} =
\bbA\frac{\partial\wt{\mathcal{E}}_\ldg}{\partial\wt{\Q}}\bbB
\,\,\,\,\,\,\text{and}\,\,\,\,\,\,
\div\left[\frac{\partial\mathcal{E}_\ldg}{\partial\nabla\Q}\right]
= 
\bbA\,\div\left[\frac{\partial\wt{\mathcal{E}}_\ldg}{\partial\wt{\nabla}\wt{\Q}}\right]\bbB.
\end{eqnarray*}
We then conclude that system \eqref{eq:qten}--\eqref{eq:inc} is 
invariant under orthogonal transformations. 
In particular, we say that
a solution $\bv$, $\Q$ is invariant under an orthogonal transformation 
$\bbB$ if
\begin{eqnarray}
\wt{\bv}(x) = \bv(x), && \text{or equivalently,}\,\,\,\,\,\,
\bv(\bbB x) = \bbB\bv(x),\\
\wt{\Q}(x) = \Q(x), && \text{or equivalently,}\,\,\,\,\,\,
\Q(\bbB x) = \bbB \Q(x)\bbA.
\end{eqnarray}

Next we investigate the symmetry property of the drag force
$\cF_\gamma$ on the particle which according to \eqref{eq:velo} is defined as
\begin{equation}
\cF_\gamma := \int_{\partial\bB}\T\nu\,d\sigma_x,
\,\,\,\,\,\,
\text{where}\,\,\,\T = -p\I + \T_{\text{SV}}^\text{v} + \T^\text{el}.
\end{equation}
By the invariance property just demonstrated, we have
\begin{eqnarray}
\cF_\gamma 
= \int_{\partial\bB}\bbA\wt{\T}\bbB\nu_x\,d\tilde{\sigma}_y
= \bbA\int_{\partial\wt{\bB}}\wt{\T}\wt{\nu}_y\,d\tilde{\sigma}_y
=:\bbA\wt{\cF}_\gamma,\,\,\,\,\,\,
\text{where
$\wt{\T} = -\wt{p}\I + \wt{\T}_{\text{SV}}^\text{v} + \wt{\T}^\text{el}$.}
\end{eqnarray}
Note that in the above, we have used the fact that the area element does not 
change under orthogonal transformation. 
Now if the solution is invariant under $\bbB$,
then $\wt{\T}(y) = \T(y)$. Hence
\[
\wt{\cF}_\gamma 
= \int_{\partial\wt{\bB}}\wt{\T}\wt{\nu}_y\,d\tilde{\sigma}_y
= \int_{\partial\bB}\T\nu_x\,d\sigma_x
= \cF_\gamma.
\]
This leads to 
\begin{equation}
\cF_\gamma = \bbA \cF_\gamma,\,\,\,\,\,\,\text{or equivalently,}\,\,\,\,\,\,
\cF_\gamma = \bbB \cF_\gamma.
\end{equation}

Next we make use of the above to analyze the symmetries of the drag force
$\cJ_\gamma$.
For convenience, we use $\cF_\gamma(\bvs,\Q_*)$ to denote the dependence 
of $\cF_\gamma$ on $\bvs$ and $\Q_*$.
Note that for our linear system \eqref{eq:concise}--\eqref{f:explicit},
with $\bv_b=0$, given $\Q_*$, its solution $\bv$ and hence $\cF_\gamma$ is 
linear in the far-field velocity field $\bv_*$. To take advantage of this, 
for concreteness, we let $\n=e_3$ in $\Q_*$ \eqref{QABL2}. 
Now we decompose $\bv_*$ as,
\begin{equation}
\bvs = \bvs^\parallel + \bvs^\perp,\quad
\text{where}\,\,\,\bvs^\parallel \parallel e_3,\,\,\,
\bvs^\perp \perp e_3,\,\,\,
\text{i.e. ${\bv_*}^\perp$ lies in the $xy$-plane.}
\end{equation}
Then the solution $\bv$ and the drag force $\cF_\gamma$ can be decomposed as
\[
\bv = \bv^\parallel + \bv^\perp,\quad
\cF_\gamma(\bvs,\Q_*) 
= \cF_\gamma(\bvs^\parallel,\Q_*) + \cF_\gamma(\bvs^\perp,\Q_*)
\]
where $\bv^\parallel, \bv^\perp$ solve 
\eqref{eq:concise}--\eqref{f:explicit} with
$\bvs^\parallel, \bvs^\perp$ as the far-field velocity field. 
Furthermore, we have
\begin{equation}
\R\cF_\gamma(\bvs, \Q_*)
=\cF_\gamma(\R\bvs, \R\Q_*\R^{-1})
\end{equation}
for any orthogonal transformation $\R$ of $\R^3$.
The following arguments reveal the symmetry properties of $\cF_\gamma$ with
respect to $\bvs^\parallel$, $\bvs^\perp$.
\begin{itemize}
\item	Since $\Q_*$ and ${\bv_*}^\parallel$ are 
invariant under any orthogonal transformation of $\R^3$ that leaves 
$e_3$ fixed, we can infer that $\bv^\parallel$ satisfies the same property. 
Hence the drag force $\cF_\gamma(\bvs^\parallel, \Q_*)$ must be parallel
to $e_3$. By linearity, we have for some constant $\gamma_\parallel$ that
\begin{equation}
\cF_\gamma(\bvs^\parallel,\Q_*) 
= \gamma_\parallel\langle \bvs^\parallel, e_3\rangle e_3
= \gamma_\parallel\big(e_3\otimes e_3)\bvs.
\end{equation}

\item	For $\bvs^\perp$, we have in general that
\begin{equation}
\R'\cF_\gamma(\bvs^\perp,\Q_*)
= \cF_\gamma(\R'\bvs^\perp,\R'\Q_*{\R'}^{-1})
= \cF_\gamma(\bvs^\perp,\Q_*)
\end{equation}
for any reflection $\R'$ of $\R^3$ that leaves the plane spanned by 
$e_3$ and $\bvs^\perp$ fixed. Hence we must have
$\cF_\gamma(\bvs^\perp,\Q_*)\in\text{Span}\{e_3,\bvs^\perp\}$. By linearity,
we have for some constant $\gamma_\perp$ and vector ${\bf C}$
from the $xy$-plane that
\[
\cF_\gamma(\bvs^\perp,\Q_*)
= \gamma_\perp\bvs^\perp + \langle {\bf C}, \bvs^\perp\rangle e_3.
\]
Now for any orthogonal transformation $\R''$ of $\R^3$ that leaves 
$e_3$ fixed, using
$$\R''\cF_\gamma(\bvs^\perp,\Q_*) 
= \cF_\gamma(\R''\bvs^\perp,\R''\Q_*{\R''}^{-1}),$$
we have
\[
\gamma_\perp\R''\bvs^\perp + \langle {\bf C}, \bvs^\perp\rangle e_3
=
\gamma_\perp\R''\bvs^\perp + \langle {\bf C}, \R''\bvs^\perp\rangle e_3.
\]
Hence $\langle {\bf C}, \bvs^\perp\rangle = 
\langle {\bf C}, \R''\bvs^\perp\rangle$ for any $\R''$.
Thus ${\bf C}$ must be zero.
We then conclude that
\begin{equation}
\cF_\gamma(\bvs^\perp,\Q_*) = \gamma_\perp\bvs^\perp 
=\gamma_\perp\big(\I - e_3\otimes e_3\big)\bvs.
\end{equation}
\end{itemize}

Combining the above, we finally have the following same formula as
\cite[(6.9)]{stark2001physics}:
\begin{equation}
\cF_\gamma(\bvs,\Q_*) 
=
\cF_\gamma(\bvs^\parallel,\Q_*)  + 
\cF_\gamma(\bvs^\perp,\Q_*)
=
\Big[\gamma_\parallel\big(e_3\otimes e_3\big)
+ \gamma_\perp\big(\I - e_3\otimes e_3\big)
\Big]\bvs.
\end{equation}
This is also consistent with the fact that if we replace $\n$ by $-\n$, 
$\Q_*$ remains unchanged and so does the the overall system.

Even though our system of equation is a reduced version of the original model, it does not seem easy to write down an asymptotic formula for the coefficients $\gamma_\perp$ and 
$\gamma_\parallel$. This is because of the presence of $\frac{1}{\text{Er}}$ in the stress tensor \eqref{eq:velnd} which cannot be easily computed in the limit of vanishing $\text{Er}$. See \cite{stark2001physics} for a discussion of some analysis and conjectures about the drag force.

\subsection{Analytical calculations}
We demonstrate here analytical calculations of \eqref{cB1}, \eqref{cB2}
and \eqref{cJ}. In principle, we can give analytical expressions for
all the terms as they only involve homogeneous functions
of negative integral degrees. More precisely, we have
\begin{eqnarray*}
\E(x)\,\,\text{\eqref{GreenE}}: && \frac{\delta_{ij}}{r}, \,\frac{x_ix_j}{r^3},\\
\F(x)\,\,\text{\eqref{GreenF}}: && \frac{\delta_{ij}}{r^3}, \,\frac{x_ix_j}{r^5}\\
\Q\,\,\text{\eqref{QABL}}: && \frac{\delta_{ij}}{r},\,\frac{\delta_{ij}}{r^3},\,\frac{x_ix_j}{r^5}
\end{eqnarray*}
and $\bv = \E_S\bv_*$ with $\E_S = \I - 6\pi a\E(x) + \frac{a^3}{4}\F(x)$
\eqref{StokesFormGenVec}.
Furthermore, the terms $\cC_\gamma$ and $\cD_\gamma$ involve
multiplications between the following matrices
\[
\A,\,\,\W,\,\,\Q_*,\,\,\Q,\,\,\,\bv_*\cdot\nabla\Q,\,\, \bv\cdot\nabla\Q.
\]

For convenience, we introduce the following conventions:
\begin{eqnarray}
[A, B, C, \ldots] & = & 
\text{arbitrary linear combinations and products between $A$, $B$, $C$,}
\nonumber\\
&&\text{and their powers;}
\label{notation.lin.mult}
\\
\{A, B, C, \ldots\} & = &
\text{linear combinations between $A$, $B$, $C$.}
\label{notation.lin.comb}
\end{eqnarray}
From Appendix \ref{MD2EF}, we have
\begin{eqnarray}
\M_\gamma:D^2\E(x)
&\in& \frac{1}{r^3}\big[1,\langle \n, \hx\rangle\big]
\Big\{\I,\,\,\,
\n\otimes \n,\,\,\,
\n\otimes\hx,\,\,\,
\hx\otimes \n,\,\,\,
\hx\otimes\hx
\Big\}\\
\M_\gamma:D^2\F(x)
&\in& \frac{1}{r^5}\big[1,\langle \n, \hx\rangle\big]
\Big\{\I,\,\,\,
\n\otimes \n,\,\,\,
\n\otimes\hx,\,\,\,
\hx\otimes \n,\,\,\,
\hx\otimes\hx
\Big\}
\end{eqnarray}
so that
\begin{eqnarray}
\M_\gamma:D^2\E(x)\bv_*
&\in& \frac{1}{r^3}
\Big[1,\,\,\,
\langle \n, \hx\rangle,\,\,\,
\langle \n, \bv_*\rangle,\,\,\,
\langle \hx, \bv_*\rangle
\Big]
\Big\{\bv_*,\,\,\,\n,\,\,\,\hx\Big\}\\
\M_\gamma:D^2\F(x)\bv_*
&\in& \frac{1}{r^5}
\Big[1,\,\,\,
\langle \n, \hx\rangle,\,\,\,
\langle \n, \bv_*\rangle,\,\,\,
\langle \hx, \bv_*\rangle
\Big]
\Big\{\bv_*,\,\,\,\n,\,\,\,\hx\Big\}.
\end{eqnarray}
From Appendix \ref{CAgammaForm} and \ref{CDgammaForm}, we have
\begin{eqnarray}
\div\cA_\gamma
&\in&\left\{\frac{1}{r^3},\frac{1}{r^4}\right\}
\Big[1,\,\,\,
\langle \hx,\n\rangle,\,\,\,
\langle \hx,\bv_*\rangle,\,\,\,
\langle \n,\bv_*\rangle
\Big]\Big\{
\hx,\,\,\n,\,\,\bv_*
\Big\},\\
\cC_\gamma
&\in&
\left[\frac{1}{r^3},\ldots,\frac{1}{r^6}\right]
\Big[1,\,\,
\langle \hx,\n\rangle,\,\,
\langle \hx,\bv_*\rangle,\,\,
\langle \n,\bv_*\rangle
\Big]\times\\
&&
\Big\{\I,\,\,\,\n\otimes \n,\,\,\,
\n\otimes\bv_* ,\,\,\,\bv_*\otimes \n,\,\,\,
\n\otimes \hx,\,\,\,\hx\otimes \n,\,\,\,
\bv_*\otimes \hx,\,\,\,\hx\otimes \bv_*,\,\,\,
\hx\otimes \hx
\Big\},\nonumber\\
\cD_\gamma
&\in&
\left[\frac{1}{r^4},\ldots,\frac{1}{r^9}\right]
\Big[1,\,\,\,
\langle \hx,\n\rangle,\,\,\,
\langle \hx,\bv_*\rangle,\,\,\,
\langle \n,\bv_*\rangle
\Big]\Big\{
\hx,\,\,\n,\,\,\bv_*
\Big\}.
\end{eqnarray}
Note that all the terms in $\M_\gamma:D^2\E$, 
$\M_\gamma:D^2\F$ and $\div\cA_\gamma$ have already been computed explicitly
in the Appendix \ref{MD2EF} and \ref{CAgammaForm}. From Appendix 
\ref{CDgammaForm}, all the coefficients appearing in 
$\cC_\gamma$ and $\cD_\gamma$ are amenable for symbolic computations. 

Using the above, we can give the following representations of the bulk and
boundary terms appearing in \eqref{cB1}, \eqref{cB2} and \eqref{cJ}.
\begin{itemize}
\item
$\displaystyle \int_\Omega 
\E(x-y)\Big[6\pi a\M_\gamma:D^2\E(y)\bv_* + \div\cA_\gamma(y) \Big]\,dy$.
By \eqref{bulk.no.log}, we have
\begin{eqnarray*}
&&\int_\Omega
\E(x-y)\Big[6\pi a\M_\gamma:D^2\E(y)\bv_* + \div\cA_\gamma(y) \Big]\,dy\\
&=&
\frac{1}{|x|}\int_{\mathbb{S}^2}\HH(\hat{x},\hy)
\Big[6\pi a\M_\gamma:D^2\E(\hy)\bv_* 
+ \div\cA_\gamma(\hy)\Big]\,d\sigma(\hy) 
+ O\left(\frac{1}{|x|^2}\right)\\
& = & 
\frac{1}{|x|}\int_{\mathbb{S}^2}\HH(\hat{x},\hy)
\Big[A_1(\hy)\hy + A_2(\hy) \n + A_3(\hy)\bv_*\Big]
\,d\sigma(\hy) 
+ O\left(\frac{1}{|x|^2}\right),
\end{eqnarray*}
where $\HH$ given by \eqref{HH.formG}, is recorded here,
\begin{equation}\label{HH.formE}
\HH(\hat{x},\hy)
= \lim_{\epsilon\to0}\left[
\int_{\epsilon}^\infty \frac{\E(\hat{x}-r\hy)}{r}\,dr
+\E(\hat{x})\ln(\epsilon)
\right],
\end{equation}
and the $A_i(\hy)$ are appropriate polynomial functions of
$\langle \hy,\n\rangle, \,\langle \hy,\bv_*\rangle,\,\langle \n,\bv_*\rangle$.

\item $\displaystyle
\int_\Omega \left[
-\frac{a^3}{4}\M_\gamma:D^2\F(y)\bv_*
+ \div\cC_\gamma(y) + \cD_\gamma(y)\right]\,dy
$.
The terms $D^2\F$ and $\cD_\gamma$ decay at least $r^{-4}$ and hence are integrable.
Using their forms, we have
\begin{eqnarray}
&&\int_\Omega \left[
-\frac{a^3}{4}\M_\gamma:D^2\F(y)\bv_* + \cD_\gamma(y)\right]\,dy
= 
\int_{\mathbb{S}^2}
\Big[A_4(\hy)\hy + A_5(\hy) \n + A_6(\hy)\bv_*\Big]\,d\sigma(\hy).\nonumber\\
\end{eqnarray}
For $\div\cC_\gamma = O\left(\frac{1}{r^3}\right)$, it can be conveniently
represented using Divergence Theorem:
\begin{equation}
\int_\Omega \div\cC_\gamma(y)\,dy
= \int_{\mathbb{S}^2}\Big\langle \cC_\gamma(\hy), \hy\Big\rangle\,d\sigma(\hy)
= \int_{\mathbb{S}^2}
\Big[A_7(\hy)\hy + A_8(\hy) \n + A_9(\hy)\bv_*\Big]\,d\sigma(\hy).
\end{equation}

\item $\displaystyle
\int_\Omega
(\E_S-\I) \Big[ -\M_\gamma:D^2\E_S(y)\bv_* 
+ \div\cA_\gamma(y)
+ \div\cC_\gamma(\E_S(y)\bv_*, y)
+ \cD_\gamma(\E_S(y)\bv_*, y)\Big]\,dy$.
As $\E_S-\I$ decays at least as $r^{-1}$, the following integral is integrable:
\begin{eqnarray}
&&\int_\Omega
(\E_S-\I) \Big[ -\M_\gamma:D^2\E_S(y)\bv_*
+ \div\cA_\gamma(\bv_*, y)
+ \cD_\gamma(\E_S(y)\bv_*, y)\Big]\,dy\\
& = & 
\int_{\mathbb{S}^2}
\Big[A_{10}(\hy)\hy + A_{11}(\hy) \n + A_{12}(\hy)\bv_*\Big]\,d\sigma(\hy),
\end{eqnarray}
where we recall the form of $\E_S$ from \eqref{StokesFormGenVec}.
The remaining term with $\div\cC_\gamma$ can also be dealt with using
Divergence Theorem:
\begin{eqnarray}
&&\int_\Omega (\E_S-\I) \div\cC_\gamma(\E_S(y)\bv_*, y) \,dy\nonumber\\
& = & 
\int_{\mathbb{S}^2}\Big\langle (\E_S(\hy)-\I), 
\cC_\gamma(\hy)\hy\Big\rangle\,d\sigma(\hy)
- \int_\Omega\Big\langle\nabla\E_S(y), \cC_\gamma(y)\Big\rangle\,dy
\nonumber\\
& = & 
\int_{\mathbb{S}^2}
\Big[A_{13}(\hy)\hy + A_{14}(\hy) \n + A_{15}(\hy)\bv_*\Big]\,d\sigma(\hy).
\end{eqnarray}
\end{itemize}

Before moving on to the next section, we consider one ``simplistic model'' in which 
$\Q$ is taken to be uniform in space, i.e. it equals its end state $\Q_*$. Such an approximation was in fact used in some works, see for example \cite{stark2001stokes,kos2018elementary}.
In this case, all the terms $\cA_\gamma$, $\cC_\gamma$ and $\cD_\gamma$
vanish as they involve either $\nabla\Q$ or $\Q-\Q_*$. Then we have
\begin{eqnarray}
\cI_\gamma(x)
&=&
\frac{1}{|x|}
\int_{\mathbb{S}^2}
\HH(\hat{x},\omega)
\Big[6\pi a\M_\gamma:D^2\E(\omega)\bv_*\Big]\,d\sigma_\omega\nonumber
\\&&+ \E(x)\int_\Omega \left[ -\frac{a^3}{4}\M_\gamma:D^2\F(y)\bv_* \right]\,dy,\\
\cJ_\gamma 
&=& \int_\Omega
(\E_S-\I) \Big[ -\M_\gamma:D^2\E_S(y)\bv_*
\Big] \,dy.
\end{eqnarray}
As demonstrated numerically in Section \ref{NumQQstar}, we see that the 
actual velocity flow $\bv$ does depend on the overall structure of $\Q$, 
not just its end state $\Q_*$. 

\section{Numerical simulations}\label{Sec:NumSim}
Here we provide numerical simulations to illustrate our analysis. 
The simulations are performed using a commercial finite elements software 
package COMSOL \cite{comsol}. 
For validation, we used this package to compute the classical Stokes flow. The results are benchmarked
against \emph{analytical solutions} in a \emph{finite domain}, more
precisely in the annulus $a \leq r \leq R$. For details, we 
refer to Appendix \ref{StokesBdDom}.

In the following, we record our numerical results for the anisotropic
Stokes system \eqref{eq:qten2s}-\eqref{eq:inc1s}. Some remarks are in order. 
\begin{enumerate}
\item
For simplicity, we assume that only $\gamma_1$ and 
$\gamma_2$ are nonzero. They are fixed to be $\gamma_1 =1$ and $\gamma_2 = 0.9$, except in 
Section \ref{NumGamma} where we allow 
them to vary. 

\item
Our analytical results show that 
$\bv = \bv_* + O\left(\frac{1}{r}\right)$. To better illustrate this, we will
plot \emph{rescaled} versions of components of $\bv$. More precisely, let $V = |\bv_*|$. If $\bv_*$ is along $e_i$, then
we will plot the following quantity for $a < r < R:$
\begin{equation}\label{rescale-r}
r\left(\frac{\bv_i}{V}-1\right)\,\,\,\text{and}\,\,\,
r\left(\frac{\bv_j}{V}\right)\,\,\,\text{for $j\neq i$}
\end{equation}
along various two-dimensional planes.

\item
All our numerical results shows that the rescaled quantity \eqref{rescale-r}
remains bounded in magnitude. Inevitably, due to boundary effect, the 
numerical solution is consistent with the true solution only for 
$a < r \ll R$. This is made more precise in Appendix \ref{StokesBdDom}.
\end{enumerate}

\subsection{$\bv_* = e_1$, $\n = e_3$}\label{NumQQstar}
In this section, we choose $\bv_*$ and $\n$ to be non-parallel to each other.
Besides plotting various (rescaled) components of $\bv_i$'s, we aim to 
illustrate the clear differences between the solution $\bv$ of 
our anisotropic Stokes system when $\Q$ is set to be $\Q_*$ and given 
by \eqref{QABL}. 
\begin{center}
(a)\includegraphics [height=2in]
{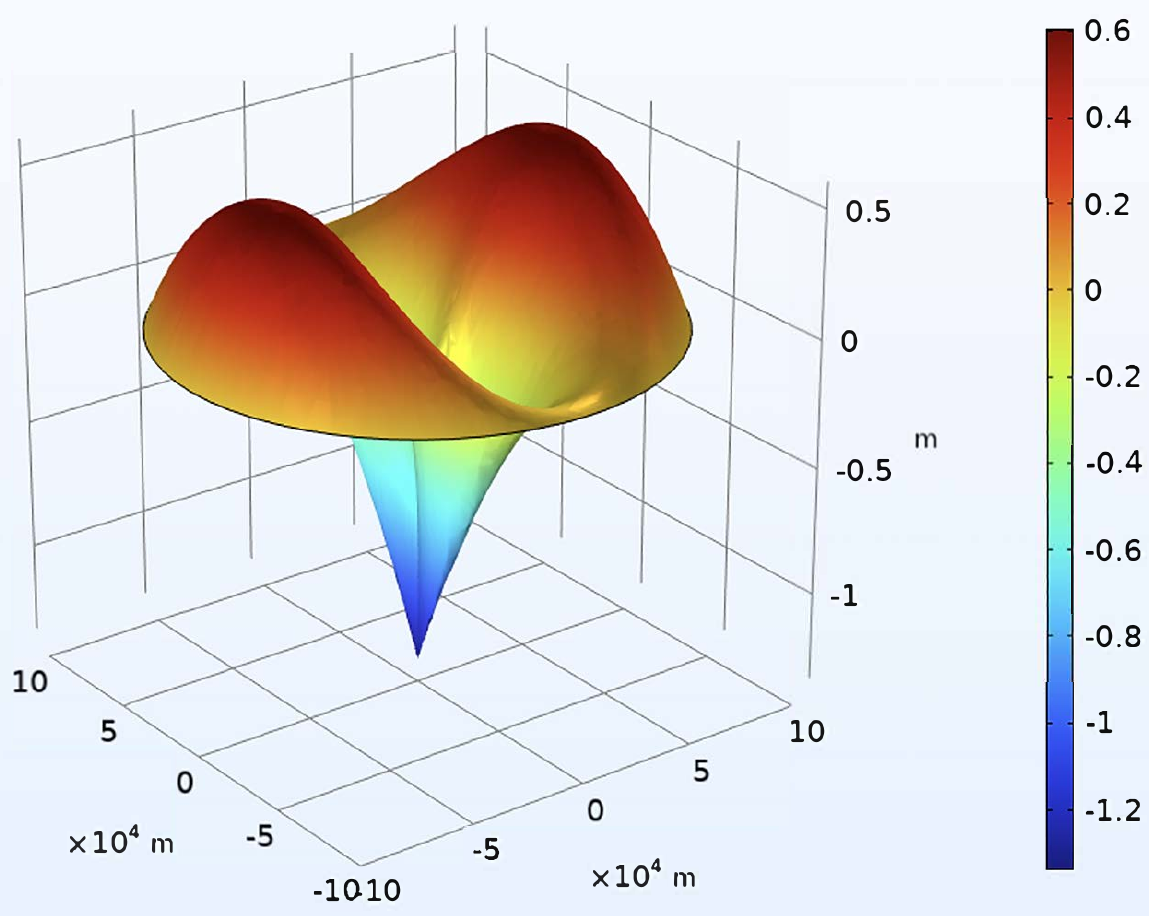}
\,\,\,
(b)\includegraphics [height=2in]
{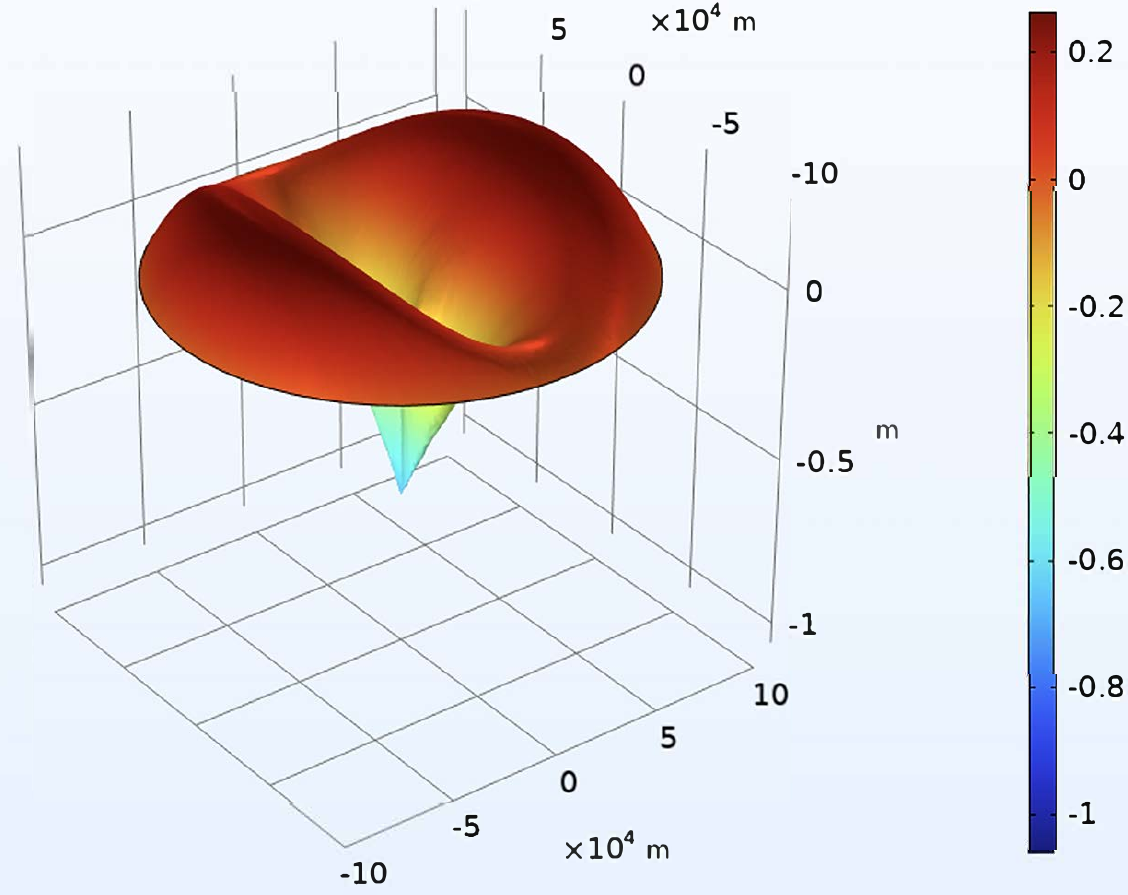}
\\
Figure 1. 3D, rescaled plot of $\bv_1$ in the $yz$-plane.\\
(a) $\Q$ is set to be $\Q_*$;
(b) $\Q$ is given by \eqref{QABL}.
\end{center}

\begin{center}
(a)\includegraphics [height=2in]
{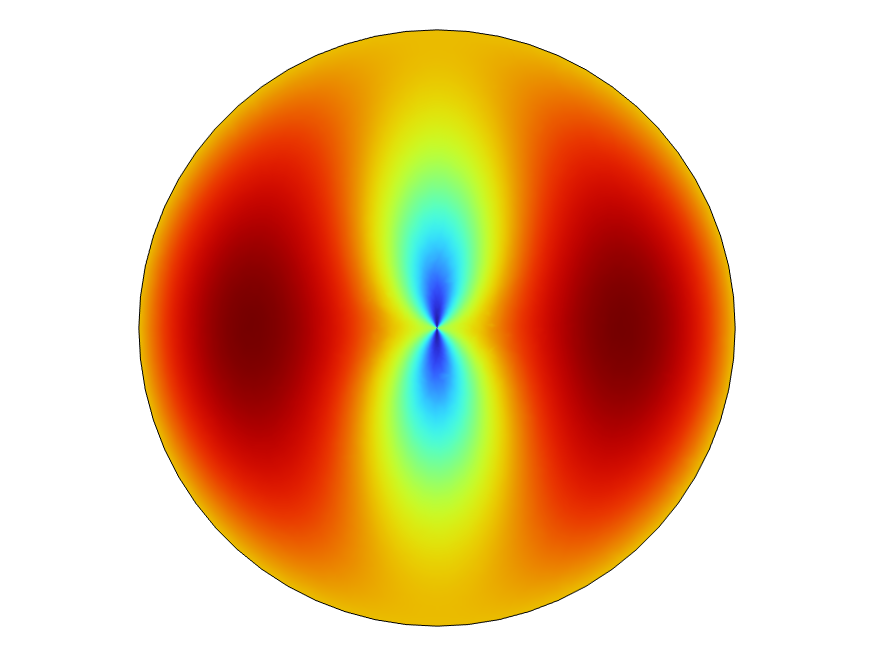}
\,\,\,
(b)\includegraphics [height=2in]
{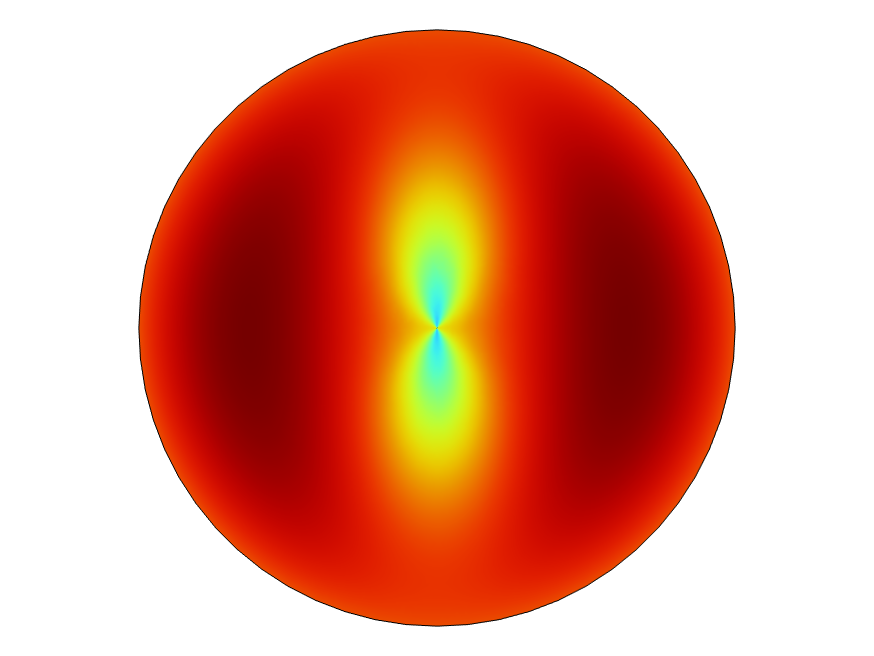}
\\
Figure 2. 2D, rescaled plot of $\bv_1$ in the $yz$-plane.\\
(a) $\Q$ is set to be $\Q_*$;
(b) $\Q$ is given by \eqref{QABL}.
\end{center}

The following figures are 2D zoomed plots for different $\bv_i$'s.
\begin{center}
(a)\includegraphics [height=2in]
{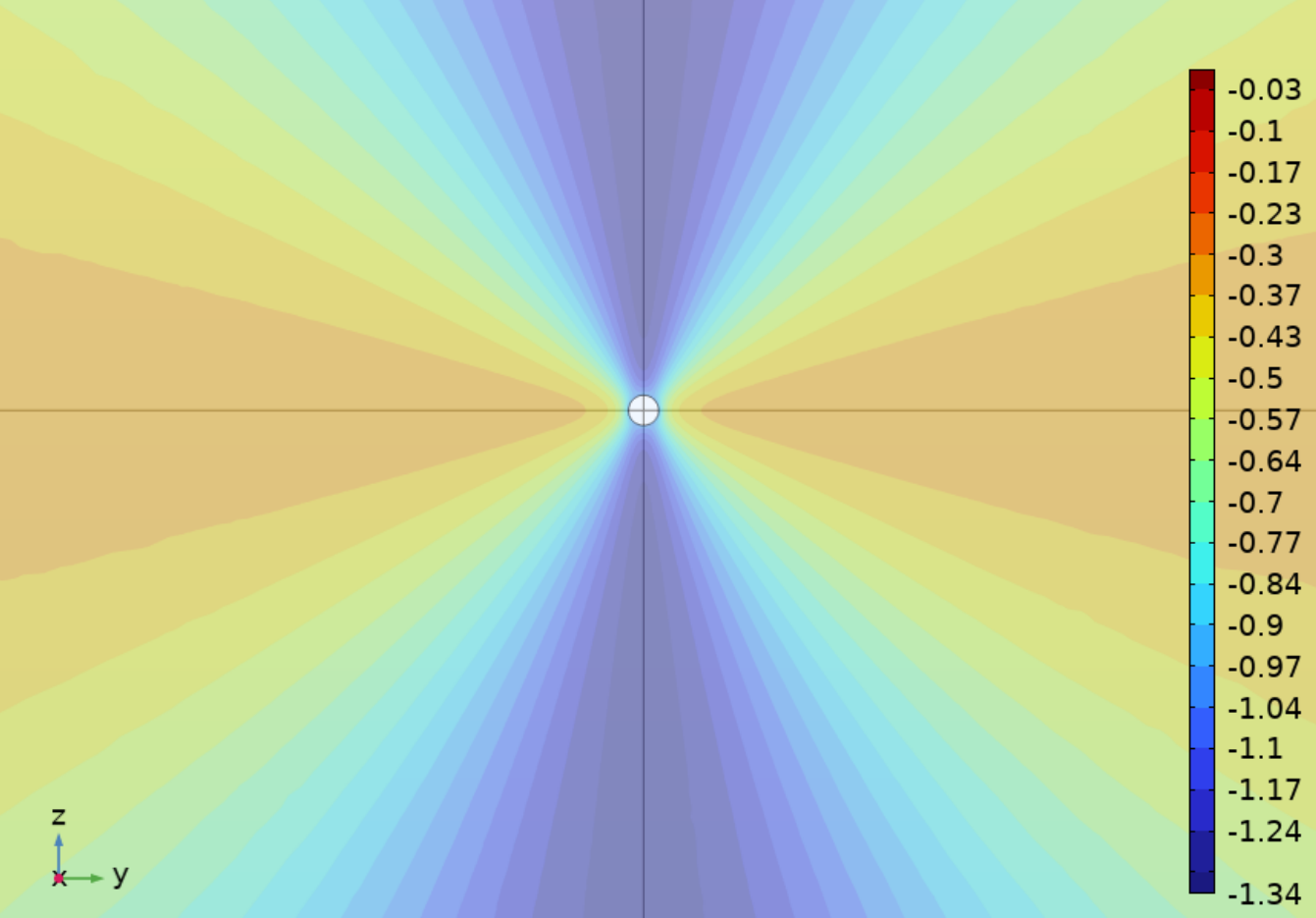}
\,\,\,
(b)\includegraphics [height=2in]
{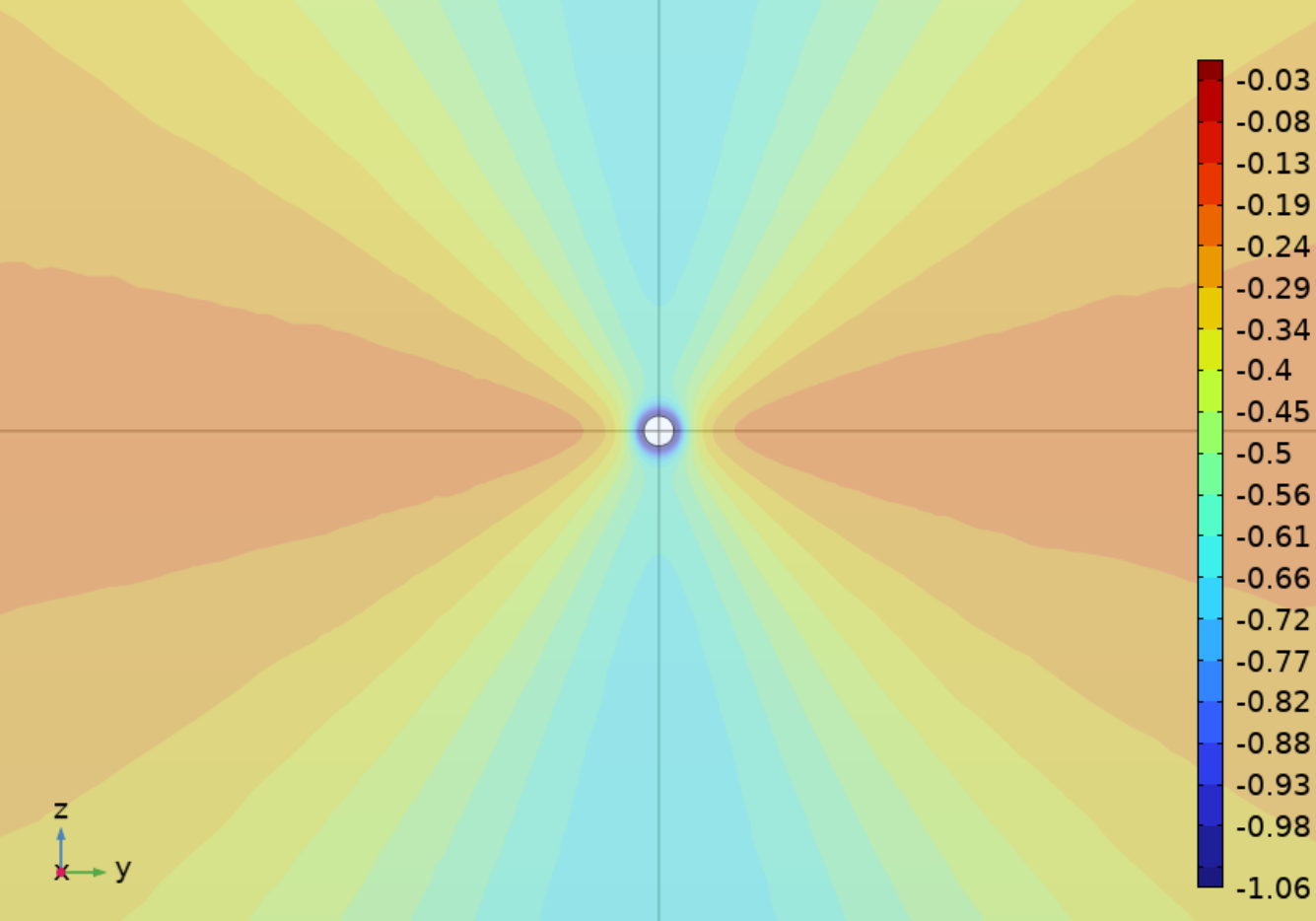}
\\
Figure 3. Zoomed, 2D, rescaled plot of $\bv_1$ in the $yz$-plane.\\
(a) $\Q$ is set to be $\Q_*$;
(b) $\Q$ is given by \eqref{QABL}.
\end{center}

\begin{center}
(a)\includegraphics [height=2in]
{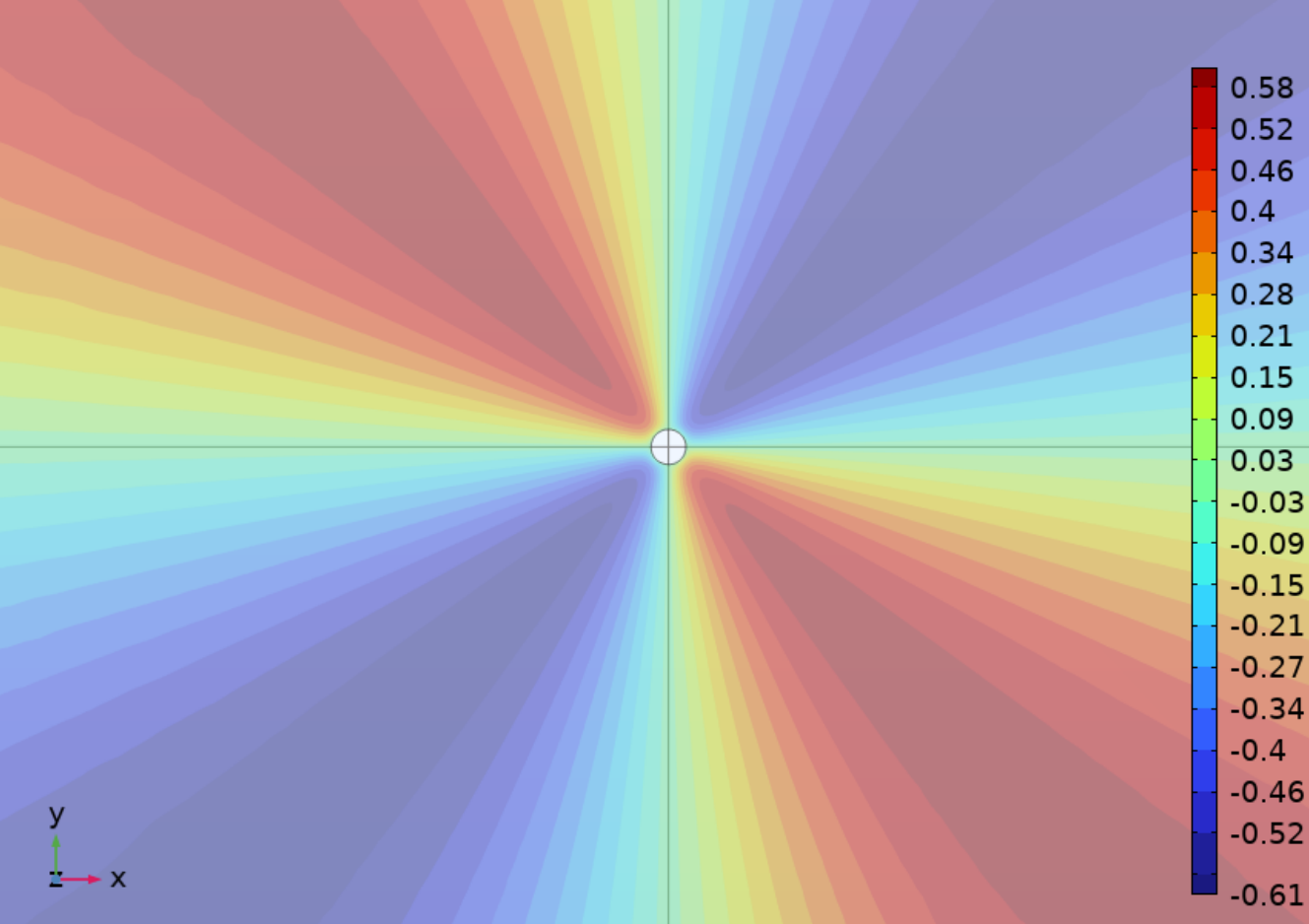}
\,\,\,
(b)\includegraphics [height=2in]
{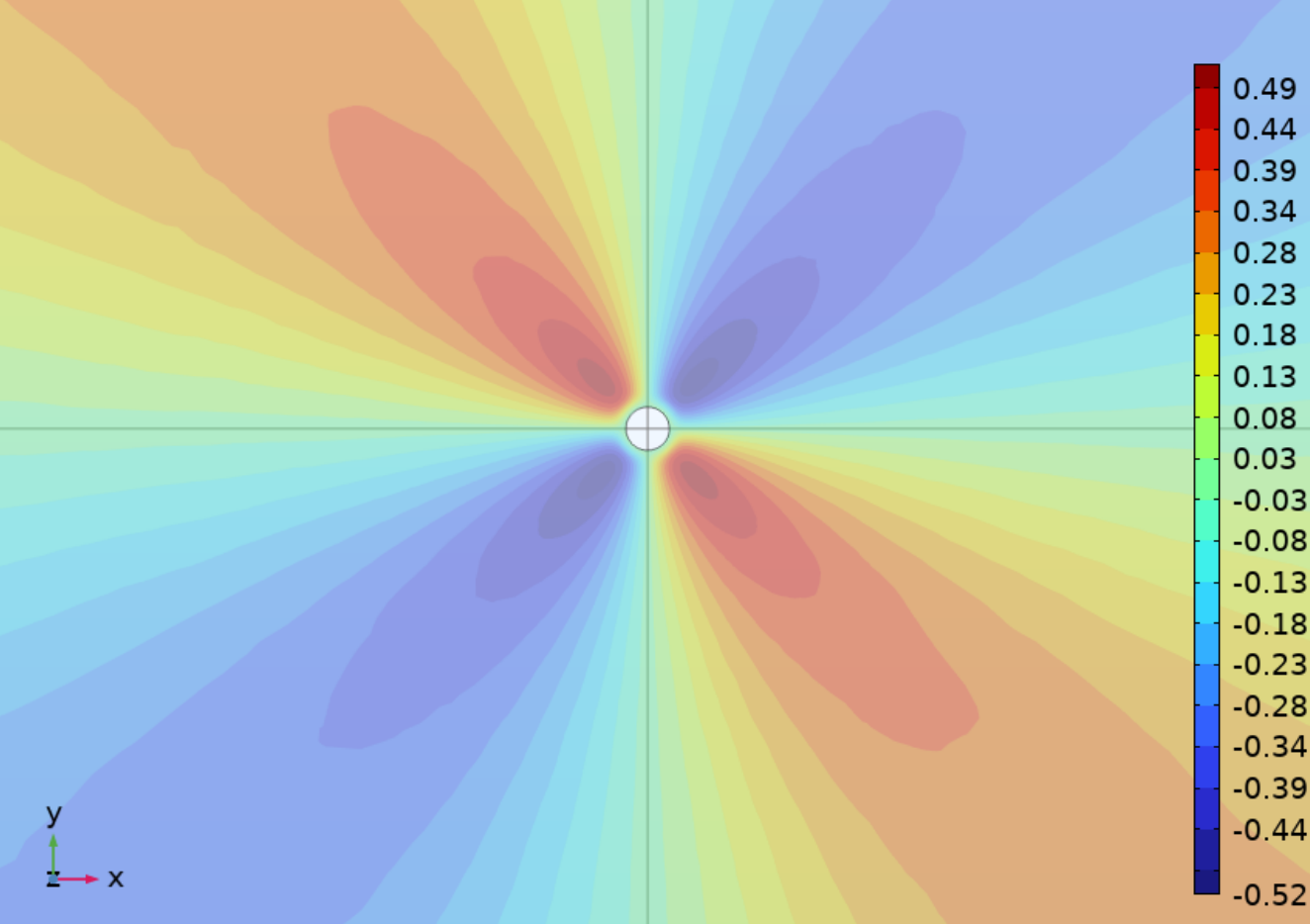}
\\
Figure 4. Zoomed, 2D, rescaled plot of $\bv_2$ in the $xy$-plane.\\
(a) $\Q$ is set to be $\Q_*$;
(b) $\Q$ is given by \eqref{QABL}.
\end{center}
\begin{center}
(a)\includegraphics [height=2in]
{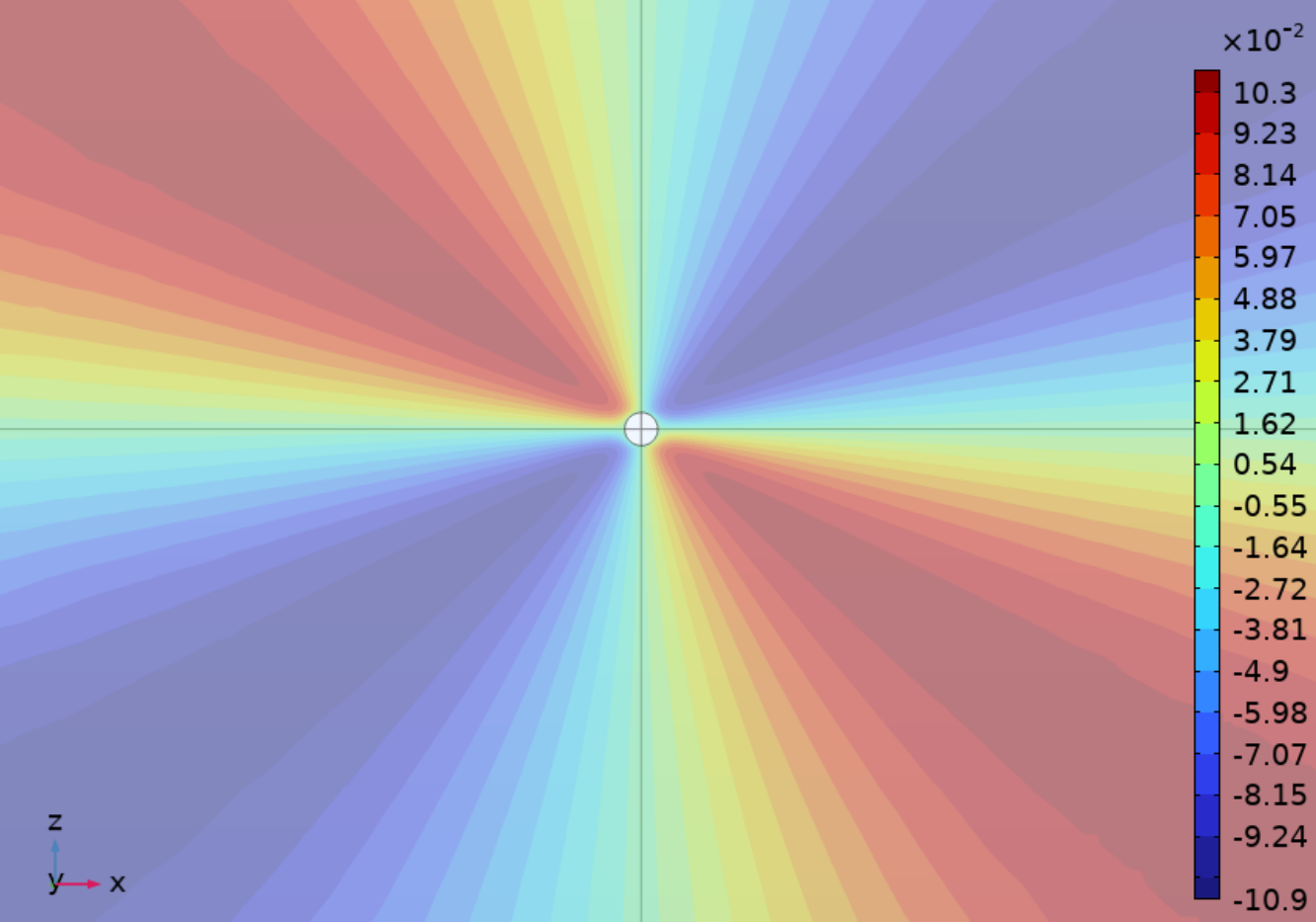}
\,\,\,
(b)\includegraphics [height=2in]
{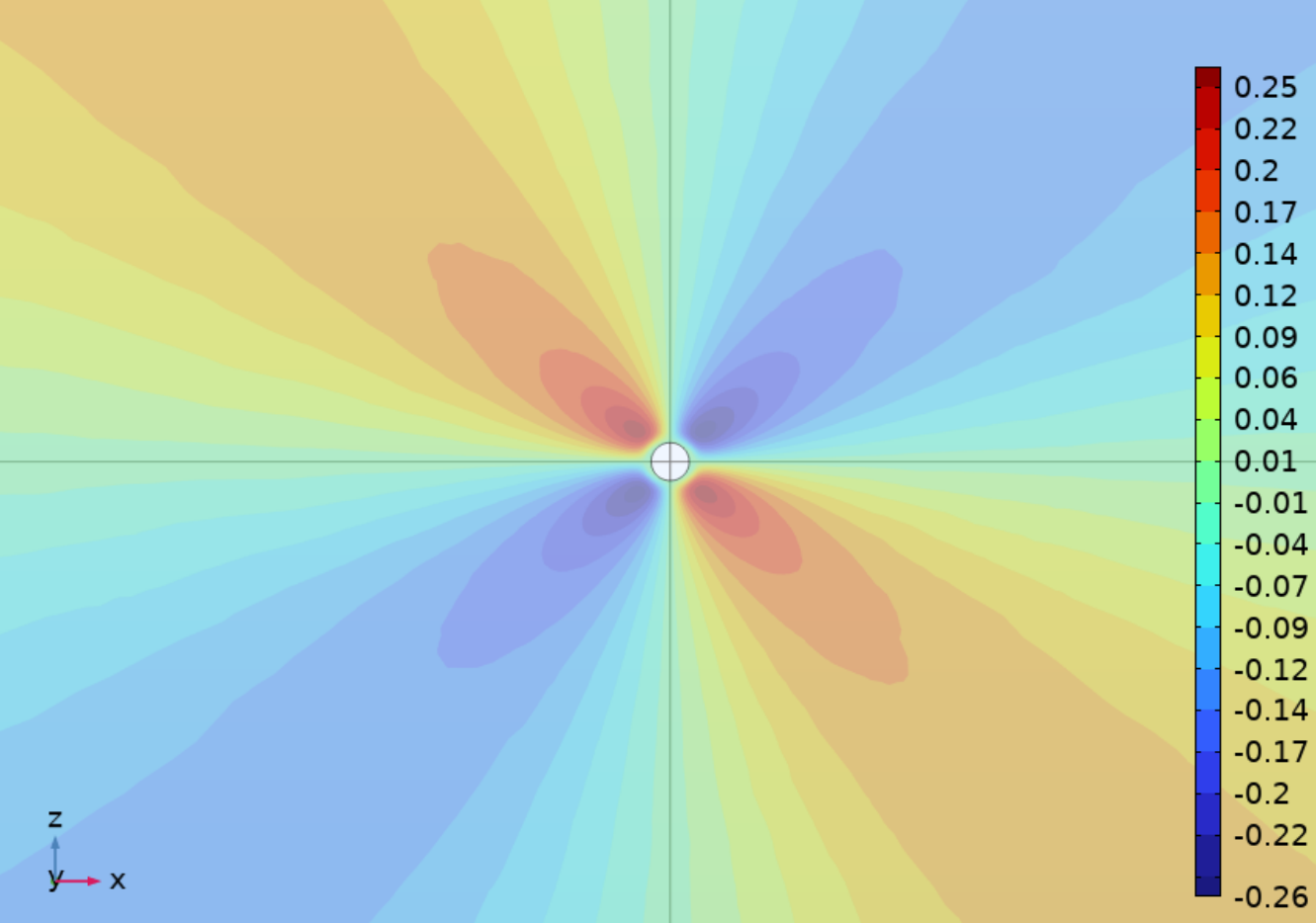}
\\
Figure 5. Zoomed, 2D, rescaled plot of $\bv_3$ in the $xz$-plane.\\
(a) $\Q$ is set to be $\Q_*$;
(b) $\Q$ is given by \eqref{QABL}.
\end{center}

Note that in Figures 1, 2, and 3, there is reflection symmetry for
$\bv_1$ in the $yz$-plane with respect to both the $y$- and $z$-axes. 
This is due to the fact that the $yz$-plane is perpendicular to the 
$xz$-plane, the plane spanned by $\bv_*$ and $\n$. 

\subsection{$\bv_* = e_3$, $\n = e_3$}\label{NumRotSym}
In this section, the $\bv_*$ and $n$ are parallel to each other, both 
pointing in the direction of $z$-axis. In this case, $\bv_3$ should be rotational
symmetric with respect to the $z$-axis. This is clearly demonstrated
in the plot for $\bv_3$ in the $xy$-plane -- see Figure 6.

\begin{center}
(a)\includegraphics [height=2in]
{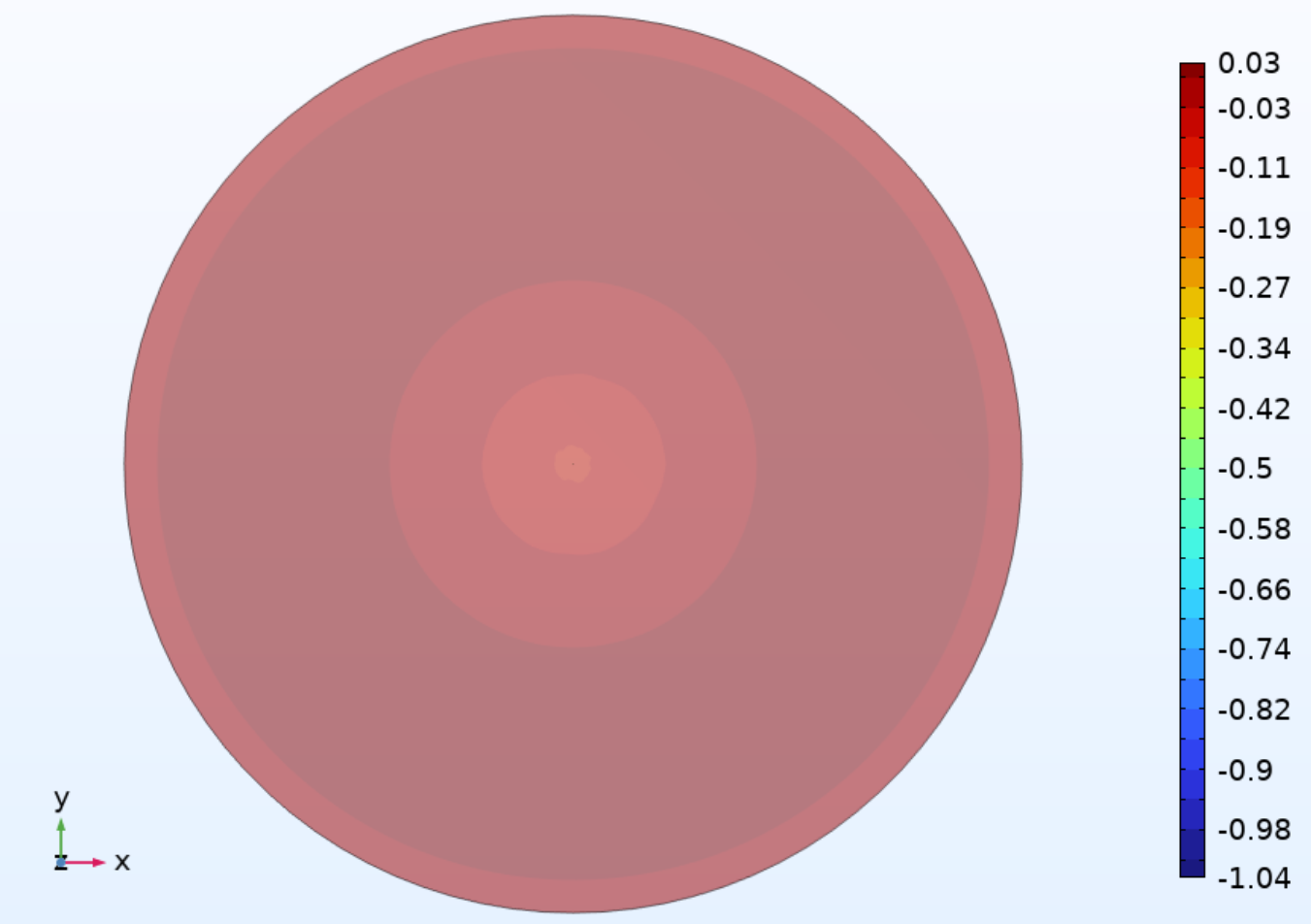}
\,\,\,
(b)\includegraphics [height=2in]
{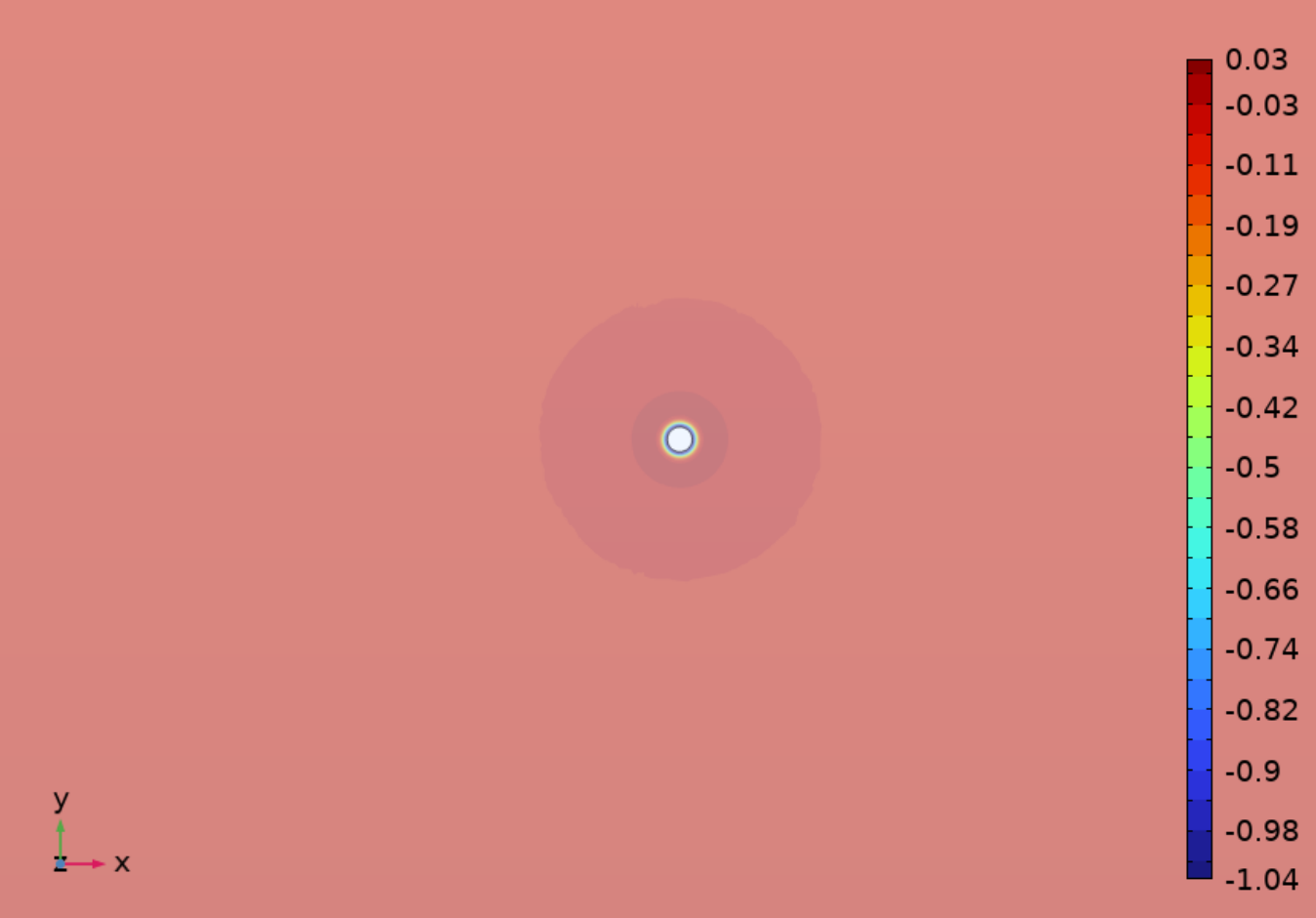}
\\
Figure 6. 2D, rescaled plot of $\bv_3$ in the $xy$-plane:\\
(a) whole computational domain;
(b) zoomed version.
\end{center}

\begin{center}
(a)\includegraphics [height=2in]
{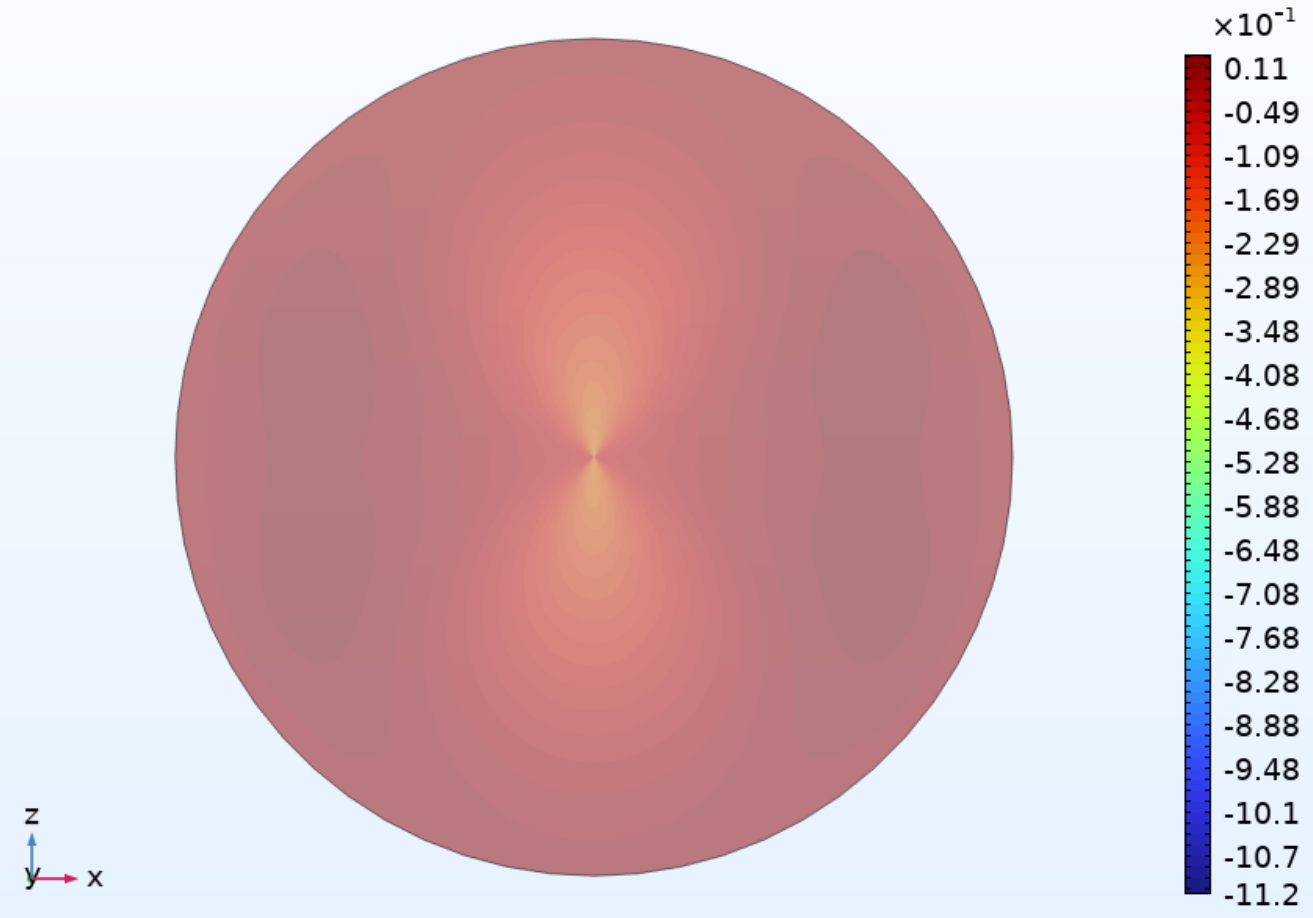}
\,\,\,
(b)\includegraphics [height=2in]
{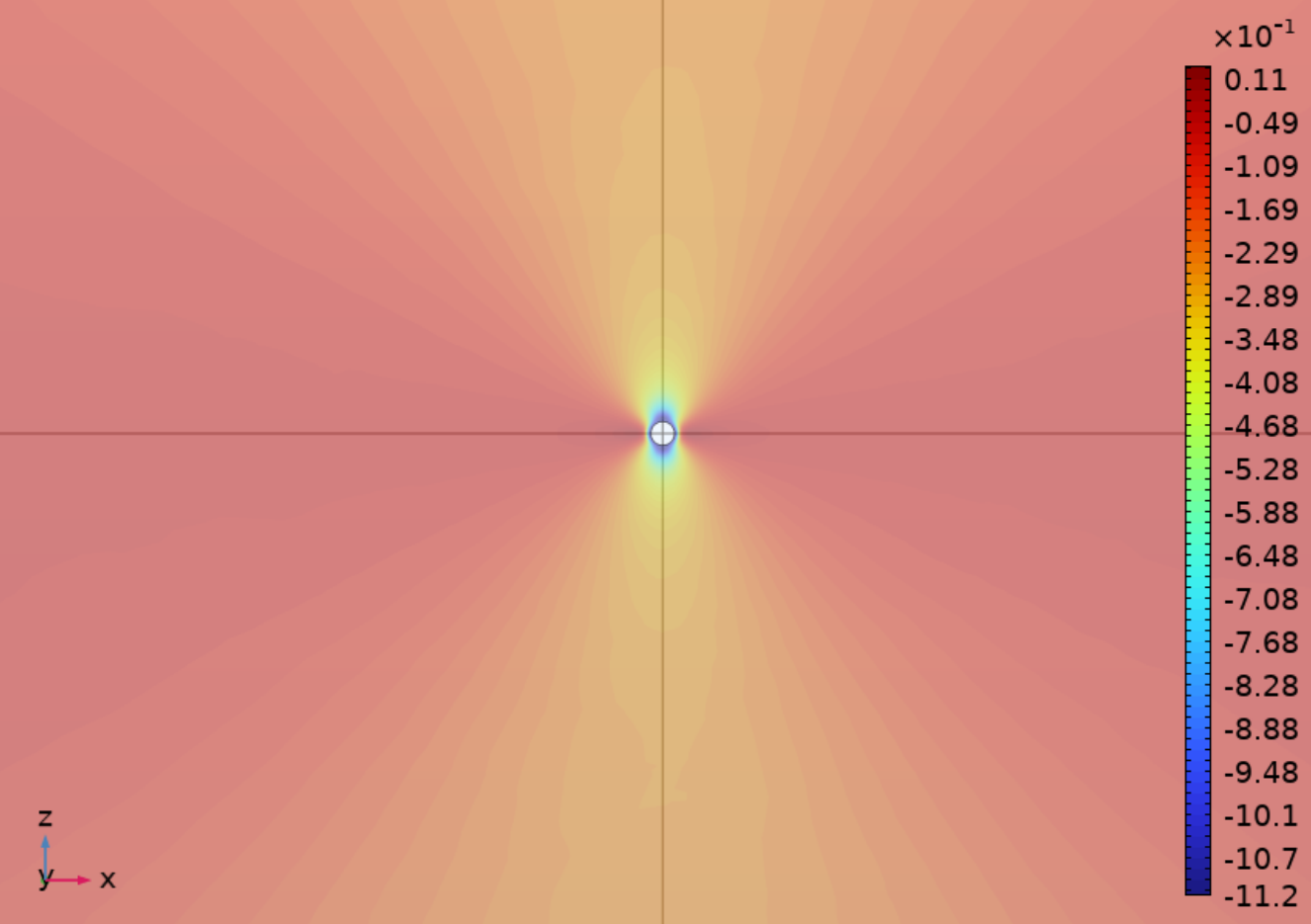}
\\
Figure 7. 2D, rescaled plot of $\bv_3$ in the $xz$-plane:\\
(a) whole computational domain;
(b) zoomed version.
\end{center}

\begin{center}
(a)\includegraphics [height=2in]
{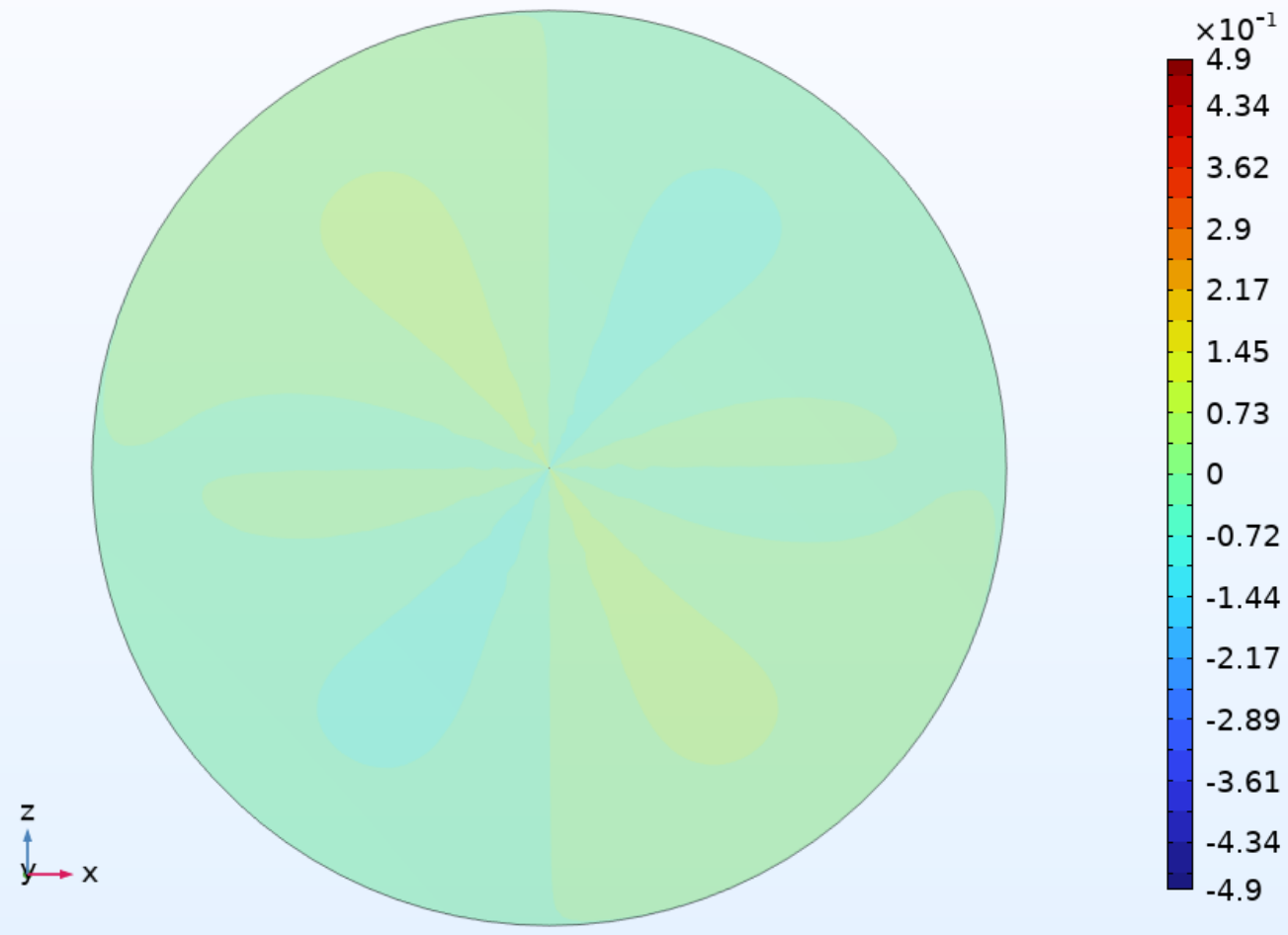}
\,\,\,
(b)\includegraphics [height=2in]
{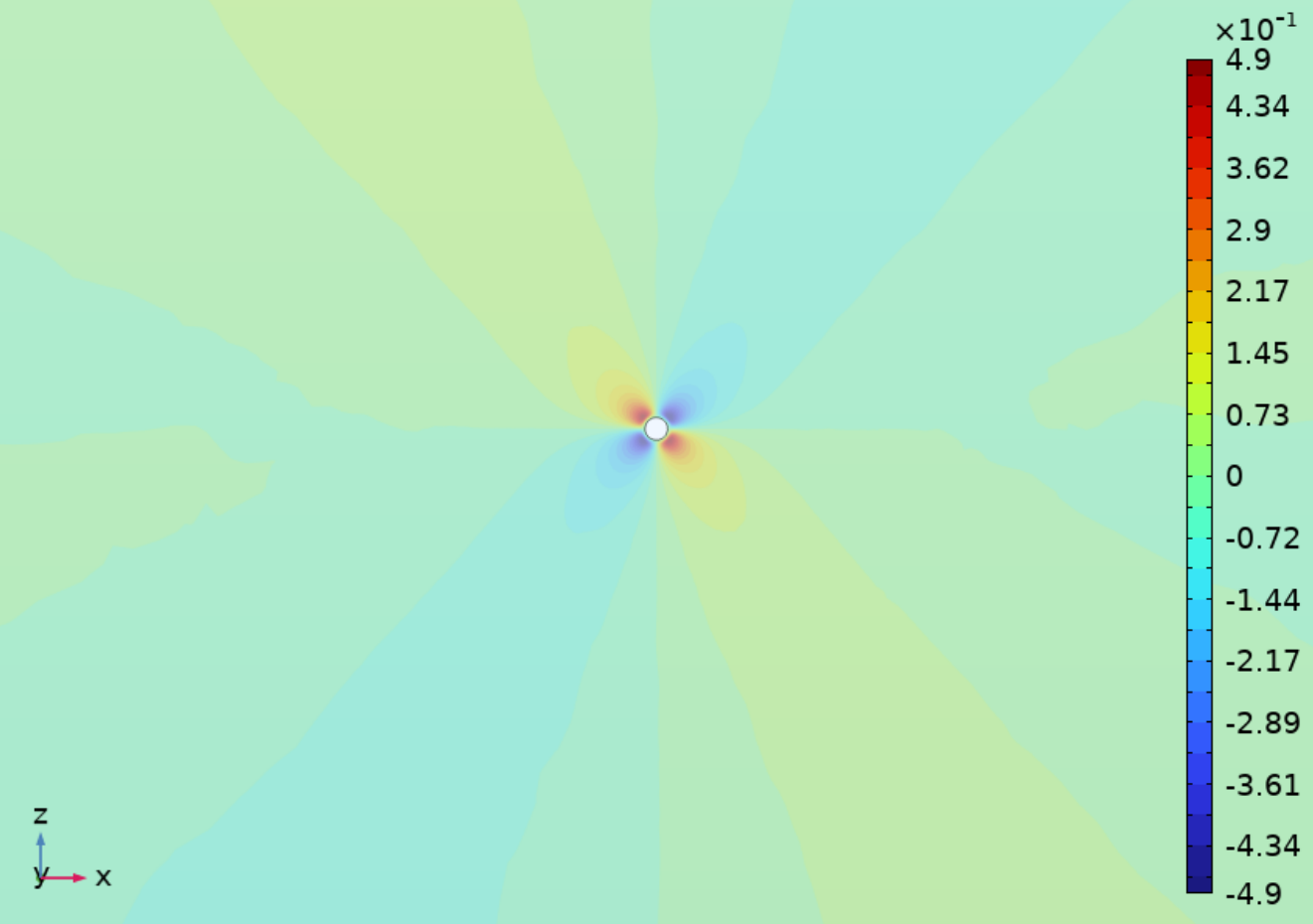}
\\
Figure 8. 2D, rescaled plot of $\bv_1$ in the $xz$-plane:\\
(a) whole computational domain;
(b) zoomed version.
\end{center}

\begin{center}
(a)\includegraphics [height=2in]
{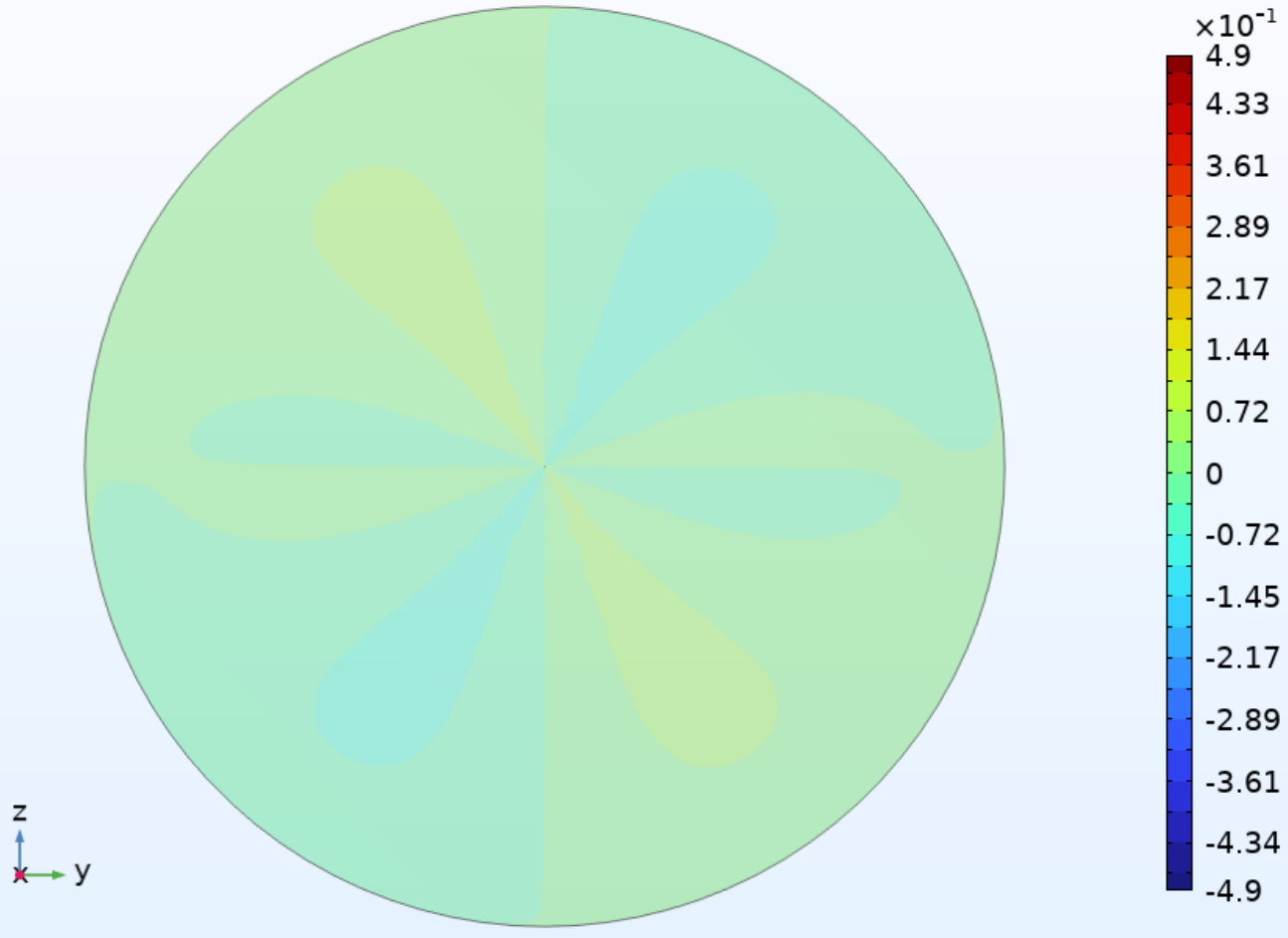}
\,\,\,
(b)\includegraphics [height=2in]
{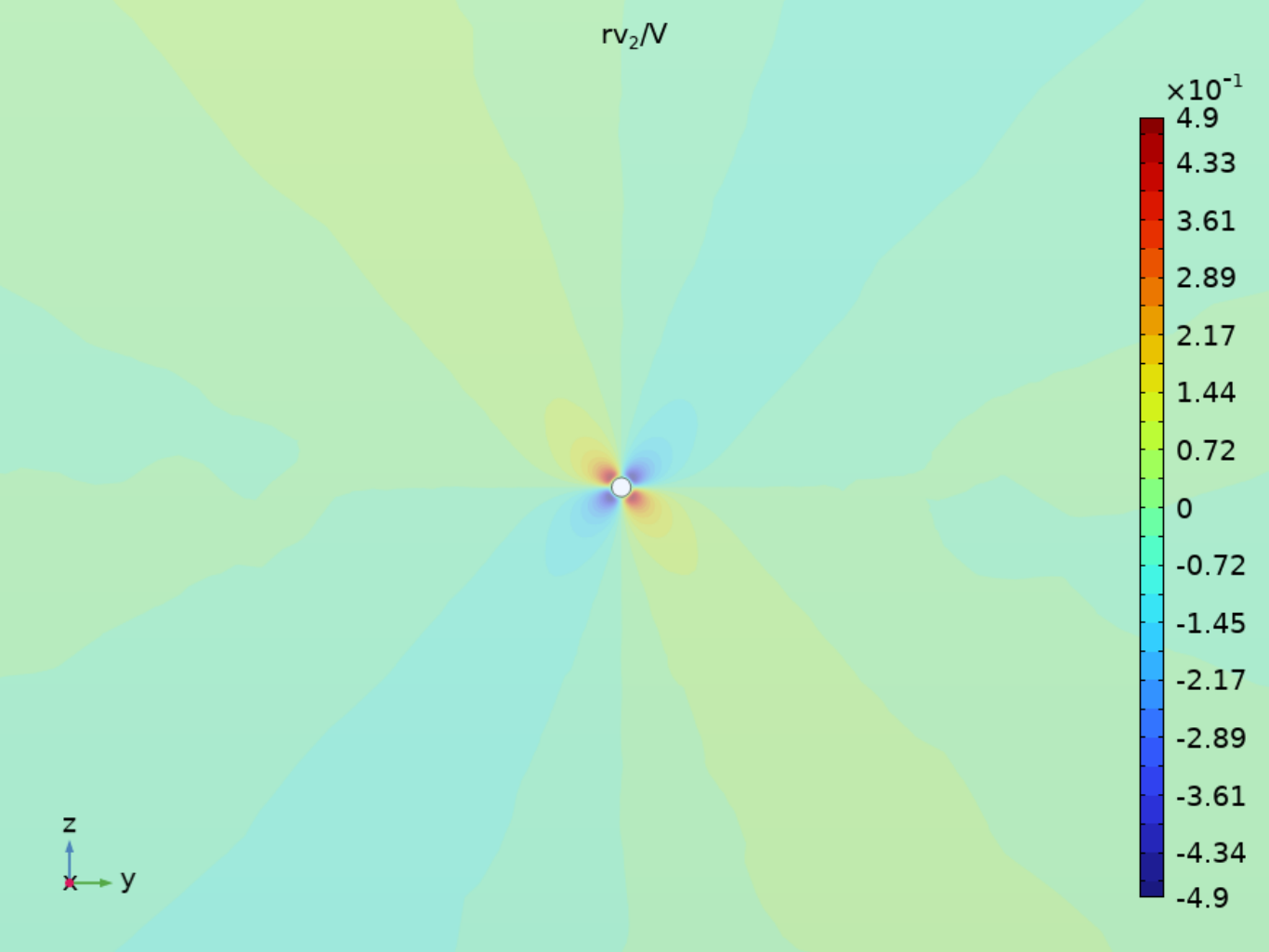}
\\
Figure 9. 2D, rescaled plot of $\bv_2$ in the $yz$-plane:\\
(a) whole computational domain;
(b) zoomed version.
\end{center}
Note that by symmetry the behavior of 
$\bv_1$ in the $xz$-plane and 
$\bv_2$ in the $yz$-plane should be ``identical'', as illustrated by 
Figures 8 and 9.

\subsection{$\bv_* = e_3$, $\n=e_3$}\label{NumGamma}
In this section, we demonstrate the decomposition \eqref{v-v0-varphig}.
This is validated by plotting
\begin{equation}\label{rescale-diff}
r\left(\frac{\bv_i - {\bv_0}_i}{V\gamma_1}\right)
=
r\left(\frac{{\overline{\varphi}_\gamma}_i}{V\gamma_1}\right)
= rO\left(\frac{\overline{\varphi}_\gamma}{|\gamma|}\right)
,\,\,\,
\text{for $i=1,2,3$.}
\end{equation}
with three samples of $\gamma_1$ and $\gamma_2$:
\[
\text{(a) $\gamma_1 = 0.2$, $\gamma_2 = 0.18$;
(b) $\gamma_1 = 0.1$, $\gamma_2 = 0.09$;
(c) $\gamma_1 = 0.05$, $\gamma_2 = 0.045$.
}
\]
As $|\gamma|\longrightarrow0$, we expect 
$\displaystyle
r\left(\frac{\overline{\varphi}_\gamma}{|\gamma|}\right)
$ converging to some fixed function.
This is clearly demonstrated in the following plots, 
in terms of both the pattern and order of magnitude.

\begin{center}
(a)\includegraphics [height=1.3in]
{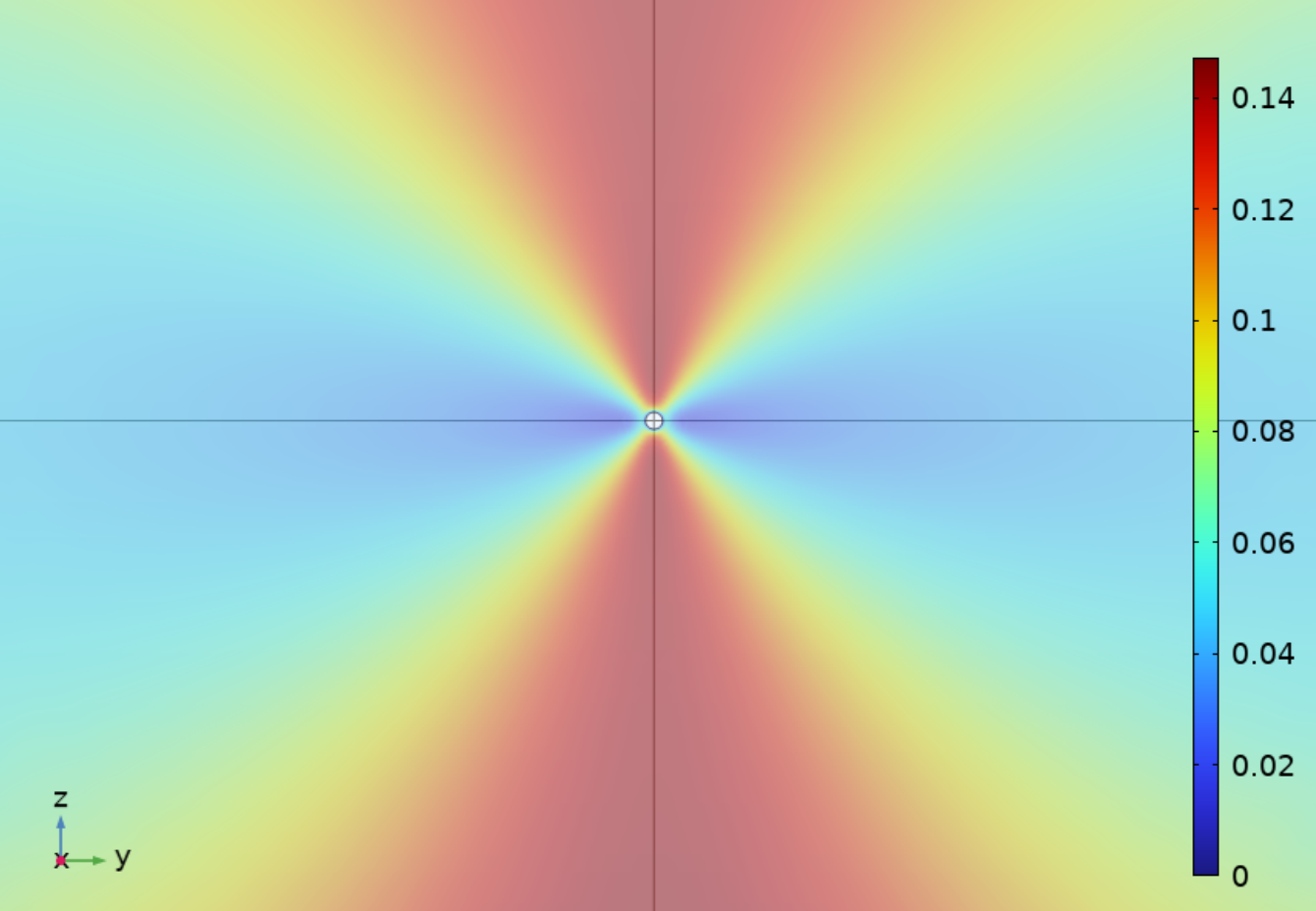}
(b)\includegraphics [height=1.3in]
{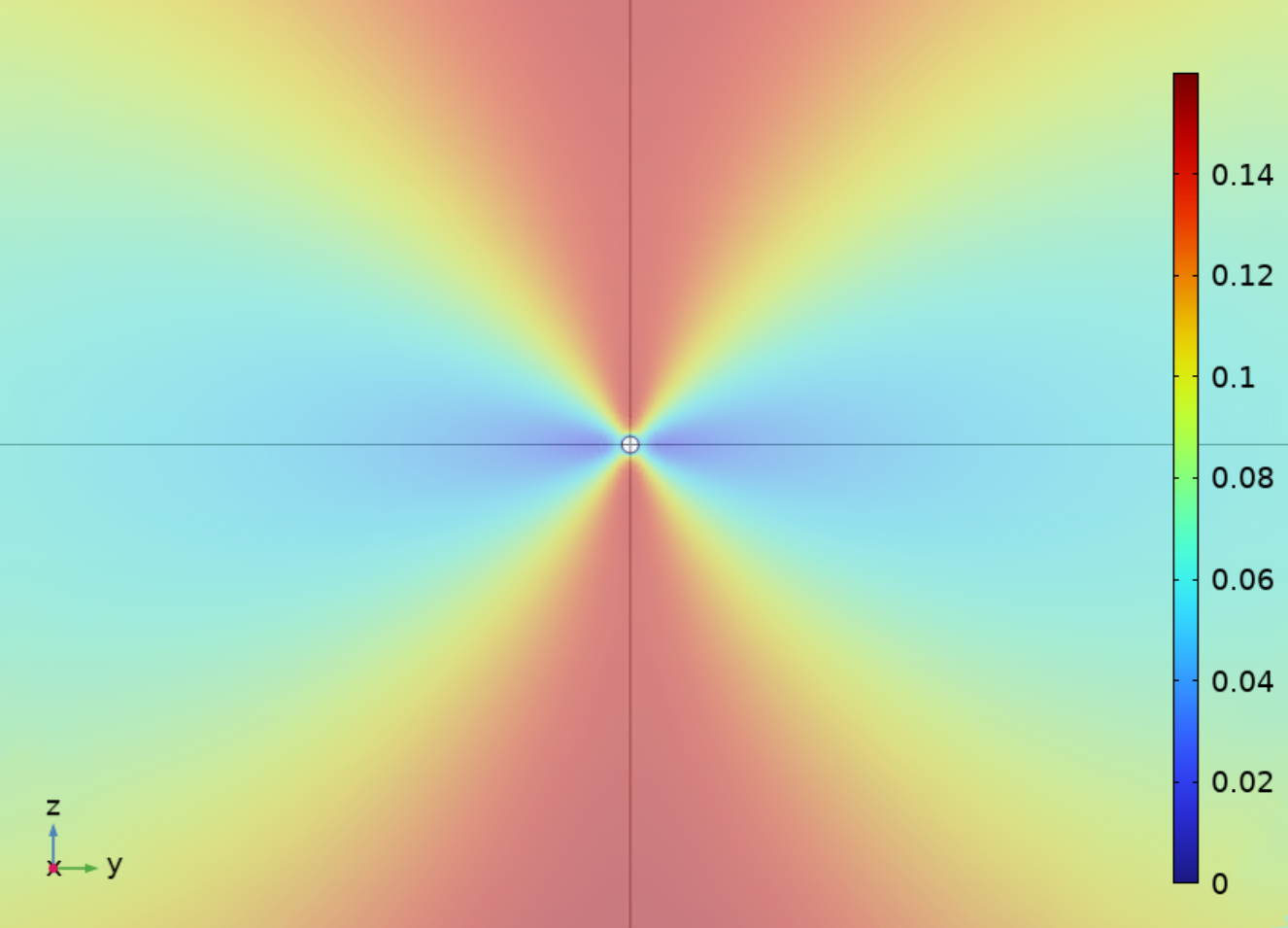}
(c)\includegraphics [height=1.3in]
{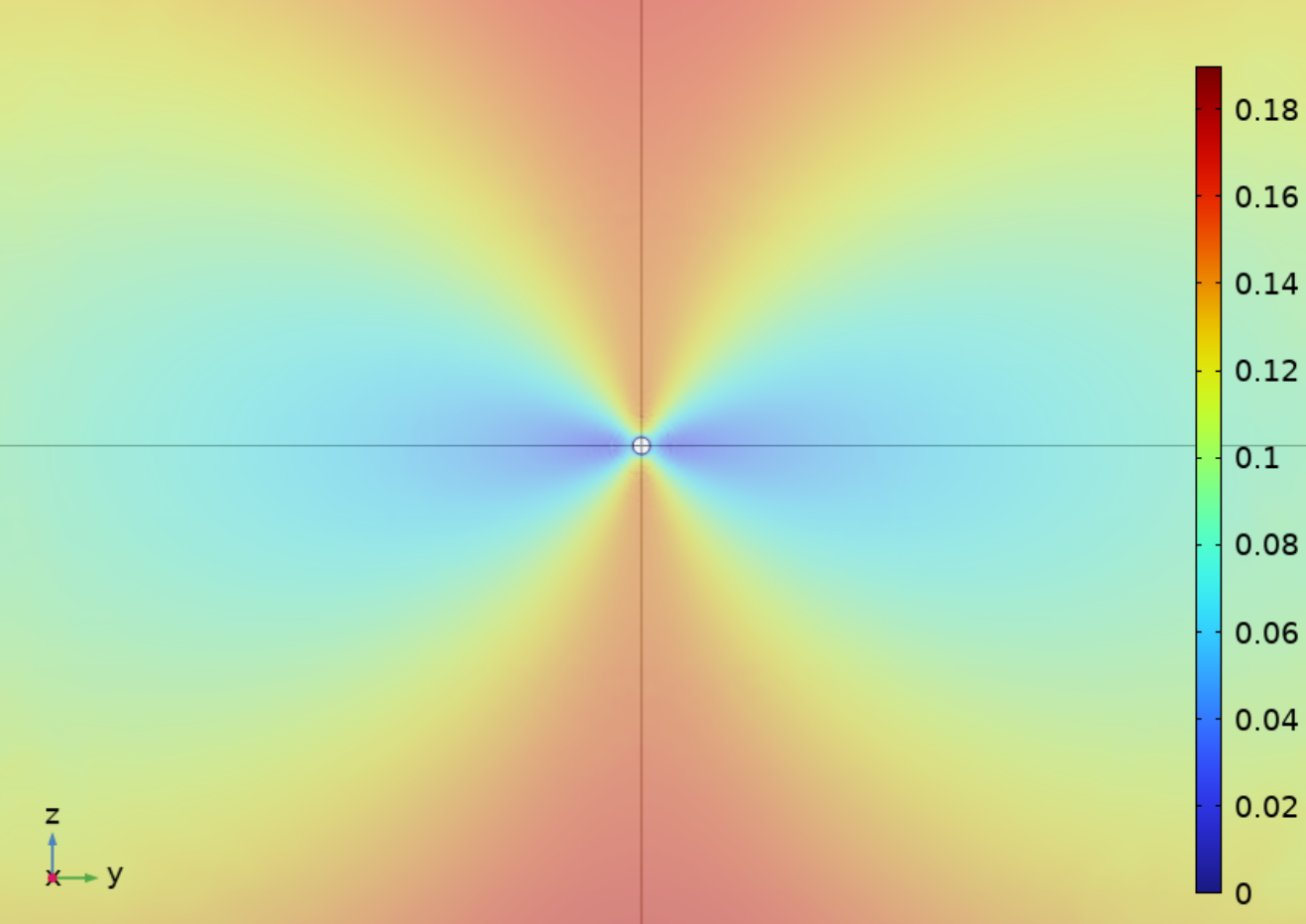}
\\
Figure 10. 2D, rescaled plot of $\bv_1 - {\bv_0}_1$ in the $yz$-plane:\\
(a) $\gamma_1 = 0.2$, $\gamma_2 = 0.18$;
(b) $\gamma_1 = 0.1$, $\gamma_2 = 0.09$;
(c) $\gamma_1 = 0.05$, $\gamma_2 = 0.045$.
\end{center}
\begin{center}
(a)\includegraphics [height=1.3in]
{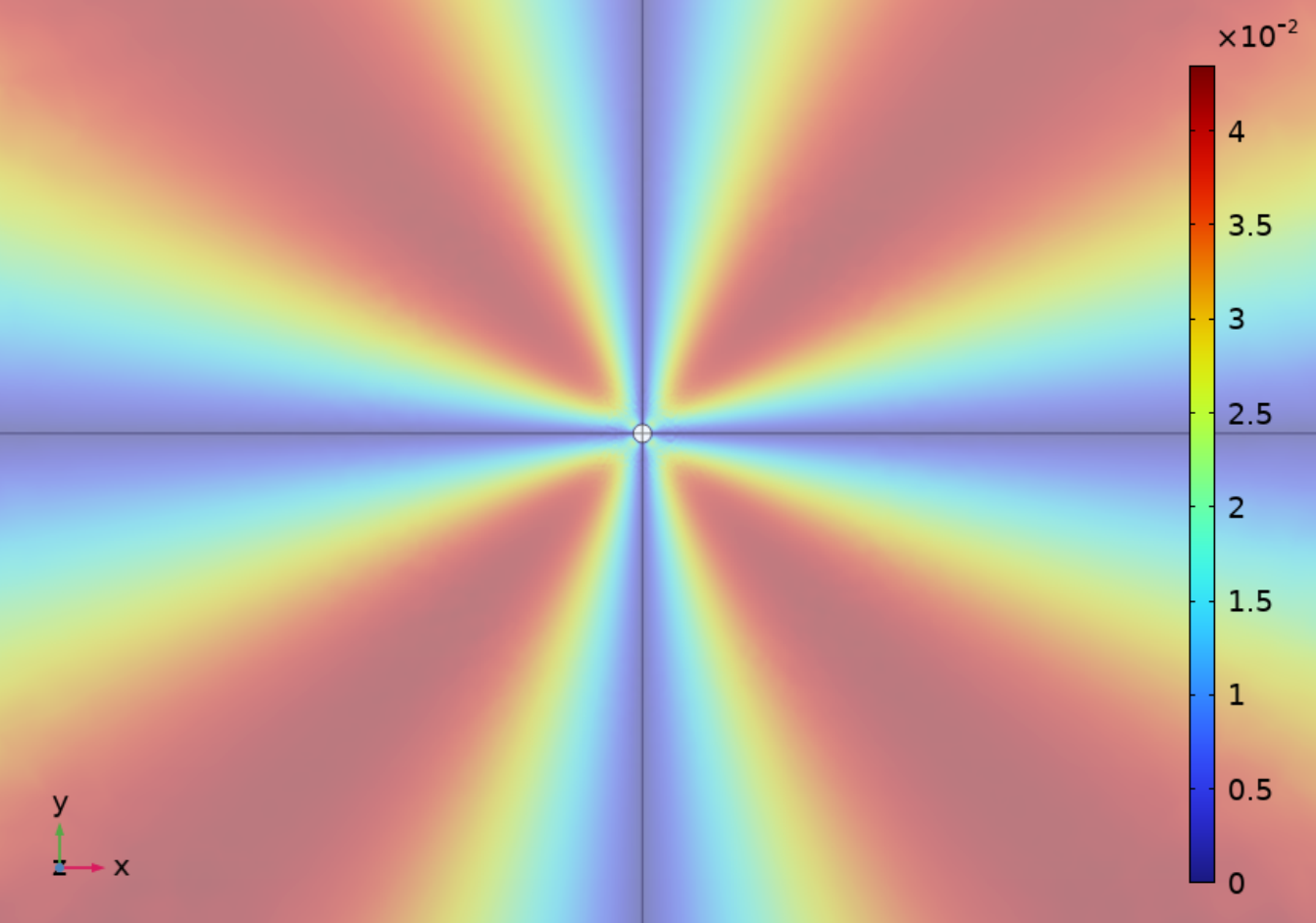}
(b)\includegraphics [height=1.3in]
{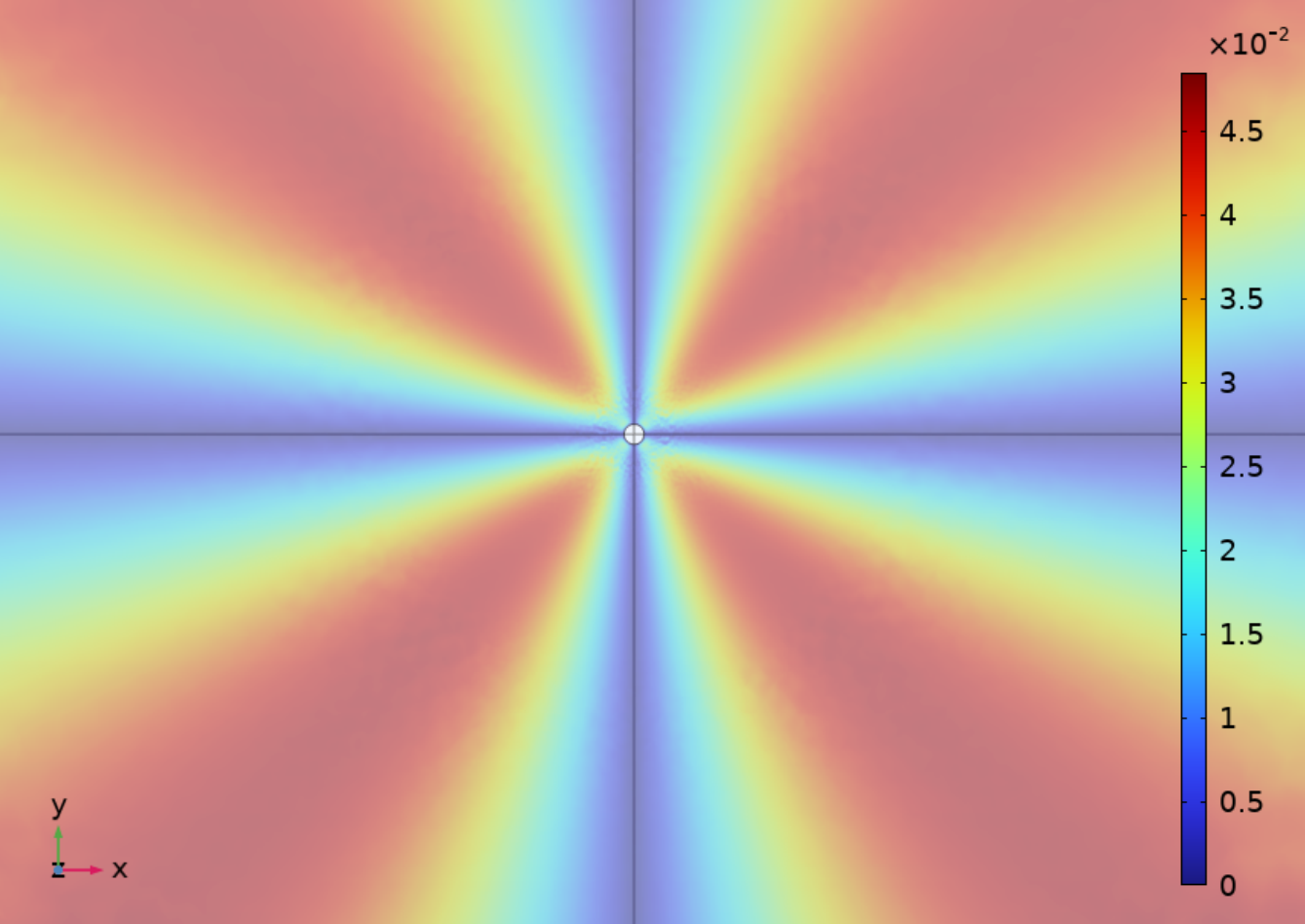}
(c)\includegraphics [height=1.3in]
{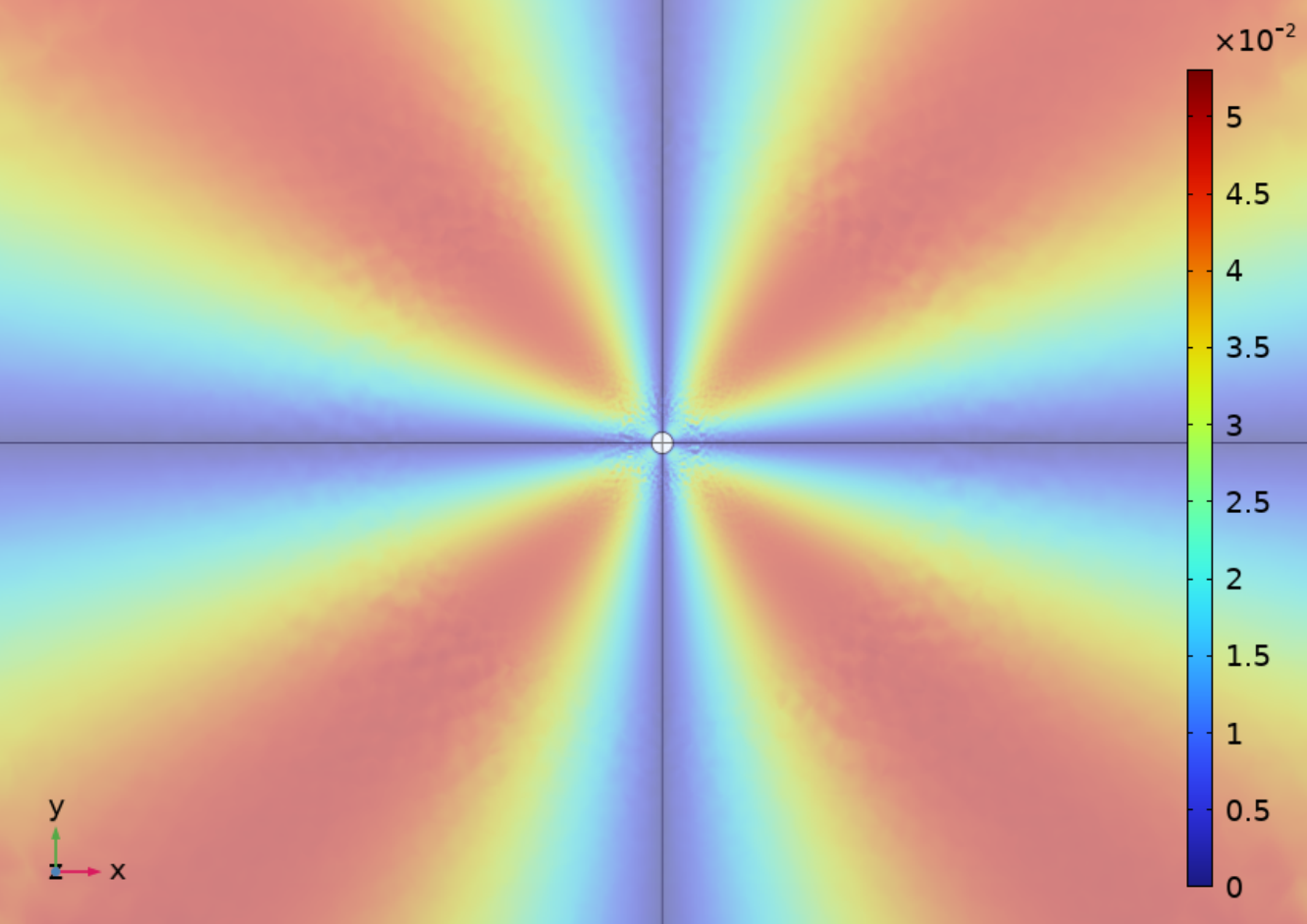}
\\
Figure 11. 2D, rescaled plot of $\bv_2 - {\bv_0}_2$ in the $xy$-plane:\\
(a) $\gamma_1 = 0.2$, $\gamma_2 = 0.18$;
(b) $\gamma_1 = 0.1$, $\gamma_2 = 0.09$;
(c) $\gamma_1 = 0.05$, $\gamma_2 = 0.045$.
\end{center}
\begin{center}
(a)\includegraphics [height=1.3in]
{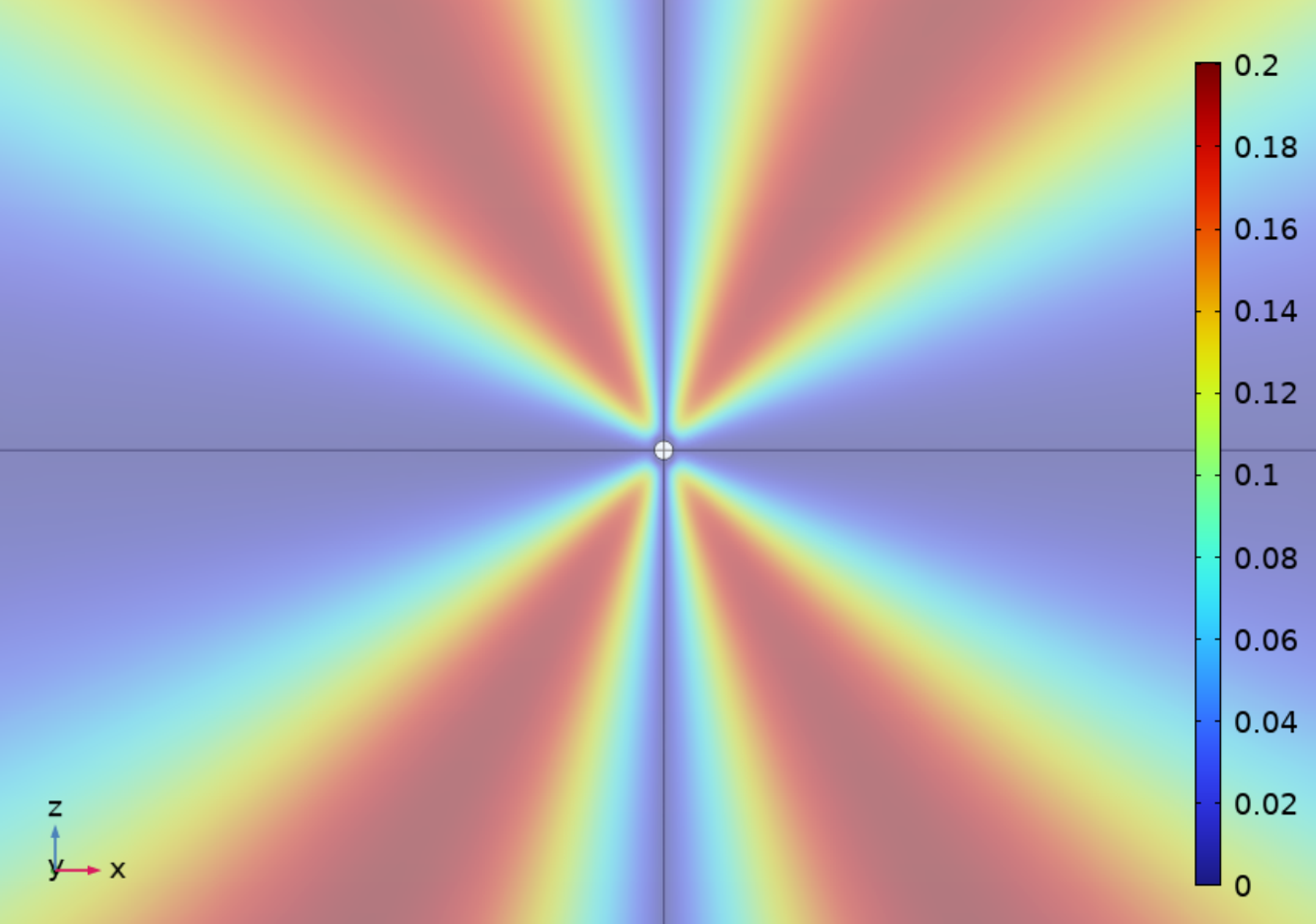}
(b)\includegraphics [height=1.3in]
{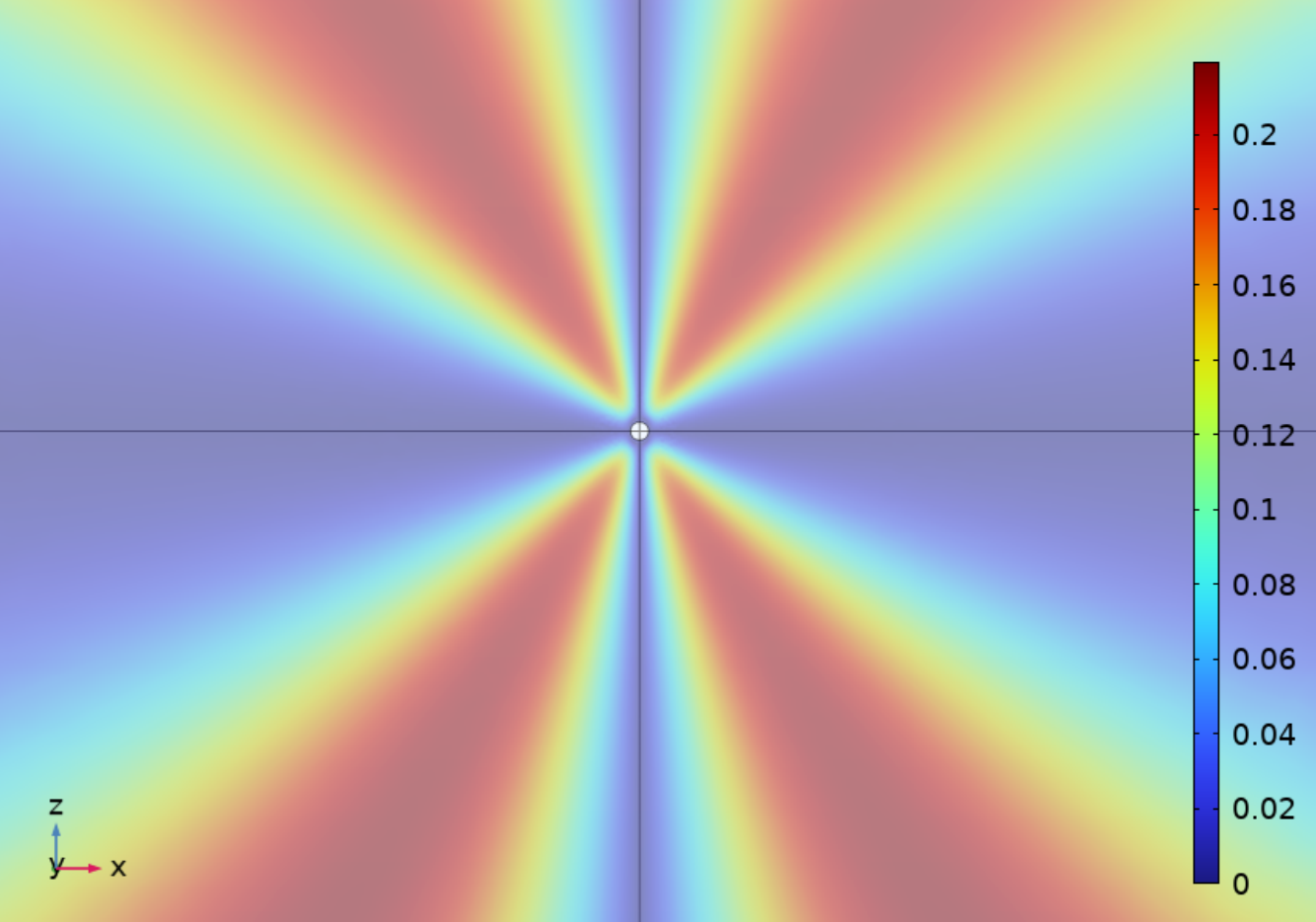}
(c)\includegraphics [height=1.3in]
{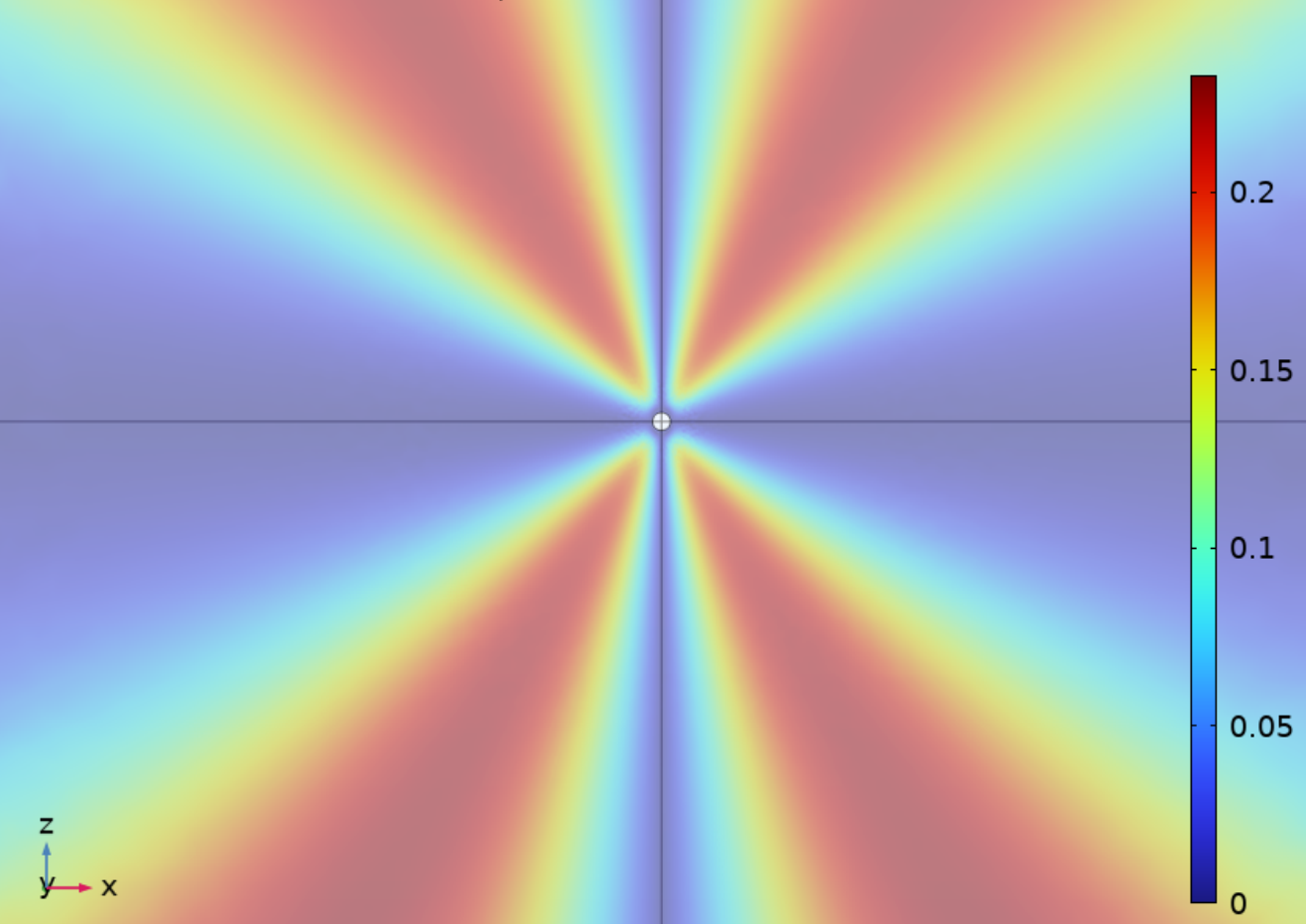}
\\
Figure 12. 2D, rescaled plot of $\bv_3 - {\bv_0}_3$ in the $xz$-plane:\\
(a) $\gamma_1 = 0.2$, $\gamma_2 = 0.18$;
(b) $\gamma_1 = 0.1$, $\gamma_2 = 0.09$;
(c) $\gamma_1 = 0.05$, $\gamma_2 = 0.045$.
\end{center}

\section{Acknowledgements}
The work of DG was partially supported by the NSF grant DMS-2106551. The authors would like to acknowledge useful discussions with Leonid Berlyand and Mykhailo Potomkin.

\section{Declaration of Interests}
The authors report no conflict of interest.

\appendix
\section{Far-field behavior of Stokes system}\label{AppFFZ}
\renewcommand{\G}{{\mathsf{G}}}
Consider equation \eqref{eq:concise}
and the representation of its solution in \eqref{bulk.int.soln}.
The asymptotics of the boundary integrals are given in
step {\bf (II)} of the proof in Section \ref{Exist.Soln}.
Here, we will analyze the asymptotics of the bulk integral
\[
I(x) = \int_{\mathbb{R}^3\backslash\B_1(0)} \G(x-y)\f(y)\,dy.
\]
The property of homogeneous of degree $-1$ \eqref{G.Fct} for $\G=\G_\gamma$ 
plays an important role in our analysis.
The precise asymptotics naturally also depends on the far-field behavior of
$\f$. We present these results in the following cases.

\subsection{Case I} 
We assume that $\|\f\|_{L^\infty(\mathbb{R}^3\backslash\B_1(0))} 
< \infty$ and
\begin{equation}
\int_{\R^3\backslash\bB_1(0)}|\f(y)|\,dy < \infty.
\end{equation}
We compute
\begin{eqnarray}
&&\int_{\mathbb{R}^3\backslash\bB_1(0)} \G(x-y)\f(y)\,dy\nonumber\\
&=&
\G(x)\left(\int_{\mathbb{R}^3\backslash\bB_1(0)} \f(y)\,dy\right)
+\int_{\mathbb{R}^3\backslash\bB_1(0)} (\G(x-y)-\G(x))\f(y)\,dy\nonumber\\
&=&
\G(x)\left(\int_{\mathbb{R}^3\backslash\bB_1(0)} \f(y)\,dy\right)
+\frac{1}{|x|}\int_{\mathbb{R}^3\backslash\bB_1(0)}
\left(\G(\hat{x}-\frac{y}{|x|})-
\G(\hat{x})\right)\f(y)\,dy\nonumber\\
&=&
\G(x)\left(\int_{\mathbb{R}^3\backslash\bB_1(0)} \f(y)\,dy\right)
+ o\left(\frac{1}{|x|}\right).
\end{eqnarray}

In order to characterize the $o\left(\frac{1}{|x|}\right)$ term, we assume that
\begin{equation}
\f(y) \lesssim \frac{1}{|y|^4}.
\end{equation}
Let $L \gg 1$. Then we have
\begin{eqnarray}
&&\left|\frac{1}{|x|}\int_{\mathbb{R}^3\backslash\bB_1(0)}
\left(\G(\hat{x}-\frac{y}{|x|})-
\G(\hat{x})\right)\f(y)\,dy\right|\nonumber\\
&=&
\left|\frac{1}{|x|}
\int_{1 \leq |y| < \frac{|x|}{L}}
+\int_{\frac{|x|}{L} \leq |y| < L|x|}
+\int_{L|x| \leq |y|}
\left(\frac{\G(\widehat{\hat{x}-\frac{y}{|x|}})}
{\left|\hat{x}-\frac{y}{|x|}\right|}-
\G(\hat{x})\right)\f(y)\,dy\right|
\nonumber\\
& \lesssim & 
\frac{1}{|x|}\int_{1 \leq |y| < \frac{|x|}{L}}
\frac{|y|}{|x|}\frac{1}{|y|^4}\,dy 
+ \frac{1}{|x|}\int_{\frac{|x|}{L} \leq |y| < L|x|}
\frac{1}{|y|^4}\,dy
+ \frac{1}{|x|}\int_{L|x| \leq |y|}
\frac{|x|}{|y|}\frac{1}{|y|^4}\,dy\nonumber\\
& \lesssim & 
\frac{1}{|x|^2}\int_{1}^{\frac{|x|}{L}} \frac{r^2}{r^3}\,dr
+\frac{1}{|x|}\int_{\frac{|x|}{L}}^{L|x|} \frac{r^2}{r^4}\,dr
+\int_{L|x|}^\infty\frac{r^2}{r^5}\,dr\nonumber\\
& \lesssim & 
\frac{\log|x|}{|x|^2} + \frac{1}{|x|^2}\nonumber\\
& \lesssim & 
\frac{\log|x|}{|x|^2}\label{logx2}.
\end{eqnarray}
Note that in the above, we have used the fact that 
$\frac{1}{\left|\hat{x}-\frac{y}{|x|}\right|}$ is an
integrable singularity at $y\sim x$ in $\mathbb{R}^3$.

\subsection{Case II} 
Next, we assume that $f$ has the following large spatial asymptotic behavior:
\begin{equation}\label{ass:fln3}
|\f| \leq O\left(\frac{A+ B\ln r}{r^k}\right)
\,\,\,\text{for some $k\geq 3$ and all $r \gg 1$.}
\end{equation}
Then we compute,
\begin{eqnarray}
&&\left|\int_{\mathbb{R}^3\backslash\bB_1(0)} \G(x-y)\f(y)\,dy\right|\nonumber\\
&\lesssim&
\int_{\mathbb{R}^3\backslash\bB_1(0)}\frac{A + B\ln|y|}{|y|^k}\frac{1}{|x-y|}\,dy
\nonumber\\
& \lesssim & 
\frac{1}{|x|}\int_{\mathbb{R}^3\backslash\B_{\frac{1}{|x|}}(0)}
\frac{A + B\ln|z| + B\ln|x|}{|x|^k|z|^k}\frac{1}{|\widehat{x}-z|}\,|x|^3 dz
\quad\text{(where $\widehat{x} = \frac{x}{|x|}$ and $z=\frac{y}{|x|}$)}
\nonumber\\
& \lesssim & 
\frac{1}{|x|^{k-2}}\left[
\int_1^\infty \frac{A + B\ln|r| + B\ln|x|}{r^{k+1}}\,r^2dr + 
\int_{\frac{1}{|x|}}^1 \frac{A + B\ln|r| + B\ln|x|}{r^{k}}\,r^2dr 
\right]
\nonumber\\
&&\quad\text{(note that the singularity at $\widehat{x}$ is integrable)}
\nonumber\\
& \lesssim & \left\{
\begin{array}{ll}
\frac{A\ln|x| + B\ln|x| + B(\ln|x|)^2}{|x|}, & \text{if $k=3$,}\\
\frac{A + B\ln|x|}{|x|}, & \text{if $k>3$.}
\end{array}
\right.\label{bulk.abs.int.est}
\end{eqnarray}
Similarly,
\begin{eqnarray}
\left|\int_{\mathbb{R}^3\backslash\bB_1(0)} - \q(x-y)\cdot\f(y)\,dy\right|
&\lesssim&
\int_{\mathbb{R}^3\backslash\bB_1(0)}\frac{A + B\ln|y|}{|y|^k}
\frac{1}{|x-y|^2}\,dy\nonumber\\
& \lesssim & \left\{
\begin{array}{ll}
\frac{A\ln|x| + B\ln|x| + B(\ln|x|)^2}{|x|^2}, & \text{if $k=3$,}\\
\frac{A + B\ln|x|}{|x|^2}, & \text{if $k>3$.}
\end{array}
\right.
\end{eqnarray}
Note that for $k=3$, even if $B=0$, a $\ln|x|$-term can appear
in the bulk integrals.

\subsection{Case III} 
In order to do a more careful analysis of the case $k=3$ of which we are 
concerned most, we further assume that 
$\f$ is homogeneous of degree $-3$, i.e.
\begin{equation}\label{ass:f-3}
\f(\lambda x) = \lambda^{-3}\f(x)\,\,\,\text{so that}\,\,\,
\f(x) = \frac{\f(\hat{x})}{|x|^3}.
\end{equation}

Note that 
the integrand $\G(x-\cdot)\f(\cdot) \sim \frac{1}{|y|^4}$ 
is integrable $L^1(\R^3\backslash \bB_1(0))$. Hence we can use Fubini's Theorem
to compute $dy$ iteratively as $r^2\,d\omega\,dr$, where
$r = |y|$ and $\omega=\hat{y}$. To this end, 
we have
\begin{eqnarray*}
\int_{\R^3\backslash \bB_1(0)} \G(x-y)\f(y)\,dy
&=&\int_{\mathbb{S}^2}
\left[\int_1^\infty\G(x-r\omega)\frac{\f(\omega)}{r^3}r^2\,dr\right]
\,d\omega\\
&=&\frac{1}{|x|}\int_{\mathbb{S}^2}
\left[\int_1^\infty\G(\hat{x}-\frac{r}{|x|}\omega)\frac{1}{r}\,dr\right]
\f(\omega)\,d\omega\\
&=&\frac{1}{|x|}\int_{\mathbb{S}^2}
\left[\int_{\frac{1}{|x|}}^\infty\G(\hat{x}-r\omega)\frac{1}{r}\,dr\right]
\f(\omega)\,d\omega
\end{eqnarray*}

For $|x| > 2$, we write
\begin{eqnarray*}
\int_{\frac{1}{|x|}}^\infty\G(\hat{x}-r\omega)\frac{1}{r}\,dr
&=&
\int_{\frac12}^\infty \G(\hat{x}-r\omega)\frac{1}{r}\,dr +
\int_{\frac{1}{|x|}}^{\frac12}\G(\hat{x}-r\omega)\frac{1}{r}\,dr.
\end{eqnarray*}
Now for $\omega\neq\hat{x}$, define
\begin{equation}\label{HH1}
\HH_1(\hat{x},\omega):=\int_{\frac12}^\infty \G(\hat{x}-r\omega)\frac{1}{r}\,dr = O(1).
\end{equation}
Then we have
\begin{equation}\label{HH1.int}
\int_{\mathbb{S}^2}
\left[\int_{\frac12}^\infty\G(\hat{x}-r\omega)\frac{1}{r}\,dr\right]
\f(\omega)\,d\omega
=
\int_{\mathbb{S}^2}\HH_1(\hat{x}, \omega)\f(\omega)\,d\omega.
\end{equation}
On the other hand, we have 
$\displaystyle \int_{\frac{1}{|x|}}^{\frac12}\G(\hat{x}-r\omega)\frac{1}{r}\,dr
= O(\ln|x|)$. To analyze this term, we compute
\begin{eqnarray}
&&\int_{\mathbb{S}^2}
\left[\int_{\frac{1}{|x|}}^\frac12\G(\hat{x}-r\omega)\frac{1}{r}\,dr\right]
\f(\omega)\,d\omega\nonumber\\
&=&
\int_{\mathbb{S}^2}
\left[\int_{\frac{1}{|x|}}^\frac12\frac{\G(\hat{x})}{r}\,dr\right]
\f(\omega)\,d\omega
+
\int_{\mathbb{S}^2}
\left[\int_{\frac{1}{|x|}}^\frac12\frac{\G(\hat{x}-r\omega)-\G(\hat{x})}{r}\,dr\right]
\f(\omega)\,d\omega\nonumber\\
&=&
\G(\hat{x})\ln(\frac{|x|}{2})\int_{\mathbb{S}^2}\f(\omega)\,d\omega + 
\int_{\mathbb{S}^2} \HH_2(\hat{x},\omega)
\f(\omega)\,d\omega
-\int_{\mathbb{S}^2}
\left[\int_0^{\frac{1}{|x|}}\frac{\G(\hat{x}-r\omega)-\G(\hat{x})}{r}\,dr\right]
\f(\omega)\,d\omega\nonumber\\
\label{HH2.int}
\end{eqnarray}
where 
\begin{equation}
\HH_2(\hat{x},\omega) 
:= \int_0^\frac12
\frac{\G(\hat{x}-r\omega)-\G(\hat{x})}{r}\,dr=O(1)
\label{HH2}
\end{equation}
because 
$\displaystyle \frac{\G(\hat{x}-r\omega)-\G(\hat{x})}{r}=O(1)$
as $r\to0^+$. Lastly, we have
\begin{equation}\label{HH3.int}
\int_{\mathbb{S}^2}
\left[\int_0^{\frac{1}{|x|}}\frac{\G(\hat{x}-r\omega)-\G(\hat{x})}{r}\,dr\right]
\f(\omega)\,d\omega=O(\frac{1}{|x|}).
\end{equation}

Combining \eqref{HH1}, \eqref{HH2.int}, \eqref{HH2}, we define
\begin{eqnarray}
\HH(\hat{x},\omega)
&:=&\HH_1(\hat{x},\omega)+\HH_2(\hat{x},\omega) - \G(\hat{x})\ln(2)\nonumber\\
&=& \int_0^\frac12
\frac{\G(\hat{x}-r\omega)-\G(\hat{x})}{r}\,dr
+\int_{\frac12}^\infty \frac{\G(\hat{x}-r\omega)}{r}\,dr
- \G(\hat{x})\ln(2).
\label{HH}
\end{eqnarray}
Note that the above does not depend on the choice of $\frac12$. Hence it
can also be equivalently written as
\begin{equation}\label{HH.formG}
\HH(\hat{x},\omega)
= \lim_{\epsilon\to0}\left[
\int_{\epsilon}^\infty \frac{\G(\hat{x}-r\omega)}{r}\,dr
+\G(\hat{x})\ln(\epsilon)
\right].
\end{equation}
Then we have for $|x|\gg 1$ that,
\begin{eqnarray}
\int_{\Omega} \G(x-y)\f(y)\,dy
=
\G(x)\ln|x|\int_{\mathbb{S}^2}\f(\omega)\,d\omega
+\frac{1}{|x|}\int_{\mathbb{S}^2} \HH(\hat{x},\omega) \f(\omega)\,d\omega
+ O\left(\frac{1}{|x|^2}\right).
\label{bulk.log}
\end{eqnarray}
Hence if $\displaystyle \int_{\mathbb{S}^2}\f(\omega)\,d\omega=0$, then
\begin{equation}
\int_{\Omega} \G(x-y)\f(y)\,dy
=\frac{1}{|x|}\int_{\mathbb{S}^2} \HH(\hat{x},\omega) \f(\omega)\,d\omega
+ O\left(\frac{1}{|x|^2}\right).
\label{bulk.no.log}
\end{equation}

\section{$C^{2,\alpha}$ estimates for bulk integral}\label{PotEst}
Here we prove the decay estimates for the bulk integral in 
\eqref{bulk.int.soln}
\[
\int_\Omega \G(x-y)\f(y)\,dy
\]
with $\G=\G_\gamma$ and $\f=\f_\gamma$ given by \eqref{f:explicit}.
These are needed in step {\bf (II)} of the existence proof of solution 
in Section \ref{Exist.Soln}.

As mentioned before, the key of the proof is based on the estimates of 
Newtonian potentials \cite[Chapter 4]{GT} as $\G$ is homogeneous with degree $-1$. 
For estimates near the boundary, we can refer to \eqref{Bdry.Holder}. 
Hence we will just concentrate here
on interior weighted estimates at $|\bx| > a+1$. Only the proof involving the term $\displaystyle \frac{\bF(\hat{\br})}{r^3}$ in $\f$ will be given as it is the dominant term in terms of spatial decay property. 
All the other terms can be handled similarly. The proof
is divided into $C^1$, $C^2$ and $C^{2,\alpha}$ estimates. 

\subsection{$C^1$ estimates}
These are relatively simple as $D_\bx\G(\bx)$ is integrable near $\bx=0$.
Since
$\displaystyle \int_{\mathbb{S}^2}\bF(\hat{\br})\,d\sigma = 0$, 
then as shown in Appendix \ref{AppFFZ}, we have
\begin{eqnarray*}
\left|\int_\Omega 
\G(\bx-\by)\frac{\bF(\hat{\by})}{|\by|^3} \,dy\right|
\lesssim \frac{1}{r}.
\end{eqnarray*}
Similarly, we have by the mean zero condition again that
\begin{eqnarray*}
\left|D_\bx\int_\Omega\G(\bx-\by)\frac{\bF(\hat{\by})}{|\by|^3}\,dy
\right|
\lesssim \frac{1}{r^2}.
\end{eqnarray*}
Hence
\[
\sup_{\bx\in\Omega}\left\{
|\bx|
\left|\int_\Omega 
\G(\bx-\by)\frac{\bF(\hat{\by})}{|\by|^3}\,dy\right|,
\quad
|\bx|^2
\left|D_{\bx}\int_\Omega 
\G(\bx-\by)\frac{\bF(\hat{\by})}{|\by|^3}\,dy\right|
\right\}
\lesssim 1.
\]

We remark that the mean zero condition is preserved during
each iteration as the function $\bF(\hat{\by})$ does not depend on 
the solution ${\bu}$.

\subsection{$C^2$ estimates}
We now proceed to the estimates for $D^2_\bx$ of the bulk integrals.
One difficulty is that now the singularity of 
$D_\bx^2 \G(\bx)\sim\frac{1}{|\bx|^3}$ is not integrable.
We follow the approach of \cite[Chapter 4, Lemma 4.2]{GT}.

Let $R=2|x|$. Consider 
\begin{eqnarray*}
&& D_\bx^2\int_\Omega
\G(\bx-\by)\frac{\bF(\hat{\by})}{|\by|^3}\,dy\\
&=& D_\bx\int_\Omega
D_\bx\G(\bx-\by)\frac{\bF(\hat{\by})}{|\by|^3} \,dy\\
&=& 
\lim_{\bh\to0}\int_\Omega
\left(\frac{D_\bx\G(\bx+\bh-\by) -D_\bx\G(\bx-\by)}{|\bh|}\right)\frac{\bF(\hat{\by})}{|\by|^3}
\,dy\\
& = & 
\lim_{\bh\to0}\int_{\R^3\backslash\bB_R}
\left(\frac{D_\bx\G(\bx+\bh-\by) - D_\bx\G(\bx-\by)}{|\bh|}\right)\frac{\bF(\hat{\by})}{|\by|^3}
\,dy\\
& & +
\lim_{\bh\to0}\int_{\bB_R\backslash\bB_a}
\left(\frac{D_\bx\G(\bx+\bh-\by) - D_\bx\G(\bx-\by)}{|\bh|}\right)\frac{\bF(\hat{\by})}{|\by|^3}
\,dy\\
&:=& A + B
\end{eqnarray*}

For $A$, we have
\begin{eqnarray*}
|I| &=&\left|\lim_{\bh\to0}\int_{\R^3\backslash\bB_R}
\left(\frac{D_\bx\G(\bx+\bh-\by) - D_\bx\G(\bx-\by)}{|\bh|}\right)
\frac{\bF(\hat{\by})}{|\by|^3}
\,dy\right|\\
&=&\left|\int_{\R^3\backslash\bB_R}
D^2_\bx\G(\bx-\by)\frac{\bF(\hat{\by})}{|\by|^3} \cdot
\,dy\right|\\
& \lesssim & 
\int_{\R^3\backslash\bB_R}
\frac{1}{|\by|^3}
\frac{1}{|\bx-\by|^3}
\,dy
\leq
\frac{1}{|\bx|^3}\int_{\R^3\backslash\bB_2}
\frac{1}{|\bz|^3}
\frac{1}{|\hat{\bx}-\bz|^3}
\,dz
\lesssim O\left(\frac{1}{|\bx|^3}\right).
\end{eqnarray*}

For $B$, we have
\begin{eqnarray*}
B & = &
\lim_{\bh\to0}
\int_{\bB_R\backslash\bB_a}
\left(\frac{D_\bx\G(\bx+\bh-\by) - D_\bx\G(\bx-\by)}{|\bh|}\right)
\,dy\frac{\bF(\hat{\bx})}{|\bx|^3}\\
& &+
\lim_{\bh\to0}\int_{\bB_R\backslash\bB_a}
\left(\frac{D_\bx\G(\bx+\bh-\by) - D_\bx\G(\bx-\by)}{|\bh|}\right)
\left(
\frac{\bF(\hat{\by})}{|\by|^3}
-\frac{\bF(\hat{\bx})}{|\bx|^3}
\right)
\,dy\\
& := & B_1 + B_2.
\end{eqnarray*}
For $B_1$, by Divergence Theorem, we have
\begin{eqnarray*}
B_1 & = & 
\int_{\partial(\bB_R\backslash\bB_a)}
D_\bx\G(\bx-\by)\cdot\hat{\bh}
\,d\sigma(y)\frac{\bF(\hat{\bx})}{|\bx|^3}
\end{eqnarray*}
and hence
\begin{eqnarray*}
|B_1| & \lesssim & 
\frac{1}{|\bx|^3}\int_{\partial\bB_a}\frac{1}{R^2}\,d\sigma(y)
+
\frac{1}{|\bx|^3}\int_{\partial\bB_R}\frac{1}{R^2}\,d\sigma(y)
\lesssim O\left(\frac{1}{|\bx|^3}\right).
\end{eqnarray*}
For $B_2$, we have
\begin{eqnarray*}
B_2 &=&
\int_{\bB_R\backslash\bB_a}
D^2_\bx\G(\bx-\by)
\left(
\frac{\bF(\hat{\by})}{|\by|^3}
-\frac{\bF(\hat{\bx})}{|\bx|^3}
\right)
\,dy\\
&=&
\frac{1}{|\bx|^3}\int_{\bB_{\frac{R}{|\bx|}}\backslash\bB_{{\frac{a}{|\bx|}}}}
D^2_\bx\G(\hat{\bx}-\bz)\left(\frac{\bF(\hat{\bz})}{|\bz|^3}
- \bF(\hat{\bx})
\right) 
\,dz\\
&=&
\frac{1}{|\bx|^3}\int_{\bB_2\backslash\bB_{{\frac{a}{|\bx|}}}}
D^2_\bx\G(\hat{\bx}-\bz)
\left(
\frac{\bF(\hat{\bz})- |\bz|^3\bF(\hat{\bx})
}{|\bz|^3}
\right)
\,dz.
\end{eqnarray*}
Note that the singularity at $\bz=\hat{\bx}$ is now integrable due to the presence of $\bF(\hat{\bz}) - |\bz|^3\bF(\hat{\bx})$
in the numerator. On the other hand, the singularity behavior of
$\frac{1}{|\bz|^3}$ near $|\bz| \sim \frac{a}{|\bx|} \ll 1$ is 
tempered by the mean zero condition again: 
\[
\lim_{|\bz|\to0}\int_{\mathbb S^2}
\Big[\bF(\hat{\bz}) - |\bz|^3\bF(\hat{\bx})\Big]
\,d\sigma(\bz)  = 0.
\] 
Then similar to the computation leading to \eqref{bulk.no.log}, we have
\[
|B_2| \lesssim O\left(\frac{1}{|\bx|^3}\right).
\]

Combing the estimates for $A$ and $B$ gives,
\[
\sup_{\bx\in\Omega}\left\{
|\bx|^3
\left|D^2_{\bx}\int_\Omega 
\G(\bx-\by)\frac{\bF(\hat{\by})}{|\by|^3} \,dy\right|
\right\}
\lesssim 1.
\]

\subsection{$C^{2,\alpha}$ estimates}
For the $C^\alpha$ estimates for the second derivatives of the bulk 
integrals, we follow \cite[Chapter 4, Lemma 4.4]{GT}.

Using the computation for $\bF$ from the previous section and letting again $R = 2|\bx|$,
we have
\begin{eqnarray*}
&& D_\bx^2\int_\Omega
\G(\bx-\by)\frac{\bF(\hat{\by})}{|\by|^3} \,dy\\
&=&
\int_{\R^3\backslash\bB_R}
D^2_\bx\G(\bx-\by)\frac{\bF(\hat{\by})}{|\by|^3}
\,dy
+  \int_{\partial(\bB_R\backslash\bB_a)}
D_\bx\G(\bx-\by)\cdot\hat{\bh}
\,d\sigma(y)\frac{\bF(\hat{\bx})}{|\bx|^3}\\
& & + \int_{\bB_R\backslash\bB_a}
D^2_\bx\G(\bx-\by)\left(
\frac{\bF(\hat{\by})}{|\by|^3}
-\frac{\bF(\hat{\bx})}{|\bx|^3}
\right)
\,dy\\
& =: & C_1(x) + C_2(x) + C_3(x).
\end{eqnarray*}

For $C_1$, in fact, it is differentiable in $\bx$:
\begin{eqnarray*}
\left|\lim_{|\bh|\to0}\frac{C_1(\bx+\bh) - C_1(\bx)}{|\bh|}\right|
& = & \left|\int_{\R^3\backslash\bB_R}
D^3_\bx\G(\bx-\by)
\frac{\bF(\hat{\by})}{|\by|^3}
\,dy\right|\\
& \lesssim & 
\int_{\R^3\backslash\bB_R} \frac{1}{|\by|^3|\bx-\by|^4}\,d^3y
\lesssim 
\frac{1}{|\bx|^4}\int_{\R^3\backslash\bB_2} \frac{1}{|\bz|^3|\hat{\bx}-\bz|^4}\,d^3z.
\end{eqnarray*}
Hence
\begin{equation}
|\bx|^4[C_1]_\alpha(\bx) \lesssim O(1).
\end{equation}

For $C_2$, consider the following estimates,
\begin{eqnarray*}
&&\left|\lim_{|\bh|\to0}\frac{1}{|\bh|}\int_{\partial\bB_a}
\big(D_\bx\G(\bx+\bh-\by) -D_\bx\G(\bx-\by)\big)
\,d\sigma(y)\right|\\
& = &
\left|\int_{\partial\bB_a}
D^2_\bx\G(\bx-\by)\cdot\hat{\bh}
\,d\sigma(y)\right|
\lesssim
\frac{1}{|\bx|^3}
\end{eqnarray*}
\begin{eqnarray*}
&&\left|\lim_{|\bh|\to0}\frac{1}{|\bh|}\int_{\partial\bB_R}
\big(D_\bx\G(\bx+\bh-\by) -D_\bx\G(\bx-\by)\big)
\,d\sigma(y)\right|\\
& = & 
\left|\int_{\partial\bB_R}
D^2_\bx\G(\bx-\by)\cdot\hat{\bh}
\,d\sigma(y)\right|
\lesssim
\frac{1}{|\bx|}
\end{eqnarray*}
\[
\left|\lim_{|\bh|\to0}\frac{1}{|\bh|}\left(
\frac{1}{|\bx+\bh|^3}
-\frac{1}{|\bx|^3}
\right)\right| \lesssim \frac{1}{|\bx|^4}
\]
\[
\left|\bF(\widehat{\bx+\bh}) - \bF(\widehat{\bx})\right|
\lesssim \left|\widehat{\bx+\bh} - \widehat{\bx}\right|^\alpha
= \left|\frac{\bx+\bh}{|\bx+\bh|} - \widehat{\bx}\right|^\alpha
= \left|\frac{\widehat{\bx}+\frac{\bh}{|\bx|}}
{\left|\widehat{\bx}+\frac{\bh}{|\bx|}\right|} - \widehat{\bx}\right|^\alpha
\lesssim \left(\frac{|\bh|}{|\bx|}\right)^\alpha
\]
\[
\left|
\int_{\partial(\bB_R\backslash\bB_a)}
D_\bx\G(\bx-\by)\cdot\hat{\bh}
\,d\sigma(y)
\right| \lesssim 1.
\]
Then we have
\begin{eqnarray}
|C_2(\bx+\bh) - C_2(\bx)|
\lesssim 
\frac{1}{|\bx|^3}\frac{|\bh|}{|\bx|}
+\frac{|\bh|}{|\bx|^4}
+\frac{|\bh|^\alpha}{|\bx|^\alpha}\frac{1}{|\bx|^3}
\lesssim 
\frac{|\bh|^\alpha}{|\bx|^{3+\alpha}}.
\end{eqnarray}

For $C_3$, consider
\begin{eqnarray*}
& & C_3(\bx+\bh) - C_3(\bx)\\
& = & 
\int_{\bB_R\backslash\bB_a}
D^2_\bx\G(\bx+\bh-\by)\left(
\frac{\bF(\hat{\by})}{|\by|^3}
-\frac{\bF(\widehat{\bx+\bh})}{|\bx+\bh|^3}
\right)
\,dy\\
& & -
\int_{\bB_R\backslash\bB_a}
D^2_\bx\G(\bx-\by)\left(
\frac{\bF(\hat{\by})}{|\by|^3}
-\frac{\bF(\widehat{\bx})}{|\bx|^3}
\right)
\,dy\\
&=&
\int_{\bB_\delta(\xi)}
D^2_\bx\G(\bx+\bh-\by)\left(
\frac{\bF(\hat{\by})}{|\by|^3}
-\frac{\bF(\widehat{\bx+\bh})}{|\bx+\bh|^3}
\right)
\,dy\\
& & - \int_{\bB_\delta(\xi)}
D^2_\bx\G(\bx-\by)\left(
\frac{\bF(\hat{\by})}{|\by|^3}
-\frac{\bF(\widehat{\bx})}{|\bx|^3}
\right)
\,dy\\
& & +
\int_{\left(\bB_R\backslash\bB_a\right)\backslash \bB_\delta(\xi)}
D^2_\bx\G(\bx+\bh-\by)
\,dy\left(\frac{\bF(\widehat{\bx})}{|\bx|^3}
-
\frac{\bF(\widehat{\bx+\bh})}{|\bx+\bh|^3}\right)\\
& & +
\int_{\left(\bB_R\backslash\bB_a\right)\backslash \bB_\delta(\xi)}
\left(D^2_\bx\G(\bx+\bh-\by)
-D^2_\bx\G(\bx-\by)\right)\left(
\frac{\bF(\hat{\by})}{|\by|^3}
-\frac{\bF(\widehat{\bx})}{|\bx|^3}
\right)
\,dy
\end{eqnarray*}
where $\xi = \bx + \frac12\bh$ and $\delta = |\bh|$. We estimate,
\begin{eqnarray*}
& & \left|\int_{\bB_\delta(\xi)}
D^2_\bx\G(\bx+\bh-\by)\left(
\frac{\bF(\hat{\by})}{|\by|^3}
-\frac{\bF(\widehat{\bx+\bh})}{|\bx+\bh|^3}
\right)
\,dy\right|\\
& \lesssim & 
\int_{\bB_\delta(\xi)}
\left(\frac{|\by-\bx|}{|\bx|^4}
+ \frac{1}{|\bx|^3}\frac{|\by-\bx|^\alpha}{|\bx|^\alpha}\right)
\frac{1}{|\by-\bx|^3}\,d^3y\\
& \lesssim & 
\frac{1}{|\bx|^{3+\alpha}}
\int_{\bB_\delta(\xi)}\frac{|\by-\bx|^\alpha}{|\by-\bx|^3}\,d^3y
\lesssim \frac{|\bh|^\alpha}{|\bx|^{3+\alpha}},
\end{eqnarray*}
and similarly,
\[
\left|\int_{\bB_\delta(\xi)}
D^2_\bx\G(\bx-\by)\left(\frac{\bF(\hat{\by})}{|\by|^3}
-\frac{\bF(\widehat{\bx})}{|\bx|^3}
\right)
\,dy\right|
\lesssim \frac{|\bh|^\alpha}{|\bx|^{3+\alpha}},
\]
Furthermore, we have by Divergence Theorem
\begin{eqnarray*}
\left|\int_{\left(\bB_R\backslash\bB_a\right)\backslash \bB_\delta(\xi)}
D^2_\bx\G(\bx+\bh-\by)
\,dy\right|
& = & 
\left|\int_{\partial\left(\left(\bB_R\backslash\bB_a\right)\backslash \bB_\delta(\xi)\right)}
D_\bx\G(\bx+\bh-\by)\cdot\nu
\,dy\right|\\
& \lesssim & 
\int_{\partial\bB_R} +\int_{\partial\bB_a}
+\int_{\partial\bB_\delta(\xi)}
|D_\bx\G(\bx+\bh-\by)|\,d\sigma(y)\\
& \lesssim & 1
\end{eqnarray*}
so that
\begin{eqnarray*}
\left|
\int_{\left(\bB_R\backslash\bB_a\right)\backslash \bB_\delta(\xi)}
D^2_\bx\G(\bx+\bh-\by)
\,dy\left(\frac{\bF(\widehat{\bx})}{|\bx|^3}
- \frac{\bF(\widehat{\bx+\bh})}{|\bx+\bh|^3}\right)\right|
\lesssim 
\frac{|\bh|^\alpha}{|\bx|^{3+\alpha}}.
\end{eqnarray*}
For the last term,
\begin{eqnarray*}
&&\int_{\left(\bB_R\backslash\bB_a\right)\backslash \bB_\delta(\xi)}
\left(D^2_\bx\G(\bx+\bh-\by)
-D^2_\bx\G(\bx-\by)\right)\left(
\frac{\bF(\hat{\by})}{|\by|^3}
-\frac{\bF(\widehat{\bx})}{|\bx|^3}
\right)
\,dy\\
&=&\int_{\bB_2(\xi)\backslash \bB_\delta(\xi)}
\left(D^2_\bx\G(\bx+\bh-\by)
-D^2_\bx\G(\bx-\by)\right)
\left(
\frac{\bF(\hat{\by})}{|\by|^3}
-\frac{\bF(\widehat{\bx})}{|\bx|^3}
\right)
\,dy\\
&&+\int_{\left(\bB_R\backslash\bB_a\right)\backslash \bB_2(\xi)}
\left(D^2_\bx\G(\bx+\bh-\by)
-D^2_\bx\G(\bx-\by)\right)\left(
\frac{\bF(\hat{\by})}{|\by|^3}
-\frac{\bF(\widehat{\bx})}{|\bx|^3}
\right)
\,dy\\
&=:& D_1 + D_2.
\end{eqnarray*}
For $D_1$, we estimate
\begin{eqnarray*}
|D_1|
&=& |\bh|
\left|\int_{\bB_2(\xi)\backslash \bB_\delta(\xi)}
\frac{(D^2_\bx\G(\bx+\bh-\by)
-D^2_\bx\G(\bx-\by)}{|\bh|}\left(
\frac{\bF(\hat{\by})}{|\by|^3}
-\frac{\bF(\widehat{\bx})}{|\bx|^3}
\right)
\,dy\right|\\
& \lesssim & 
|\bh|
\left|\int_{\bB_2(\xi)\backslash \bB_\delta(\xi)}
D^3_\bx\G(\tilde{\bx}+\bh-\by)\left(
\frac{\bF(\hat{\by})}{|\by|^3}
-\frac{\bF(\widehat{\bx})}{|\bx|^3}
\right)
\,dy\right|\\
& & \text{(where $\tilde{\bx}$ is some point between $\bx$ and $\bx+\bh$)}\\
&\lesssim &
\frac{|\bh|}{|\bx|^{3+\alpha}}
\int_{\bB_2(\xi)\backslash \bB_\delta(\xi)}
\frac{|\bx-\by|^\alpha}{|\tilde{\bx}-\by|^4}\,dy
\lesssim
\frac{|\bh|^\alpha}{|\bx|^{3+\alpha}}.
\end{eqnarray*}
For $J_2$, we have
\begin{eqnarray*}
D_2 
& = & 
|\bh|\int_{\left(\bB_R\backslash\bB_a\right)\backslash \bB_2(\xi)}
\cdot
\frac{D^2_\bx\G(\bx+\bh-\by)
-D^2_\bx\G(\bx-\by)}{|\bh|}\left(
\frac{\bF(\hat{\by})}{|\by|^3}
-\frac{\bF(\widehat{\bx})}{|\bx|^3}
\right)
\,dy\\
& \approx & 
|\bh|\int_{\left(\bB_R\backslash\bB_a\right)\backslash \bB_2(\xi)}
D^3_\bx\G(\tilde{\bx}-\by)\left(
\frac{\bF(\hat{\by})}{|\by|^3}
-\frac{\bF(\widehat{\bx})}{|\bx|^3}
\right)
\,dy\\
& = &
\frac{|\bh|}{|\bx|^4}
\int_{\left(\bB_{\frac{R}{|\bx|}}\backslash\bB_{\frac{a}{|\bx|}}\right)
\big\backslash \frac{\bB_2(\xi)}{|\bx|}}
D^3_\bx\G(\hat{\bx}-\bz)\left(
\frac{\bF(\hat{\bz})}{|\bz|^3}
- \bF(\widehat{\bx})
\right)
\,dz\\
& = &
\frac{|\bh|}{|\bx|^4}
\int_{
\bB_\frac12(\frac{\xi}{|\bx|})
\backslash
\bB_{\frac{2}{|\bx|}}(\frac{\xi}{|\bx|})
}
D^3_\bx\G(\hat{\bx}-\bz)\left(
\frac{\bF(\hat{\bz})}{|\bz|^3}
- \bF(\widehat{\bx})
\right)
\,dz\\
& & +
\frac{|\bh|}{|\bx|^4}
\int_{
\left(\bB_{\frac{R}{|\bx|}}\backslash\bB_{\frac{a}{|\bx|}}\right)
\backslash
\bB_\frac12(\frac{\xi}{|\bx|})
}
D^3_\bx\G(\hat{\bx}-\bz)\left(
\frac{\bF(\hat{\bz})}{|\bz|^3}
- \bF(\widehat{\bx})
\right)
\,dz\\
& := & 
D_{21} + D_{22}
\end{eqnarray*}
For $D_{21}$,
\begin{eqnarray*}
|D_{21}| 
& \lesssim & 
\frac{|\bh|}{|\bx|^4}
\int_{
\bB_\frac12(\frac{\xi}{|\bx|})
\backslash
\bB_{\frac{2}{|\bx|}}(\frac{\xi}{|\bx|})
}\frac{|\bz - \hat{\bx}|^\alpha}{|\bz - \hat{\bx}|^4}\,dz
\lesssim 
\frac{|\bh|}{|\bx|^4}|\bx|^{1-\alpha}
=
\frac{|\bh|}{|\bx|^{3+\alpha}}
\end{eqnarray*}
For $D_{22}$, we again make use of the mean zero condition for $\bF$ to 
have that
\[
|D_{22}| \lesssim \frac{|\bh|}{|\bx|^4}.
\]
Combining the above, we conclude that
\begin{equation}
|\bx|^{3+\alpha}\left[D_\bx^2\int_\Omega
\G(\bx-\by)\frac{\bF(\hat{\by})}{|\by|^3} \,dy\right]_{\alpha}
\lesssim 1.
\end{equation}

\section{Calculation of
$\M_\gamma:D^2\E,\,\,\,\M_\gamma:D^2\F,\,\,\,\div\cA_\gamma,\,\,\,\cC_\gamma$ and 
$\cD_\gamma$}
Before proceeding, we recall the conventions
\eqref{notation.lin.mult} and \eqref{notation.lin.comb}.

\subsection{Formula for $\M_\gamma:D^2\bF$}\label{MD2EF}
Here we list the formulas for $\M_\gamma:D^2\bF$ where $\bF$ is some 
divergence free vector field with homogeneous degree $-k$, for
$k=1,3,5$. 

As a preparation, following \eqref{DDgr} for
$g=g(r) = \frac{1}{r},\,\frac{1}{r^3},\,\frac{1}{r^5}$, we compute,
\begin{eqnarray}
\partial_{k}g(r) & = & g'(r)\frac{x_k}{r},\\
\partial_{kl}g(r) &=&
\frac{g'(r)}{r}\delta_{kl}
+ \left(g''(r)-\frac{g'(r)}{r}\right)\frac{x_kx_l}{r^2},\\
\partial_{kl}\big(g(r)\delta_{ij}\big) &=&
\left(\frac{g'(r)}{r}\delta_{kl} 
+ \left(g''(r)-\frac{g'(r)}{r}\right)\frac{x_kx_l}{r^2}\right)\delta_{ij},\\
\partial_{k}\big(g(r)x_ix_j\big)
& = & g(r)\big(\delta_{ik}x_j + \delta_{jk}x_i\big)
+\frac{g'(r)}{r}x_ix_jx_k,\\
\partial_{kl}\big(g(r)x_ix_j\big)
&=&
g(r)\big(\delta_{ik}\delta_{jl} + \delta_{jk}\delta_{il}\big)\nonumber\\
& & +rg'(r)\Big(\frac{
\delta_{il}x_jx_k + \delta_{jl}x_ix_k + \delta_{kl}x_ix_j
+ \delta_{ik}x_jx_l + \delta_{jk}x_ix_l
}{r^2}
\Big)\nonumber\\
&&+ (r^2g''(r) - rg'(r))\frac{x_ix_jx_kx_l}{r^4}.
\end{eqnarray}
Hence,
\begin{eqnarray*}
\partial_{kl}\left(\frac{\delta_{ij}}{r}\right)
& = & \frac{1}{r^3}\big(-\delta_{kl} + 3\hx_k\hx_l\big)\delta_{ij},\\
\partial_{kl}\left(\frac{\delta_{ij}}{r^3}\right)
& = & \frac{1}{r^5}\big(-3\delta_{kl} + 15\hx_k\hx_l\big)\delta_{ij},\\
\partial_{kl}\left(\frac{x_ix_j}{r^3}\right)
&=&
\frac{1}{r^3}\Big[\big(\delta_{ik}\delta_{jl} + \delta_{jk}\delta_{il}\big)\nonumber\\
&& \hspace{20pt}-3\big(
\delta_{il}\hx_j\hx_k + \delta_{jl}\hx_i\hx_k + \delta_{kl}\hx_i\hx_j
+ \delta_{ik}\hx_j\hx_l + \delta_{jk}\hx_i\hx_l
\big)\nonumber\\
&& \hspace{20pt}+ 15\hx_i\hx_j\hx_k\hx_l\Big],\\
\partial_{kl}\left(\frac{x_ix_j}{r^5}\right)
&=&
\frac{1}{r^5}\Big[\big(\delta_{ik}\delta_{jl} + \delta_{jk}\delta_{il}\big)\nonumber\\
&& \hspace{20pt}-5\big(
\delta_{il}\hx_j\hx_k + \delta_{jl}\hx_i\hx_k + \delta_{kl}\hx_i\hx_j
+ \delta_{ik}\hx_j\hx_l + \delta_{jk}\hx_i\hx_l
\big)\nonumber\\
&& \hspace{20pt}+ 35\hx_i\hx_j\hx_k\hx_l\Big].
\end{eqnarray*}

The above are applied to $\E, \F$ and $\Q$.
\begin{itemize}
\item	
$\displaystyle \E(\br) = \frac{1}{8\pi}\left[
\frac{\id}{r} + \frac{\br\otimes\br}{r^3}
\right]$.
\begin{eqnarray}
\partial_{k}\E_{ij}(x) 
&=& \frac{1}{8\pi r^2}\Big[-\delta_{ij}\hat{x}_k + \delta_{ik}\hat{x}_j + 
\delta_{jk}\hat{x}_i
- 3\hat{x}_i\hat{x}_j\hat{x}_k\Big]
\label{DE}\\
\partial_{kl}\E_{ij}(x) 
& = &  
\frac{1}{8\pi r^3}\Big[
\left(-\delta_{ij}\delta_{kl} + \delta_{ik}\delta_{jl} 
+ \delta_{jk}\delta_{il}\right)\nonumber\\
&&
\hspace{40pt}-3\left(
-\delta_{ij}\hat{x}_k\hat{x}_l +\delta_{ik}\hat{x}_j\hat{x}_l 
+\delta_{jk}\hat{x}_i\hat{x}_l
+\delta_{il}\hat{x}_j\hat{x}_k +\delta_{jl}\hat{x}_i\hat{x}_k 
+\delta_{kl}\hat{x}_i\hat{x}_j
\right)\nonumber\\
&&\hspace{40pt} + 15\hat{x}_i\hat{x}_j\hat{x}_k\hat{x}_l\Big].
\label{DDE}
\end{eqnarray}

\item	$\displaystyle \F(x) = \frac{3x\otimes x - r^2\I}{r^5}$.
\begin{eqnarray}
\partial_{k}\F_{ij}(x)
& = & \frac{1}{r^4}\Big[
3\big(\delta_{ij}\hat{x}_k + \delta_{ik}\hat{x}_j + \delta_{jk}\hat{x}_i\big)
- 15\hat{x}_i\hat{x}_j\hat{x}_k
\Big]\\
\partial_{kl}\F_{ij}(x)
&=& 
\frac{1}{r^5}\Big[
3\big(\delta_{ik}\delta_{jl} + \delta_{jk}\delta_{il}
+ \delta_{kl}\delta_{ij}\big)\nonumber\\
&&\hspace{20pt}
-15\big(
\delta_{ij}\hx_k\hx_l
+\delta_{il}\hx_j\hx_k
+\delta_{jl}\hx_i\hx_k
+\delta_{kl}\hx_i\hx_j
+\delta_{ik}\hx_j\hx_l
+\delta_{jk}\hx_i\hx_l
\big)\nonumber\\
&&\hspace{20pt}
+ 105 \hx_i\hx_j\hx_k\hx_l
\Big]
\end{eqnarray}

\item $\displaystyle \Q(x) = 
\left(1 - \frac{w}{1+w}\frac{1}{r}\right)\Q_* +\frac{w}{3+w}\frac{1}{r^3}\Q_b
= 
\left(1 - \frac{w}{1+w}\frac{1}{r}\right)\Q_* +\frac{ws_*}{3(3+w)}\F(x)$,
\\
where $\displaystyle {\Q_*}_{ij} = \n_i\n_j - \frac13\delta_{ij}$.
\begin{eqnarray}
\partial_k\Q_{ij}(x)
& = & 
\frac{w}{(1+w)r^2}\hx_k{\Q_*}_{ij}\nonumber
\\&&+ \frac{ws_*}{(3+w)r^4}
\Big[
\delta_{ij}\hat{x}_k + \delta_{ik}\hat{x}_j + \delta_{jk}\hat{x}_i
- 5\hat{x}_i\hat{x}_j\hat{x}_k
\Big],
\\
\partial_{kl}\Q_{ij}(x)
& = & 
- \frac{w}{(1+w)r^3}\Big[-\delta_{kl} 
+ 3\hx_k\hx_l\Big]{\Q_*}_{ij}\nonumber\\
&&+
\frac{ws_*}{(3+w)r^5}\Big[
\big(\delta_{ik}\delta_{jl} + \delta_{jk}\delta_{il}
+ \delta_{kl}\delta_{ij}\big)\nonumber\\
&&\hspace{20pt}
-5\big(
\delta_{ij}\hx_k\hx_l
+\delta_{il}\hx_j\hx_k
+\delta_{jl}\hx_i\hx_k
+\delta_{kl}\hx_i\hx_j
+\delta_{ik}\hx_j\hx_l
+\delta_{jk}\hx_i\hx_l
\big)\nonumber\\
&&\hspace{20pt}
+ 35 \hx_i\hx_j\hx_k\hx_l
\Big].
\end{eqnarray}

\end{itemize}

We now recall the operation \eqref{uEqnM1}, \eqref{uEqnM2} and
\eqref{uEqnM3} for $\B_\gamma$:
\begin{eqnarray*}
\Big(\div\big(\B_\gamma[\nabla\bv]\big)\Big)_i  
= \Big(\M_\gamma:D^2\bv\Big)_i
= \sum_{p = 1,\ldots 7, 10, 11}\gamma_p\M^p_{i,j;k,l}\partial_{kl}\bv_j.
\end{eqnarray*}
Given any matrix valued function $\mathbf{F}(x) 
= \{\mathbf{F}_{ij}(x)\}_{1\leq i,j\leq 3}$, upon introducing the
contraction,
\begin{equation}\label{MD2EContract}
\big(\M^p:D^2{\mathbf F})_{ij} = \M^p_{i,m;,k,l}\partial_{kl}{\mathbf F}_{mj},
\end{equation}
then we have
\begin{eqnarray}
\big(\M^p:D^2(\bF\bv_*)\big)_i
= \big(\M^p_{i,m;k,l}\partial_{kl}\bF_{mj}\big){\bv_*}_j.
\end{eqnarray}
Note that we will only consider those ${\mathbf F}$ such that
${\mathbf F}\bv_*$ is divergence
free for all $\bv_*$, i.e.
\begin{equation}
\partial_i{\mathbf F}_{ik} = 0\,\,\,\text{for all $k$.}
\end{equation}
This condition is indeed satisfied by
$\bF = \E$ \eqref{GreenE}, $\E_S$ \eqref{StokesFormGenVec}, 
and hence $\F$ \eqref{GreenF}.

With the above, we proceed to find:
\begin{itemize}
\item{\bf Formula for $(\M^p:D^2\E)_{ij}(x) 
= \M^p_{i,m;k,l}\partial_{kl}\E_{mj}(x)$}.
From \eqref{MLinComb}, we need to compute the following:
\begin{eqnarray*}
&&\Big\{
\delta_{kl}\delta_{im},\,\,\,
\delta_{kl}\n_i\n_m,\,\,\,
\delta_{im}\n_k\n_l,\,\,\,
\delta_{ik}\n_l\n_m,\,\,\,
\n_i\n_m\n_k\n_l\Big\}\times \partial_{kl}\E_{mj}(x)\\
& = & 
\Big\{
\delta_{kl}\delta_{im},\,\,\,
\delta_{kl}\n_i\n_m,\,\,\,
\delta_{im}\n_k\n_l,\,\,\,
\delta_{ik}\n_l\n_m,\,\,\,
\n_i\n_m\n_k\n_l\Big\}\times \\
&&
\frac{1}{8\pi r^3}\Big[\left(-\delta_{mj}\delta_{kl} + \delta_{mk}\delta_{jl}
+ \delta_{jk}\delta_{ml}\right)\\
&&
\hspace{35pt}-3 \big(
-\delta_{mj}\hx_k\hx_l +\delta_{mk}\hx_j\hx_l +\delta_{jk}\hx_m\hx_l
+\delta_{ml}\hx_j\hx_k +\delta_{jl}\hx_m\hx_k +\delta_{kl}\hx_m\hx_j
\big)\\
&&\hspace{35pt} + 15\hx_m\hx_j\hx_k\hx_l\Big].
\end{eqnarray*}
We tabulate the result in the following:
\begin{eqnarray*}
\delta_{kl}\delta_{im}\partial_{kl}\E_{mj}(x)
& = & \frac{1}{8\pi r^3}\Big[2\delta_{ij} - 6\hx_i\hx_j\Big]\\
\delta_{kl}\n_i\n_m\partial_{kl}\E_{mj}(x)
& = & \frac{1}{8\pi r^3}\Big[
2\n_i\n_j - 6\n_i\hx_j\langle n, \hx\rangle\Big]\\
\delta_{im}\n_k\n_l\partial_{kl}\E_{mj}(x)
& = & 
\frac{1}{8\pi r^3}\Big[-\delta_{ij} + 2\n_i\n_j
+3\delta_{ij}\langle \n,\hx\rangle^2
-6(\n_i\hx_j + \n_j\hx_i)\langle \n,\hx\rangle
\nonumber\\
&&\hspace{35pt}
-3\hx_i\hx_j
+ 15\hx_i\hx_j\langle \n,\hx\rangle^2\Big]
\\
\delta_{ik}\n_l\n_m\partial_{kl}\E_{mj}(x)
& = &\frac{1}{8\pi r^3}\Big[
\delta_{ij} - 3\delta_{ij}\langle \n, \hx\rangle^2 
-6\n_ix_j\langle \n,\hx\rangle
-3\hx_i\hx_j + 15\hx_i\hx_j\langle a,\hx\rangle^2\Big]\\
\n_i\n_m\n_k\n_l\partial_{kl}\E_{mj}(x)
& = & 
\frac{1}{8\pi r^3}\Big[\n_i\n_j - 3\n_i\n_j\langle \n, \hx\rangle^2
-9 \n_i\hx_j\langle \n,\hx\rangle 
+ 15\n_i\hx_j\langle \n,\hx\rangle^3\Big].
\end{eqnarray*}

With the above, we have explicitly,
\begin{eqnarray*}
\M^1:D^2\E(x)&=&
\frac{s_*^2}{8\pi r^3}
\Big[
-\big(1 - 3{\langle \n,\hx\rangle^2}\big)\I
+ 2\n\otimes \n 
-3\langle \n,\hx\rangle
\big(\n\otimes \hx + \hx\otimes \n\big)
\Big]\\
\M^2:D^2\E(x)&=&
\frac{s_*}{8\pi r^3}
\Big[
\big(1 - 3\langle \n,\hx\rangle^2\big)\I 
+3\big(1-5\langle \n,\hx\rangle^2\big)\hx\otimes\hx
+ 6\langle \n,\hx\rangle \hx\otimes \n
\Big]\\
\M^3:D^2\E(x)&=&
\frac{1}{8\pi r^3}\frac{s_*^2}{6}
\Big[
\big(1 - 3\langle \n,\hx\rangle^2\big)\I 
+3\big(1-5\langle \n,\hx\rangle^2\big)\hx\otimes\hx
+ 6\langle \n,\hx\rangle \hx\otimes \n
\Big]\\
\M^4:D^2\E(x)&=&
\frac{s_*}{8\pi r^3}\Big[
-\frac23\I + 2\n\otimes \n
-\big(1 - 15\langle \n,\hx\rangle^2\big)\hx\otimes\hx
\nonumber\\
&&\hspace{35pt}
-9 \langle \n,\hx\rangle \n\otimes \hx
-3 \langle \n,\hx\rangle \hx\otimes \n
\Big]
\\
\M^5:D^2\E(x)&=&
\frac{1}{8\pi r^3}\frac{s_*^2}{6}\Big[
\frac43\I + 4\n\otimes \n
-10\big(1- 3\langle \n,\hx\rangle^2\big)\hx\otimes\hx
\nonumber\\
&&\hspace{45pt}
-18 \langle \n,\hx\rangle \n\otimes \hx
-6 \langle \n,\hx\rangle \hx\otimes \n
\Big]\\
\M^6:D^2\E(x)&=&
\frac{s_*^2}{8\pi r^3}\Big[
-\big(\frac13 -  \langle \n,\hx\rangle^2\big) \I
+ \big(1 - 3 \langle \n,\hx\rangle^2\big) \n\otimes \n
\nonumber\\
&&\hspace{35pt}
+ \big(1 -10 \langle \n,\hx\rangle^2\big) \hx\otimes\hx
-7 \langle \n,\hx\rangle \n\otimes\hx
+ 15 \langle \n,\hx\rangle^3 n\otimes\hx
\Big]\\
\M^7:D^2\E(x)&=&
\frac{s_*^2}{8\pi r^3}\Big[\frac23\I -2\hx\otimes\hx \Big]\\
\M^{10}:D^2\E(x)&=&
\frac{1}{8\pi r^3}\frac{2s_*^3}{3}\Big[
\big(1- 3 \langle \n,\hx\rangle^2\big) \n\otimes \n
-9 \langle \n,\hx\rangle \n\otimes\hx
+ 15 \langle \n,\hx\rangle^3 \n\otimes\hx
\Big]\\
\M^{11}:D^2\E(x)&=&
\frac{1}{8\pi r^3}\frac{2s_*^4}{3}
\Big[
-\big(\frac13 -  \langle \n,\hx\rangle^2\big) \I
+ \big(1 - 3 \langle \n,\hx\rangle^2\big) \n\otimes \n
\nonumber\\
&&\hspace{52pt}
+ \big(1 -10 \langle \n,\hx\rangle^2\big) \hx\otimes\hx
-7 \langle \n,\hx\rangle \n\otimes\hx
+ 15 \langle \n,\hx\rangle^3 \n\otimes\hx
\Big]
.
\end{eqnarray*}
We note that $\M^3 = \frac{s_*}{6}\M^2$, $\M^{11} = \frac{2s_*^2}{3}\M^6$.

\item{\bf Formula for $\M_\gamma:D^2\F(x)$}.
Similarly, we need to compute the following:
\begin{eqnarray*}
&& \Big\{
\delta_{kl}\delta_{im},\,\,\,
\delta_{kl}\n_i\n_m,\,\,\,
\delta_{im}\n_k\n_l,\,\,\,
\delta_{ik}\n_l\n_m,\,\,\,
\n_i\n_m\n_k\n_l\Big\}
\times \partial_{kl}\F_{mj}(x)\\
& = &  
\Big\{
\delta_{kl}\delta_{im},\,\,\,
\delta_{kl}\n_i\n_m,\,\,\,
\delta_{im}\n_k\n_l,\,\,\,
\delta_{ik}\n_l\n_m,\,\,\,
\n_i\n_m\n_k\n_l\Big\}\times \\
&&\frac{1}{r^5}\Big[
3\big(\delta_{mk}\delta_{jl} + \delta_{jk}\delta_{ml}
+ \delta_{kl}\delta_{mj}\big)\nonumber\\
&&\hspace{20pt}
-15\big(
\delta_{mj}\hx_k\hx_l
+\delta_{ml}\hx_j\hx_k
+\delta_{jl}\hx_m\hx_k
+\delta_{kl}\hx_m\hx_j
+\delta_{mk}\hx_j\hx_l
+\delta_{jk}\hx_m\hx_l
\big)\nonumber\\
&&\hspace{20pt}
+ 105 \hx_m\hx_j\hx_k\hx_l
\Big]
\end{eqnarray*}
We again tabulate the result in the following:
\begin{eqnarray*}
\delta_{kl}\delta_{im}\partial_{kl}\F_{mj}(x)
& = & 0\\
\delta_{kl}\n_i\n_m\partial_{kl}\F_{mj}(x)
& = & 0\\
\delta_{im}\n_k\n_l\partial_{kl}\F_{mj}(x)
& = &
\frac{1}{r^5}\Big[3\delta_{ij} 
-15\delta_{ij}\langle \n,\hx\rangle^2
+ 6\n_i\n_j -15\hx_i\hx_j
\nonumber\\
&&\hspace{20pt}
-30(\n_i\hx_j + \n_j\hx_i)\langle a,\hx\rangle
+ 105\hx_i\hx_j\langle \n,\hx\rangle^2\Big]
\\
\delta_{ik}\n_l\n_m\partial_{kl}\F_{mj}(x)
& = &\frac{1}{r^5}\Big[3\delta_{ij} - 15\delta_{ij}\langle \n, \hx\rangle^2
+6\n_i\n_j -15\hx_i\hx_j 
\nonumber\\
&&\hspace{20pt}
-30(\n_ix_j+ \n_j\hx_i)\langle a,\hx\rangle
+ 105\hx_i\hx_j\langle \n,\hx\rangle^2\Big]\\
\n_i\n_m\n_k\n_l\partial_{kl}\F_{mj}(x)
& = &
\frac{1}{r^5}\Big[9\n_i\n_j - 45\n_i\n_j\langle a, \hx\rangle^2
-45 \n_i\hx_j\langle \n,\hx\rangle
+ 105\n_i\hx_j\langle \n,\hx\rangle^3\Big].
\end{eqnarray*}

Similar to $\M_\gamma:D^2\E$, we have the following:
\begin{eqnarray*}
\M^1:D^2\F(x)&=& 0\\
\M^2:D^2\F(x)&=&
-\frac{s_*}{r^5}\Big[
3\big(1 - 5\langle \n,\hx\rangle^2\big)\I 
+ 6\n\otimes \n 
- 15(1-7\langle \n,\hx\rangle^2)
\hx\otimes\hx\\
&&\hspace{30pt}
-30\langle \n,\hx\rangle\big(\n\otimes\hx + \hx\otimes \n\big)
\Big]\\
\M^3:D^2\F(x)&=&
-\frac{s_*^2}{6r^5}\Big[
3\big(1 - 5\langle \n,\hx\rangle^2\big)\I 
+ 6\n\otimes \n 
- 15(1-7\langle \n,\hx\rangle^2)
\hx\otimes\hx\\
&&\hspace{35pt}
-30\langle \n,\hx\rangle\big(\n\otimes\hx + \hx\otimes \n\big)
\Big]\\
\M^4:D^2\F(x)&=&
\frac{s_*}{r^5}
\Big[
3\big(1 - 5\langle \n,\hx\rangle^2\big)\I
+ 6\n\otimes \n
- 15(1-7\langle n,\hx\rangle^2)
\hx\otimes\hx\\
&&\hspace{35pt}
-30\langle \n,\hx\rangle\big(\n\otimes\hx + \hx\otimes \n\big)
\Big]
\\
\M^5:D^2\F(x)&=&
\frac{s_*^2}{3r^5}
\Big[
3\big(1 - 5\langle \n,\hx\rangle^2\big)\I
+ 6\n\otimes \n
- 15(1-7\langle \n,\hx\rangle^2)
\hx\otimes\hx\\
&&\hspace{35pt}
-30\langle \n,\hx\rangle\big(\n\otimes\hx + \hx\otimes \n\big)
\Big]
\\
\M^6:D^2\F(x)&=&
\frac{s_*^2}{r^5}\Big[
-\big(1-5\langle \n,\hx\rangle^2\big)\I
+ \big(7 - 45\langle \n,\hx\rangle^2\big)\n\otimes \n
+ 5\big(1-7\langle \n,\hx\rangle^2\big)\hx\otimes\hx
\\
&&\hspace{40pt}
+35\big(-\langle \n,\hx\rangle+3\langle \n,\hx\rangle^3\big)\n\otimes\hx
+10\langle \n,\hx\rangle\hx\otimes \n
\Big]
\\
\M^7:D^2\F(x)&=&0\\
\M^{10}:D^2\F(x)&=&
\frac{2s_*^3}{3r^5}\Big[
9\big(1 - 5\langle \n,\hx\rangle^2\big)\n\otimes \n
+\big(-45\langle \n,\hx\rangle+105\langle \n,\hx\rangle^3\big)\n\otimes\hx
\Big]
\\
\M^{11}:D^2\F&=&
\frac{2s_*^4}{3r^5}
\Big[
-\big(1-5\langle \n,\hx\rangle^2\big)\I
+ \big(7 - 45\langle \n,\hx\rangle^2\big)\n\otimes \n
+ 5\big(1-7\langle \n,\hx\rangle^2\big)\hx\otimes\hx
\\
&&\hspace{40pt}
+35\big(-\langle \n,\hx\rangle+3\langle \n,\hx\rangle^3\big)\n\otimes\hx
+10\langle \n,\hx\rangle\hx\otimes \n
\Big].
\end{eqnarray*}
\end{itemize}

From the above, we can conclude that
$\M:D^2\E$ and $\M:D^2\F$ are
linear combinations of the following matrices
\[
\I,\,\,\, \n\otimes \n,\,\,\, \hx\otimes\hx,\,\,\, \n\otimes\hx,
\,\,\, \hx\otimes \n
\]
with coefficients given by $\langle n, \hx\rangle^k$ for
$k=0,1,2,3$. Concisely, we can write,
\begin{eqnarray}
&&\M_\gamma:D^2\E(x)\\
&=& \frac{1}{r^3}\Big[
f_1(\langle \n, \hx\rangle)\I
+f_2(\langle \n, \hx\rangle) \n\otimes \n
+f_3(\langle \n, \hx\rangle) \hx\otimes\hx
+f_4(\langle \n, \hx\rangle) \n\otimes\hx,
+f_5(\langle \n, \hx\rangle)\hx\otimes \n
\Big]\nonumber
\end{eqnarray}
and
\begin{eqnarray}
&&\M_\gamma:D^2\F(x)\\
&=& \frac{1}{r^5}\Big[
g_1(\langle \n, \hx\rangle)\I
+g_2(\langle \n, \hx\rangle) \n\otimes \n
+g_3(\langle \n, \hx\rangle) \hx\otimes\hx
+g_4(\langle \n, \hx\rangle) \n\otimes\hx,
+g_5(\langle \n, \hx\rangle)\hx\otimes \n
\Big]\nonumber
\end{eqnarray}
where the $f_i$ and $g_i$'s are polynomials of degree at most three.

\subsection{Formula for $\div\cA_\gamma$}\label{CAgammaForm}
By \eqref{divAform0}, \eqref{divAform} and \eqref{DDgr}, we have
for some matrix ${\bf M}_* = (m_{ij})$ that
\[
\Big(\div(\cA_\gamma)\Big)_i
= \div\left(\langle\bv_*, \nabla\frac{1}{r}\rangle {\bf M}_*\right)_i
= -m_{ij}(\delta_{kj} - 3\hx_k\hx_j){\bv_*}_k
= -m_{ij}{\bv_*}_j + 3m_{ij}\hx_j\langle\bv_*, \hx\rangle
\]
so that
\[
\div(\cA_\gamma) = -{\bf M}_*(\I - 3\hx\otimes\hx)\bv_*.
\]
Hence
\[
\div\left(\langle\bv_*, \nabla\frac{1}{r}\rangle {\bf M}_*\right)
= \left\{
\begin{array}{ll}
\big(-\langle\bv_*, \n\rangle + 
3\langle \n, \hx\rangle\langle\bv_*, \hx\rangle\big) \n,
& \text{if}\,\,\,{\bf M}_* = \n\otimes \n,\\
-\bv_* + 3\langle\bv_*, \hx\rangle\hx,
& \text{if}\,\,\,{\bf M}_* = \I.
\end{array}
\right.
\]
Using the form of $\Q_*$ from \eqref{QABL2}, 
we have
\[
{\bf M}_* 
= -\frac{w}{1+w}\left[
\left(\gamma_2s_* + \frac{\gamma_3}{3}s_*^2 + \frac{2\gamma_9}{3} s_*^3
\right)\n\otimes \n  
- \frac13\left(
\gamma_2s_*
-\frac{\gamma_3s_*^2}{3}
+ \frac{2\gamma_9s_*^3}{3}
\right)\I
\right].
\]
Hence
\begin{eqnarray}
\div(\cA_\gamma) 
&=& -\frac{w}{1+w}\left[
\left(\gamma_2s_* + \frac{\gamma_3}{3}s_*^2 + \frac{2\gamma_9}{3}s_*^3
\right)\big(
-\langle\bv_*, \n\rangle + 3\langle \n,\hx\rangle\langle\bv_*, \hx\rangle\big)\n  
\right.\nonumber\\ 
&&\hspace{50pt}\left.
- \frac13\left(
\gamma_2s_*
-\frac{\gamma_3s_*^2}{3}
+ \frac{2\gamma_9s_*^3}{3}
\right)
\big(
-\bv_* + 3\langle\bv_*, \hx\rangle\hx
\big)
\right]\nonumber\\
&=& -\frac{w}{1+w}\left[
\left(\gamma_2s_* + \frac{\gamma_3}{3}s_*^2 + \frac{2\gamma_9}{3}s_*^3
\right)
\big(-\n\otimes \n + 3\langle \n,\hx\rangle \n\otimes\hx \big)
\right.\nonumber\\
&&\hspace{50pt}\left.
- \left(
\gamma_2s_*
-\frac{\gamma_3s_*^2}{3}
+ \frac{2\gamma_9s_*^3}{3}
\right)
\big(
\hx\otimes\hx - \frac13\I
\big)
\right]\bv_*.
\end{eqnarray}

\subsection{Formula for $\cC_\gamma$ and $\cD_\gamma$}\label{CDgammaForm}
From the form of $\cC_\gamma$ and $\cD_\gamma$, it can be seen that they 
involve multiplications between the following matrices 
\[
\A,\,\,\W,\,\,\Q_*,\,\,\Q,\,\,\,\bv_*\cdot\nabla\Q,\,\, \bv\cdot\nabla\Q.
\]
\begin{itemize}
\item
For $\A, \W$, note that $\bv = \E_S\bv_*$ so that
$\partial_j\bv_i = \partial_j\big({\E_S}_{il}{\bv_*}_l\big)$ and
$\partial_i\bv_j = \partial_i\big({\E_S}_{jl}{\bv_*}_l\big)$. 
From \eqref{StokesFormGenVec}, we have
\begin{eqnarray}
{\E_S}_{il}
&\in &
\left[1, \frac{1}{r}, \frac{1}{r^3}\right]
\big\{\delta_{il}, \hx_i\hx_l\big\}\\
\partial_j{\E_S}_{il},\,\,
\partial_i{\E_S}_{jl}
&\in&
\left[\frac{1}{r^2}, \frac{1}{r^4}\right]
\Big\{
\delta_{il}\hx_j,\,\,
\delta_{ij}\hx_l,\,\,
\delta_{jl}\hx_i,\,\,
\hx_i\hx_j\hx_l
\Big\}
\end{eqnarray}
Thus,
\begin{eqnarray}
\A,\,\,\W
&\in&
\big\{
\partial_j{\E_S}_{il}{\bv_*}_l,\,\,
\partial_i{\E_S}_{jl}{\bv_*}_l
\big\}\nonumber\\
& = & 
\left[\frac{1}{r^2}, \frac{1}{r^4}\right]
\big\{
\delta_{il}\hx_j{\bv_*}_l,\,\,
\delta_{ij}\hx_l{\bv_*}_l,\,\,
\delta_{jl}\hx_i{\bv_*}_l,\,\,
\hx_i\hx_j\hx_l{\bv_*}_l
\big\}\nonumber\\
& = & 
\left[\frac{1}{r^2}, \frac{1}{r^4}\right]
\big[1,\,\,\langle \hx,\bv_*\rangle\big]
\big\{
\I,\,\,\bv_*\otimes \hx,\,\,\hx\otimes\bv_*,\,\,\hx\otimes\hx
\big\}
.
\label{CAW}
\end{eqnarray}

\item
For $\Q$, from \eqref{QABL}--\eqref{QABL3}, we have
$\Q_*,\Q\in
\big[1, \frac{1}{r}, \frac{1}{r^3}\big]
\big\{\I,\,\,\,\n\otimes \n,\,\,\,\hx\otimes \hx\big\}$. 
Hence, 
\begin{equation}
\Q_*, \Q_*^2, \Q, \Q^2 \in
\left[1, \frac{1}{r}, \ldots, \frac{1}{r^6}\right]
\big[1, \langle \hx,\n\rangle\big]
\Big\{\I,\,\,\,\n\otimes \n,\,\,\,\n\otimes \hx,\,\,\,\hx\otimes \n,\,\,\,
\hx\otimes \hx
\Big\}. 
\label{CQ}
\end{equation}
Furthermore, as
$\partial_l\Q_{ij} =
\big[\frac{1}{r^2}, \frac{1}{r^4}\big]
\{
\n_i\n_jx_l,\,\,\,\delta_{ij}x_l,\,\,\,x_ix_jx_l,\,\,\,\delta_{il}x_j,\,\,\,
\delta_{jl}x_i
\}$,
$(\bv_*\cdot\nabla\Q)_{ij} = \partial_l\Q_{ij}{\bv_*}_l$, and
$(\bv\cdot\nabla\Q)_{ij} = \partial_l\Q_{ij}\bv_l$, 
we have
\begin{eqnarray}
&&\bv_*\cdot\nabla\Q ,\,\,\,\bv\cdot\nabla\Q\nonumber\\
&\in & 
\left[\frac{1}{r^2}, \frac{1}{r^4}\right]
\big\{
\n_i\n_jx_l,\,\,\,\delta_{ij}x_l,\,\,\,x_ix_jx_l,\,\,\,\delta_{il}x_j,\,\,\,
\delta_{jl}x_i
\big\}\big\{
{\bv_*}_l, \bv_l
\big\}
\nonumber\\
&\in&
\left[\frac{1}{r^2}, \frac{1}{r^4}\right]
\big[1, \langle \hx,\bv_*\rangle \big]
\Big\{\I,\,\,\,\n\otimes \n,\,\,\,
\bv_*\otimes \hx,\,\,\,\hx\otimes \bv_*,\,\,\,
\hx\otimes \hx
\Big\}
\label{CdQ}
\end{eqnarray}

\item From the above, we have
\begin{eqnarray}
&&
\A\cdot\Q_*,\,\,\,\A\cdot\Q_*^2,\,\,\,
\A\cdot\Q,\,\,\,\A\cdot\Q^2,\,\,\,
|\Q_*|^2,\,\,\,\,|\Q|^2,\,\,\,
\Q_*\cdot(\bv_*\cdot\nabla\Q),\,\,\,
\Q\cdot(\bv\cdot\nabla\Q)
\nonumber\\
&\in&
\left[
1,\frac{1}{r^2},\ldots\frac{1}{r^6}
\right]
\big[
1,\,\,\,
\langle x,\bv_*\rangle,\,\,\,
\langle x,\n\rangle,\,\,\,
\langle \bv_*, \n\rangle
\big].
\label{CTrace}
\end{eqnarray}
\end{itemize}

Hence taking appropriate products of all the above, we have
\begin{eqnarray}
\cC_\gamma
&\in&
\left[\frac{1}{r^3},\ldots,\frac{1}{r^6}\right]
\Big[1,\,\,
\langle \hx,\n\rangle,\,\,
\langle \hx,\bv_*\rangle,\,\,
\langle \n,\bv_*\rangle
\Big]\times\nonumber\\
&&
\Big\{\I,\,\,\n\otimes \n,\,\,
\n\otimes\bv_* ,\,\,\bv_*\otimes \n,\,\,
\n\otimes \hx,\,\,\hx\otimes \n,\,\,\,
\bv_*\otimes \hx,\,\,\hx\otimes \bv_*,\,\,
\hx\otimes \hx
\Big\}.
\label{CC}
\end{eqnarray}

For $\cD_\gamma$, note that
$\Q_{ij} \in 
\left[1,\frac{1}{r},\frac{1}{r^3}\right]
\big\{\delta_{ij}, \n_i\n_j, \hx_i\hx_j\big\}$, we have
\[
\partial_k\Q_{ij} \in 
\left[\frac{1}{r^2},\frac{1}{r^4}\right]
\big\{\delta_{ij}x_k,\,\, \n_i\n_jx_k,\,\,
\delta_{ik}\hx_j,\,\, \delta_{jk}\hx_i,\,\,\hx_i\hx_j\hx_k\big\}.
\]
Hence
\begin{equation}
\cD_\gamma\in
\left[\frac{1}{r^4},\ldots,\frac{1}{r^9}\right]
\Big[1,\,\,\,
\langle \hx,\n\rangle,\,\,\,
\langle \hx,\bv_*\rangle,\,\,\,
\langle \n,\bv_*\rangle
\Big]\Big\{
\hx,\,\,\n,\,\,\bv_*
\Big\}.
\label{CD}
\end{equation}

\section{Solution of isotropic Stokes flow in bounded domain}
\label{StokesBdDom}
Even though our analysis is in the exterior domain $\R^3\backslash\bB_{a}(0),$ 
the simulation domain is assumed to be an annulus with the inner and outer radii $a$ and $R,$ respectively, with $a\ll R.$ 
As a validation of our numerical code, we verify that our simulations for the 
standard isotropic Stokes flow match with the analytical solution in this \emph{finite domain}. 
We also establish that the solution in the bounded domain retains the decay properties in an infinite domain as long as we stay away from the outer boundary. 

We first compute the analytical solution to the Stokes flow in an annular domain $\Omega_{a,R}:=\bB_R(0)\backslash \bB_a(0)$. 
In this case, we still have
\eqref{potential.flow} and the first part of \eqref{potential.form}. 
Now the boundary conditions for $u_r, u_\theta$ become:
\begin{eqnarray*}
\text{at $r=a$:} && u_r=u_\theta=0;\\
\text{at $r=R$:} && u_r=V\cos\theta,\,\,u_\theta=-V\sin\theta
\end{eqnarray*}
which are translated to:
\begin{equation}\label{f.bc}
f(a)=0,\,\,f'(a)=0,\,\,f(R)=\frac{VR^2}{2},\,\,f'(R)=VR.
\end{equation}

The Stokes equation ($-\Delta\bu + \nabla p=0$) leads to following form of $f$:
\[
f(r) = \frac{A}{r} + Br + Cr^2 + Dr^4
\]
where the coefficients are determined by the boundary conditions 
\eqref{f.bc}. 
We then have the following system of linear equations:
\begin{eqnarray*}
A+a^2B+a^3C+a^5D & = & 0,\\
-A+a^2B+2a^3C+4a^5D & = & 0,\\
2A+2R^2B+2R^3C+2R^5D & = & VR^3,\\
-A+R^2B+2R^3C+4R^5D & = & VR^3.
\end{eqnarray*}
Upon introducing $\displaystyle \lambda = \frac{a}{R}$, the solution to above
system is given by
\begin{eqnarray}
A &=& \frac{\lambda^3\left(1+\lambda+\lambda^2\right)}
{(1-\lambda)^3\left(4+7\lambda+4\lambda^2\right)}R^3V=:\eta_A(\lambda)R^3V,\\
B &=& -\frac{3\lambda\left(1+\lambda+\lambda^2+\lambda^3+\lambda^4\right)}
{(1-\lambda)^3\left(4+7\lambda+4\lambda^2\right)}RV=:\eta_B(\lambda)RV,\\
C &=& \frac{\left(4+\lambda(1+\lambda)\left(4+9\lambda^2\right)\right)}
{2(1-\lambda)^3\left(4+7\lambda+4\lambda^2\right)}V=:\eta_C(\lambda)V\\
D &=& -\frac{3\lambda\left(1+\lambda\right)}
{2(1-\lambda)^3\left(4+7\lambda+4\lambda^2\right)}\frac{V}{R^2}=:\eta_D(\lambda)\frac{V}{R^2}.
\end{eqnarray}
Then we have
\begin{eqnarray}
u_r  &=&  \frac{1}{r^2\sin\theta}\frac{\partial\Psi}{\partial\theta} = \frac{2f(r)}{r^2}\cos\theta\nonumber\\
&=& 2V\left(
\eta_A(\lambda)\left(\frac{R}{r}\right)^3 + \eta_B(\lambda)\left(\frac{R}{r}\right) + \eta_C(\lambda) + \eta_D(\lambda)\left(\frac{r}{R}\right)^2
\right)\cos\theta
\label{eq:ur}\\
u_\theta &=& \frac{-1}{r\sin\theta}\frac{\partial\Psi}{\partial r}
= -\frac{f'(r)}{r}\sin\theta\nonumber\\
&=&V\left(
\eta_A(\lambda)\left(\frac{R}{r}\right)^3 - \eta_B(\lambda)\left(\frac{R}{r}\right) - 2\eta_C(\lambda) - 4\eta_D(\lambda)\left(\frac{r}{R}\right)^2
\right)\sin\theta.
\label{eq:ut}
\end{eqnarray}

We note the self-similarity or decay structures of the solution. 
These can be used to benchmark the numerical solution. In particular, if $a\ll R$ so that $\lambda\ll1$, we have
\[
\eta_A(\lambda)=\frac{\lambda^3}{4}+O\left(\lambda^4\right),
\,\,\, \eta_B(\lambda)=-\frac{3\lambda}{4}+O\left(\lambda^2\right),
\,\,\, \eta_C(\lambda)=\frac{1}{2}+\frac{9\lambda}{8}+O\left(\lambda^2\right),
\,\,\, \eta_D(\lambda)=-\frac{3\lambda}{8}+O\left(\lambda^2\right).
\]
One can identify three distinct parameter regimes (a) $r\sim a$, (b) $a\ll r\ll R$, and (c) $r\sim R$. We will mostly interested in regime (b) as it corresponds to the flow far away from the particle, yet it is unaffected by the boundary of the computational domain. From \eqref{eq:ur} and \eqref{eq:ut}, we have
the following asymptotics:
\begin{enumerate}
\item If $r\sim a,$ then $r/R\sim\lambda,$ then
\begin{multline}
u_r  \sim V\left(1- \frac{3}{2}\left(\frac{a}{r}\right)
+\frac{1}{2}\left(\frac{a}{r}\right)^3  \right)\cos\theta
\end{multline}
and
\begin{multline}
u_\theta \sim V\left(-1+\frac{3}{4}\left(\frac{a}{r}\right)
+\frac{1}{4}\left(\frac{a}{r}\right)^3  \right)\sin\theta.
\end{multline}

\item If $a\ll r\ll R,$ then $\lambda\ll r/R\ll 1,$ then
\begin{multline}
u_r  \sim V\left(
 1- \frac{3}{2}\left(\frac{a}{r}\right) \right)\cos\theta
\end{multline}
and
\begin{multline}
u_\theta \sim V\left(
-1+\frac{3}{4}\left(\frac{a}{r}\right) \right)\sin\theta.
\end{multline}

\item If $r\sim R,$ then 
\begin{multline}
u_r  \sim V\left(
 1+\frac{3}{4}\left(\frac{a}{r}\right)\left(- 2 + 3\left(\frac{r}{R}\right)- \left(\frac{r}{R}\right)^3\right)\right)\cos{\theta}
\end{multline}
and
\begin{multline}
u_\theta  \sim V\left(
  - 1+\frac{3}{4}\left(\frac{a}{r}\right)\left(1-3\left(\frac{r}{R}\right)+ 2\left(\frac{r}{R}\right)^3\right)\right)\sin{\theta}.
\end{multline}
\end{enumerate}
Note that regimes (a) and (b) are consistent with the exact solution 
\eqref{ur.form} and \eqref{uth.form} in the exterior domain.

In the following, we compute classical Stokes' flow with
$\bv_* = Ve_1$. We will plot $\bv_1$ on the $yz$-plane demonstrating its
radially symmetric behavior. To this end, note that
$\displaystyle \bv_1 
= u_r\cos\theta - u_\theta\sin\theta$. 
On the $yz$-plane, $\sin\theta = 1$ and hence by \eqref{eq:ut} we have,
\begin{equation}
\bv_1 
= -u_\theta\Big|_{\sin\theta=1}
= -V\left(
\eta_A(\lambda)\left(\frac{R}{r}\right)^3 - \eta_B(\lambda)\left(\frac{R}{r}\right) - 2\eta_C(\lambda) - 4\eta_D(\lambda)\left(\frac{r}{R}\right)^2
\right).
\end{equation}
Our numerical results recover the above 
three asymptotics (a), (b), and (c). To illustrate this, we plot
the rescaled radial profile of $\bv_1$ for $a \leq r \leq R$,
\begin{equation}
g(r) 
:= r\left(\frac{\bv_1}{V} - 1\right)
= 
-r\left(
\eta_A(\lambda)\left(\frac{R}{r}\right)^3 - \eta_B(\lambda)\left(\frac{R}{r}\right) - 2\eta_C(\lambda) - 4\eta_D(\lambda)\left(\frac{r}{R}\right)^2
+ 1\right).
\end{equation}
The results are depicted in Figures 13 and 14.

\begin{center}
(a)\includegraphics [height=2in]
{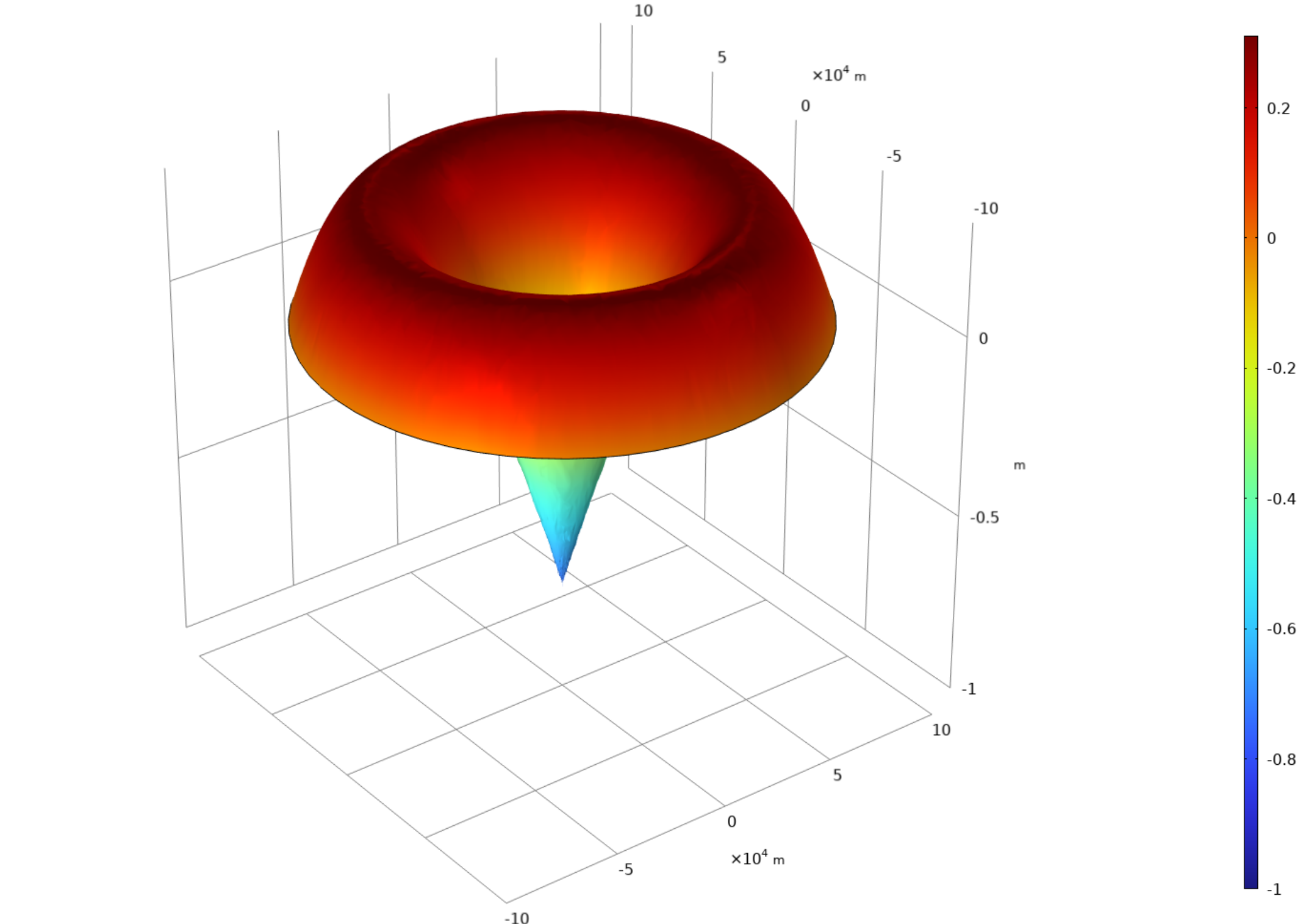}
\,\,\,
(b)\includegraphics [height=2in]
{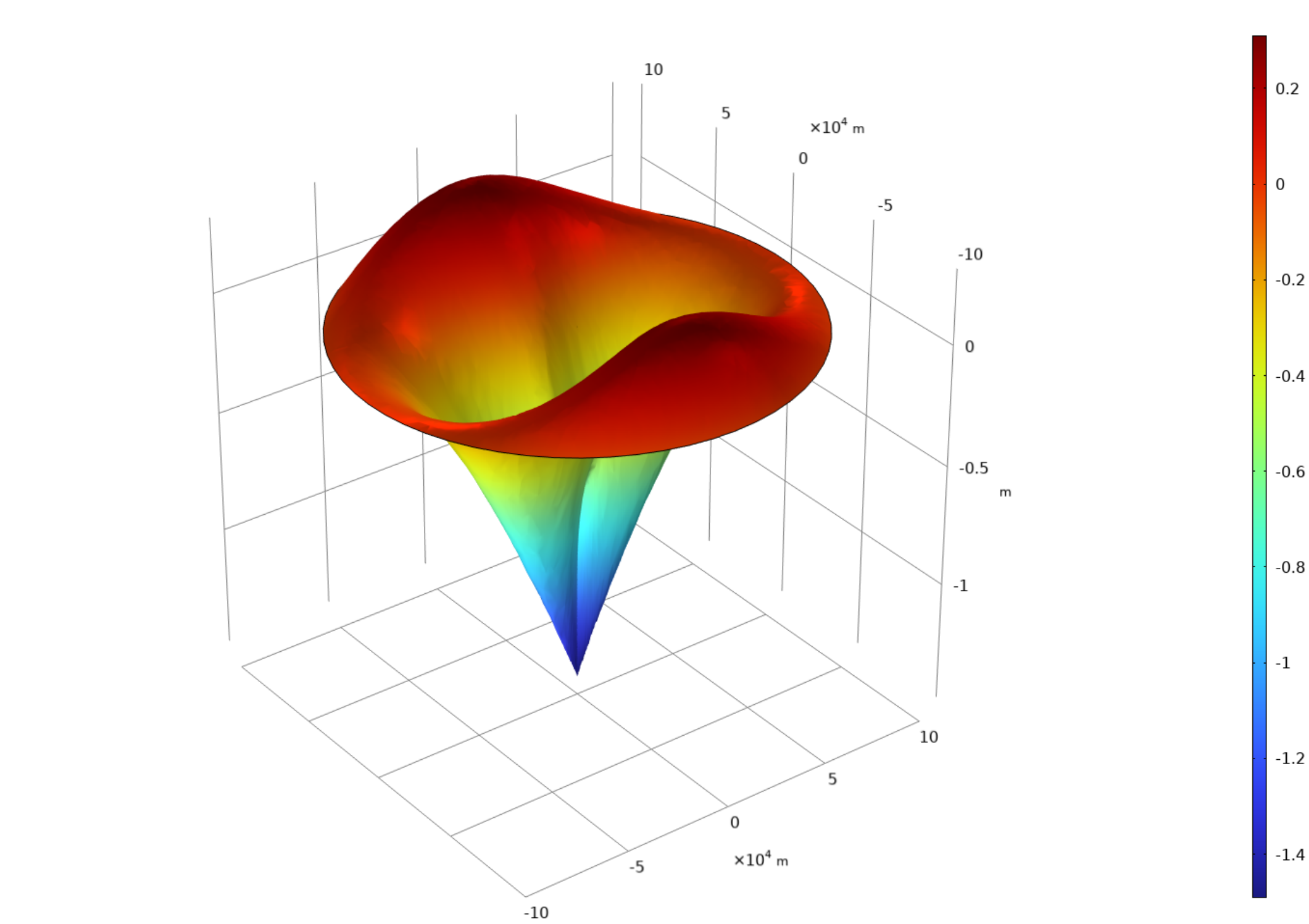}\\
Figure 13. $\displaystyle r\left(\frac{\bv_1}{V} - 1\right)$ for
classical Stokes flow:\\
(a) in $yz$-plane; 
(b) in $xy$-plane.
\end{center}
\begin{center}
(a)\includegraphics [height=2in]
{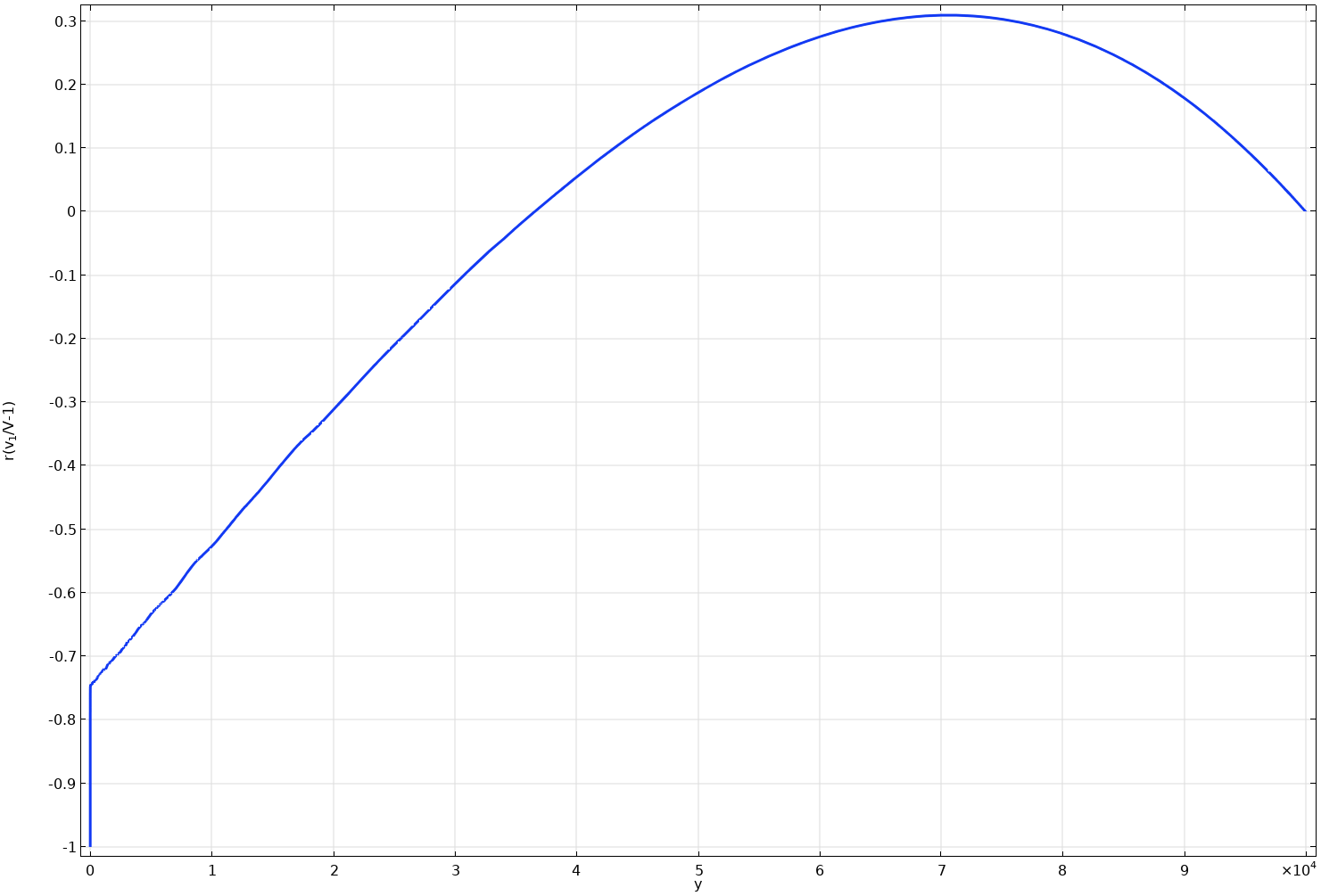}
\,\,\,
(b)\includegraphics [height=2in]
{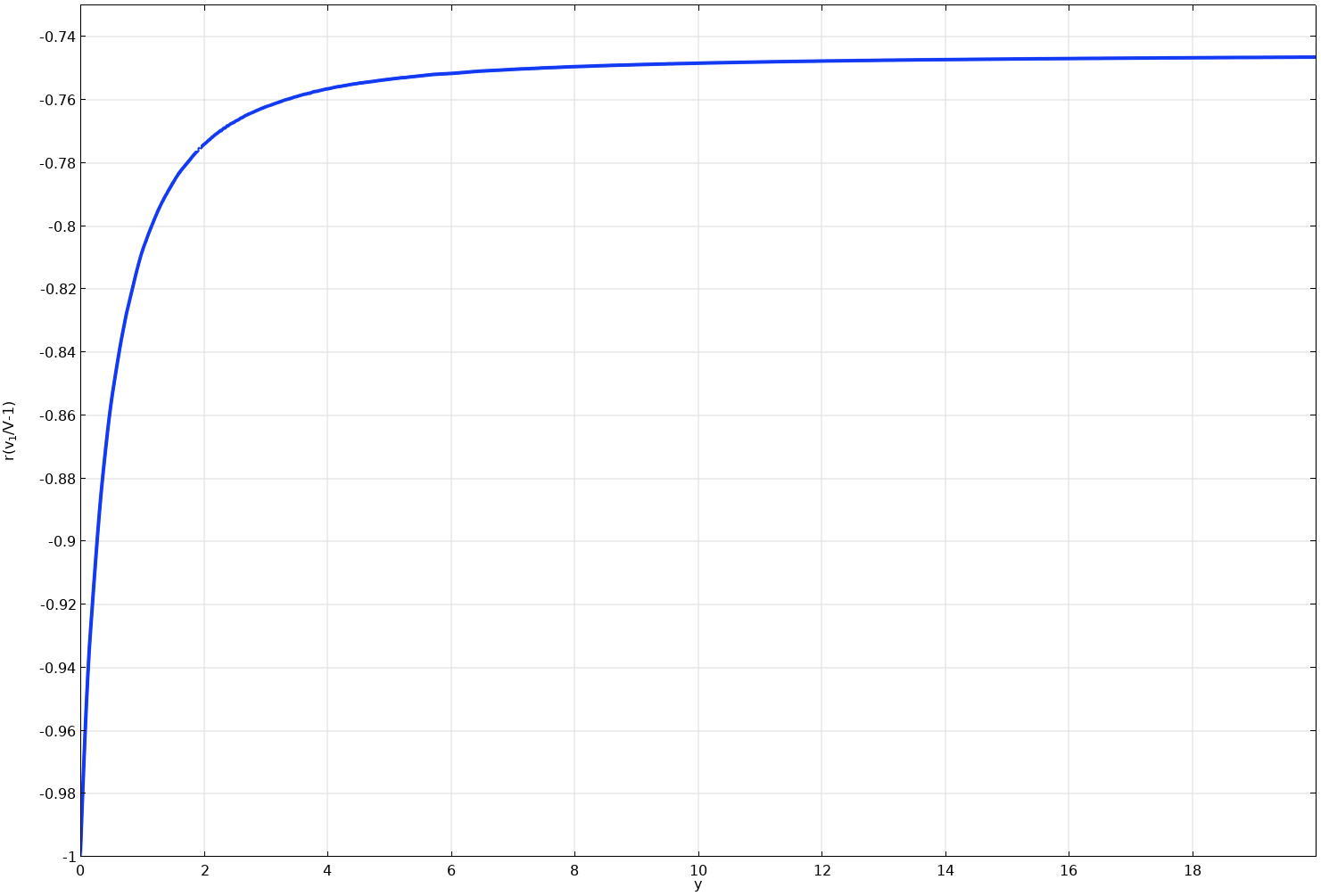}\\
Figure 14.
(Rescaled) radial profile $\bv_1$ for classical Stokes flow
(in the $yz$-plane):\\
(a) $\displaystyle g(r) = r\left(\frac{\bv_1}{V} - 1\right)$;
(b) zoomed version of (a).
\end{center}

As further demonstrations, we compare the radial behavior between 
(along the $y$-axis) the $\bv_1$ in the finite and 
infinite domain calculations. Although the overall profiles differ, due to 
the finite size effect, they do coincide very well near the core, i.e. in regimes (a) and (b) above.
\begin{center}
(a)\includegraphics [height=1.8in]
{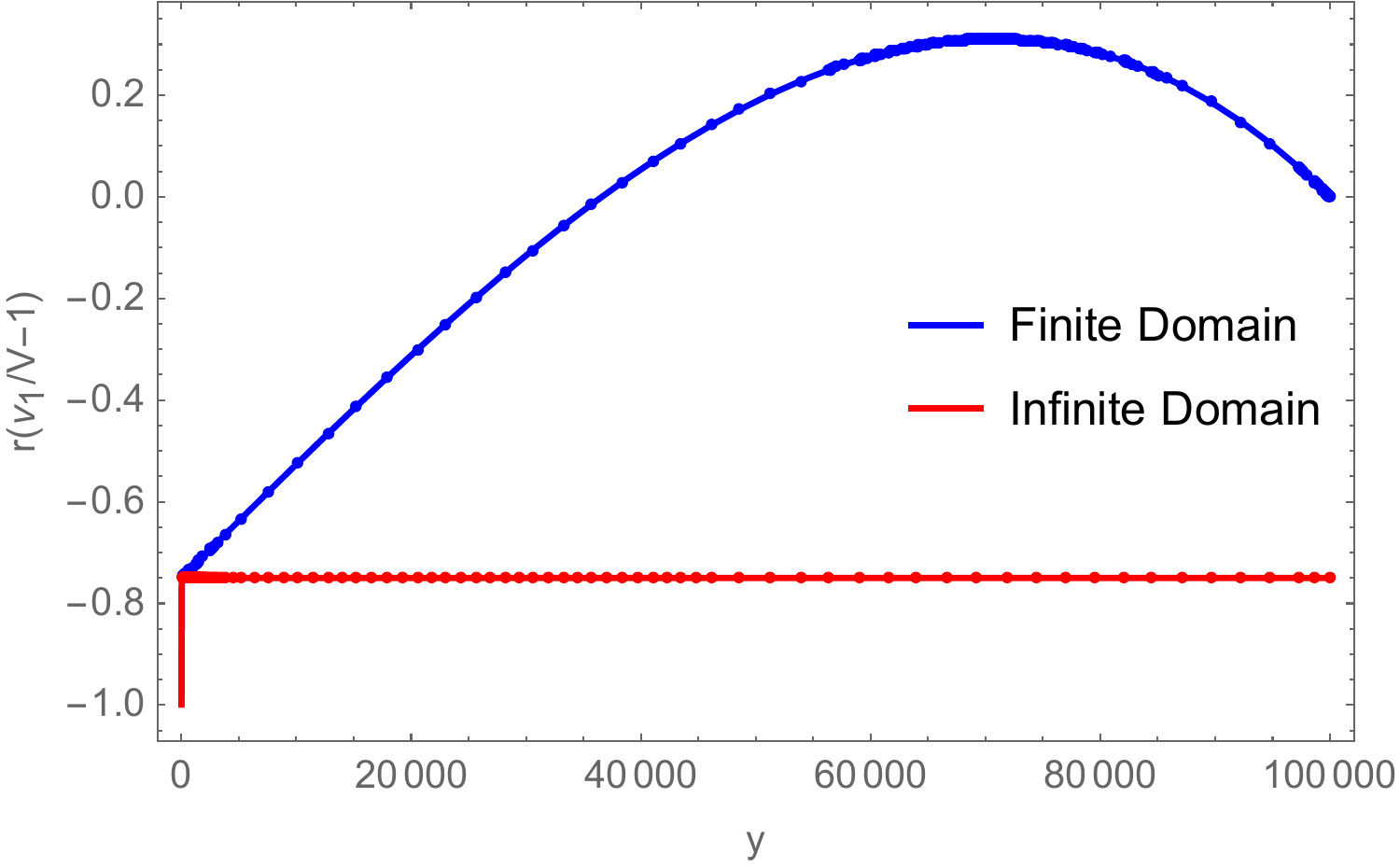}
\,\,\,
(b)\includegraphics [height=1.8in]
{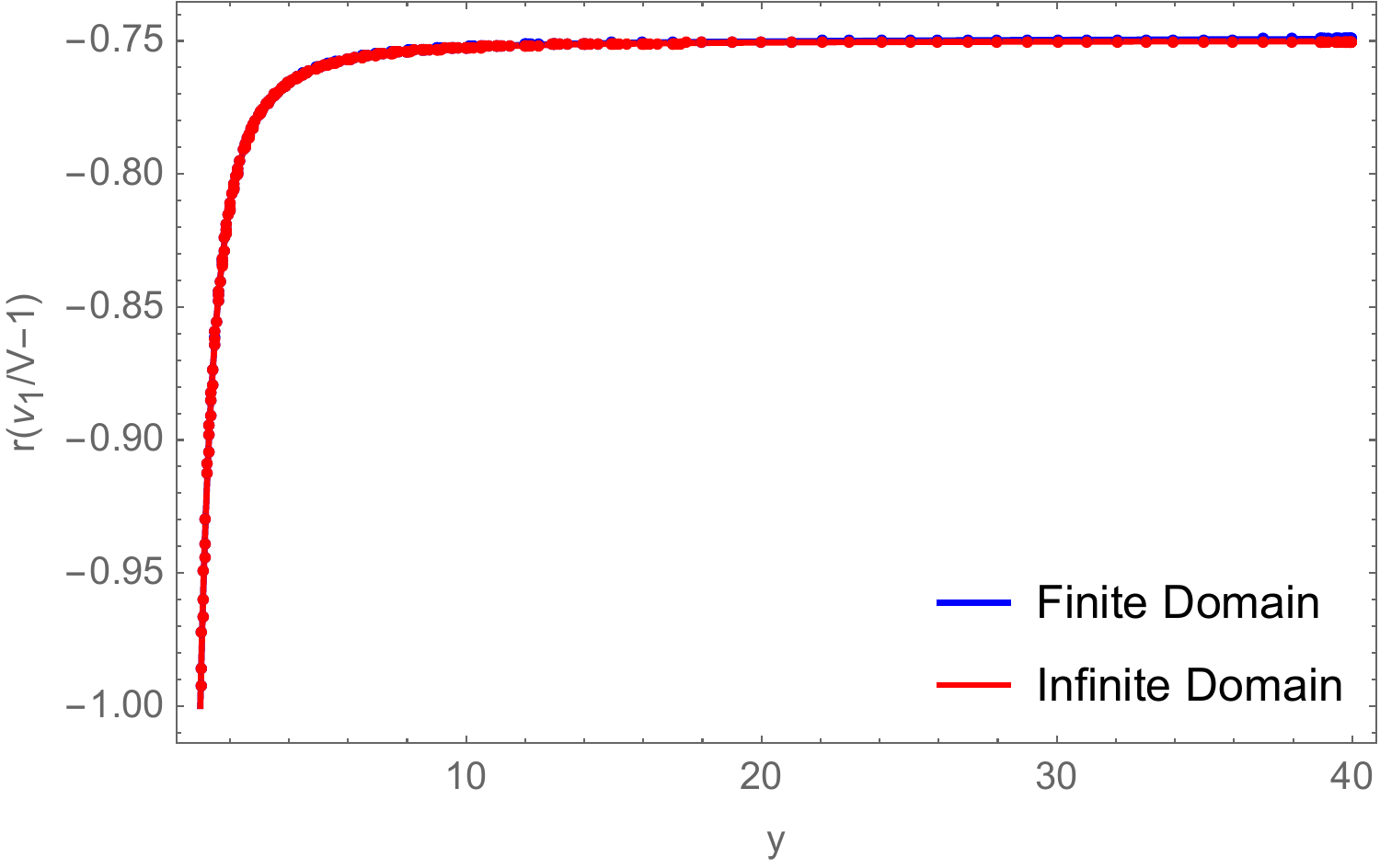}\\
Figure 15.
Comparision between the rescaled radial profiles of $\bv_1$
(along the $y$-axis) for classical Stokes flows:\\
(a) blue: finite domain; red: infinite domain; 
(b) zoomed version of (a).
\end{center}

\bibliographystyle{ieeetr}
\bibliography{stokeslet}
\end{document}